\documentclass % ~~~~~~~~~~~~~~~~~~~~~~~~~~~~~~~~~~~~~~~~~~~~~~~~~~ %
[																	% 
	10pt,															% 
    journal,														%
 	letterpaper,													%
	final,															%
	oneside,														%
	twocolumn,														%
	romanappendices													%
]																	%
{ieeetransactions}													%
\usepackage[english]{babel}											%
\usepackage[T1]{fontenc}											%
\usepackage{dsfont}													%
\usepackage{booktabs}												%
\usepackage{array}													%
\usepackage{bm}														%
\usepackage{enumerate}												%
\usepackage[table]{xcolor}											%
\usepackage{amsmath, amssymb, amsfonts}								% 
\usepackage[section,verbose]{placeins}								%
\usepackage{subfigure}												%
\usepackage{cases}													%
\usepackage{ifthen}													%
	\usepackage[pdftex]{graphicx}									% 
	\DeclareGraphicsExtensions{.pdf,.png,.jpg,.mps}					%
\usepackage{framed}													%
\usepackage[amsmath,thmmarks,framed,hyperref]{ntheorem}				%
\usepackage{eso-pic}												%
\usepackage[switch, pagewise, mathlines]{lineno}					%
\usepackage[color=red!90!black, textwidth=3cm]{todonotes}			%
\usepackage{algorithm}												%
\usepackage{algpseudocode}											% [noend]
\usepackage{acronym}												%
\let\labelindent\relax % for IEEE transactions conflicts			%
\usepackage{enumitem}												%
\usepackage{tikz}				% (note: AFTER \pdfminorversion)	%
\usepackage{pgfplots}												%
\usepgfplotslibrary{colormaps}										%
\usepgfplotslibrary{groupplots}										%
\pgfplotsset{compat=1.3}											%
\def\PgfboxplotsWidth{0.5em}										%
\def\PgfboxplotsOutliersOpacity{0.3}								%
\pgfplotsset{
    box plot/.style={
        /pgfplots/.cd,
        black,
        only marks,
        mark=-,
        mark size=\PgfboxplotsWidth,
        /pgfplots/error bars/.cd,
        y dir=plus,
        y explicit,
    },
    box plot box/.style={
        /pgfplots/error bars/draw error bar/.code 2 args={%
            \draw  ##1 -- ++(\PgfboxplotsWidth,0pt) |- ##2 -- ++(-\PgfboxplotsWidth,0pt) |- ##1 -- cycle;
        },
        /pgfplots/table/.cd,
        y index=2,
        y error expr={\thisrowno{3}-\thisrowno{2}},
        /pgfplots/box plot
    },
    box plot top whisker/.style={
        /pgfplots/error bars/draw error bar/.code 2 args={%
            \pgfkeysgetvalue{/pgfplots/error bars/error mark}%
            {\pgfplotserrorbarsmark}%
            \pgfkeysgetvalue{/pgfplots/error bars/error mark options}%
            {\pgfplotserrorbarsmarkopts}%
            \path ##1 -- ##2;
        },
        /pgfplots/table/.cd,
        y index=4,
        y error expr={\thisrowno{2}-\thisrowno{4}},
        /pgfplots/box plot
    },
    box plot bottom whisker/.style={
        /pgfplots/error bars/draw error bar/.code 2 args={%
            \pgfkeysgetvalue{/pgfplots/error bars/error mark}%
            {\pgfplotserrorbarsmark}%
            \pgfkeysgetvalue{/pgfplots/error bars/error mark options}%
            {\pgfplotserrorbarsmarkopts}%
            \path ##1 -- ##2;
        },
        /pgfplots/table/.cd,
        y index=5,
        y error expr={\thisrowno{3}-\thisrowno{5}},
        /pgfplots/box plot
    },
    box plot median/.style={
        /pgfplots/box plot
    },
	box plot outliers/.style={
        /pgfplots/.cd,
        black,
        only marks,
		opacity = \PgfboxplotsOutliersOpacity,
	}
}
\acrodef{BLUE}		[BLUE]			{Best Linear Unbiased Estimator}
\acrodef{GP}		[GP]			{Gaussian Process}
\acrodef{MAP}		[MAP]			{Maximum A Posteriori}
\acrodef{MVUE}		[MVUE]			{Minimum Variance Unbiased Estimator}
\acrodef{MMSE}		[MMSE]			{Minimum Mean Square Error}
\acrodef{MSE}		[MSE]			{Mean Square Error}
\acrodef{ML}		[ML]			{Maximum Likelihood}
\acrodef{LMMSE}		[LMMSE]			{Linear Minimum Mean Square Error}
\acrodef{LS}		[LS]			{Least Squares}
\acrodef{WSN}		[WSN]			{Wireless Sensor Network}
\acrodef{NARX}		[NARX]			{Nonlinear AutoRegressive eXogenus}
\acrodef{NCS}		[NCS]			{Networked Control System}
\acrodef{RKHS}		[RKHS]			{Reproducing Kernel Hilbert Space}
\acrodef{RN}		[RN]			{Regularization Network}
\acrodef{LTI}		[LTI]			{Linear Time Invariant}
\acrodef{DSM}		[DSM]			{Distributed Subgradient Method}
\acrodef{DCM}		[DCM]			{Distributed Control Method}
\acrodef{ADMM}		[ADMM]			{Alternating Direction Method of Multipliers}
\acrodef{FNRC}		[FNRC]			{Fast Newton-Raphson Consensus}
\acrodef{NR}		[NR]			{Newton-Raphson}
\acrodef{NRC}		[NRC]			{Newton-Raphson Consensus}
\acrodef{JC}		[JC]			{Jacobi Consensus}
\acrodef{GDC}		[GDC]			{Gradient Descent Consensus}
\acrodef{rv}		[r.v.]			{random variable}
\acrodef{rvv}		[r.\textbf{v}.]	{random vector}
\newcounter{generalCounter}
\theoremclass		{Theorem}
\theoremstyle		{plain}
\theoremheaderfont	{\normalfont\bfseries}
\theorembodyfont	{\itshape}
\newtheorem	{theorem}		[generalCounter]	{Theorem}

\newtheorem	{assumption}	[generalCounter]	{Assumption}

\newtheorem	{lemma}			[generalCounter]	{Lemma}

\theoremstyle		{nonumberplain}
\theorempostwork	{\hrule}
\theoremseparator	{}
\theoremheaderfont	{\normalfont\bfseries}
\theorembodyfont	{\normalfont}
\theoremsymbol		{\ensuremath{\diamondsuit}}

\newshadedtheorem	{shadedproof}		{Proof}
\newtheorem			{proof}				{Proof}
\pgfdeclarelayer{ultrabackground}
\pgfdeclarelayer{background}
\pgfdeclarelayer{foreground}
\pgfsetlayers{ultrabackground,background,main,foreground}
\def\TablesColumnsColor{black!4}
\newcolumntype
{g}
{
	>{\centering \columncolor{\TablesColumnsColor} \arraybackslash}
	p{0.15\textwidth}
	<{}
}
\newcolumntype
{w}
{
	>{\centering \arraybackslash}
	p{0.15\textwidth}
	<{}
}
\newcommand{\cmax}				[1]	{\left[ #1 \right]_{c}}
\newcommand{\DefinedAs}			[0]	{\mathrel{\mathop:}=}
\newcommand{\IDefinedAs}		[0]	{=\mathrel{\mathop:}}

\newcommand{\BigOOf}			[1]	{O \left( #1 \right)}
\newcommand{\MSEOf}				[1]	{\textrm{MSE} \left( #1 \right)}

\newcommand{\OnesVector}			[0]	{\mathds{1}}

\newcommand{\UniformDistribution}					[2]	{\mathcal{U} \left[ #1, #2 \right]}

\newcommand{\Reals}									[0]	{\mathbb{R}}
\newcommand{\PositiveReals}							[0]	{\mathbb{R}_{+}}

\newcommand	{\Assumption}			[0]	{Assumption}
\newcommand	{\Assumptions}			[0]	{Assumptions}
\newcommand	{\Section}				[0]	{Section}
\newcommand	{\Sections}				[0]	{Sections}
\newcommand	{\Equation}				[0]	{Equation}

\newcommand	{\Figure}				[0]	{Figure}

\newcommand	{\Lemma}				[0]	{Lemma}

\newcommand	{\Table}				[0]	{Table}

\newcommand	{\Theorem}				[0]	{Theorem}
\newcommand	{\Theorems}				[0]	{Theorems}
\newcommand	{\Algorithm}			[0]	{Algorithm}

\newcommand	{\GrayText}				[1]	{{\color{gray} \itshape #1}}
\newcommand	{\GrayMath}				[1]	{{\color{gray} #1}}
\makeatletter                                                       %
\let\IEEEproof\proof                                                %
\let\proof\@undefined                                               %
\let\endproof\@undefined                                            %
\makeatother                                                        %
\begin{document}                                                    %
%                                                                   %
% \overrideIEEEmargins                                              %
%                                                                   %
% \setpagewiselinenumbers                                             %
% \modulolinenumbers[1]                                               %
% \linenumbers                                                        %
% \newcommand{\rr}[1]{{\color{blue}{#1}}}								%
%                                                                   %
% ~~~~~~~~~~~~~~~~~~~~~~~~~~~~~~~~~~~~~~~~~~~~~~~~~~~~~~~~~~~~~~~~~ %
%                                                                   %
\title{Newton-Raphson Consensus \\ for Distributed Convex Optimization}
%                                                                   %
\author                                                             %
{                                                                   %
    Damiano Varagnolo, Filippo Zanella, Angelo Cenedese, \\ Gianluigi Pillonetto, Luca Schenato
    \thanks                                                         %
    {                                                               %
		D.\ Varagnolo is with the Department of Computer Science, Electrical and Space Engineering, Lule{\aa} University of Technology, Lule{\aa} Sweden. Email: {\tt damiano.varagnolo@ltu.se}. F.\ Zanella, A.\ Cenedese, G.\ Pillonetto and L.\ Schenato are with the Department of Information Engineering, Universit\`a di Padova, Padova, Italy. Emails: {\tt \{fzanella | angelo.cenedese | giapi | schenato \}@dei.unipd.it}. 
    }                                                               %
    \thanks                                                         %
    {                                                               %
		This work is supported by the Framework Programme for Research and Innovation Horizon 2020 under the grant agreement n.\ 636834 ``DISIRE'', the Swedish research council Norrbottens Forskningsr{\aa}d, by the University of Padova under the ``Progetto di Ateneo CPDA147754/14-New statistical learning approach for multi-agents adaptive estimation and coverage control.'', and by the Italian Ministry of Education under the grant agreement SCN 00398 ``Smart \& safe Energy-aware Assisted Living''.
		This paper is an extended and revised version of~\cite{zanella_et_al__2011__newton_raphson_consensus_for_distributed_convex_optimization,zanella_et_al__2012__multidimensional_newton_raphson_consensus_for_distributed_convex_optimization}.
    }                                                               %
}                                                                   %
%                                                                   %
\date{}                                                             %
\maketitle                                                          %
\IEEEoverridecommandlockouts                                        %
\IEEEpeerreviewmaketitle											%
%                                                                   %
% ~~~~~~~~~~~~~~~~~~~~~~~~~~~~~~~~~~~~~~~~~~~~~~~~~~~~~~~~~~~~~~~~~ %

\begin{abstract}
	We address the problem of distributed unconstrained convex optimization under separability assumptions, i.e., the framework where each agent of a network is endowed with a local private multidimensional convex cost, is subject to communication constraints, and wants to collaborate to compute the minimizer of the sum of the local costs. We propose a design methodology that combines average consensus algorithms and separation of time-scales ideas. This strategy is proved, under suitable hypotheses, to be globally convergent to the true minimizer. Intuitively, the procedure lets the agents distributedly compute and sequentially update an approximated Newton-Raphson direction by means of suitable average consensus ratios. We show with numerical simulations that the speed of convergence of this strategy is comparable with alternative optimization strategies such as the \acl{ADMM}. Finally, we propose some alternative strategies which trade-off communication and computational requirements with convergence speed.
\end{abstract}

\begin{IEEEkeywords}
    Distributed optimization, unconstrained convex optimization, consensus, multi-agent systems, Newton-Raphson methods, smooth functions.
\end{IEEEkeywords}

\acresetall

\section{Introduction}
\label{sec:introduction}

Optimization is a pervasive concept underlying many aspects of modern life~\cite{shor__1985__minimization_methods_for_non_differentiable_functions,bertsekas_et_al__2003__convex_analysis_and_optimization,boyd_vandenberghe__2004__convex_optimization}, and it also includes the management of distributed systems, i.e., artifacts composed by a multitude of interacting entities often referred to as ``agents''. Examples are transportation systems, where the agents are both the vehicles and the traffic management devices (traffic lights), and smart electrical grids, where the agents are the energy producers-consumers and the power transformers-transporters.

Here we consider the problem of distributed optimization, i.e., the class of algorithms suitable for networked systems and characterized by the absence of a centralized coordination unit~\cite{tsitsiklis__1984__problems_in_decentralized_decision_making_and_computation,bertsekas_tsitsiklis__1997__parallel_and_distributed_computation,bertsekas__1998__network_optimization__continuous_and_discrete_models}. Distributed optimization tools have received an increasing attention over the last years, concurrently with the research on networked control systems. Motivations comprise the fact that the former methods let the networks self-organize and adapt to surrounding and changing environments, and that they are necessary to manage extremely complex systems in an autonomous way with only limited human intervention. In particular we focus on unconstrained convex optimization, although there is a rich literature also on distributed constrained optimization such as Linear Programming~\cite{burger_et_al__2012__a_distributed_simplex_algorithm_for_degenerate_linear_programs_and_multi_agent_assignments}.

\subsection*{Literature review}

The literature on distributed unconstrained convex optimization is extremely vast and a first taxonomy can be based whether the strategy uses or not the Lagrangian framework, see, e.g.,~\cite[Chap.~5]{boyd_vandenberghe__2004__convex_optimization}.

Among the distributed methods exploiting Lagrangian formalism, the most widely known algorithm is \ac{ADMM}~\cite{bertsekas__1982__constrained_optimization_and_lagrange_multiplier_methods}, whose roots can be traced back to~\cite{hestenes__1969__multiplier_and_gradient_methods}. Its efficacy in several practical scenarios is undoubted, see, e.g.,~\cite{boyd_et_al__2010__distributed_optimization_and_statistical_learning_via_the_admm} and references therein. A notable size of the dedicated literature focuses on the analysis of its convergence performance and on the tuning of its parameters for optimal convergence speed, see, e.g., \cite{erseghe_et_al__2011__fast_consensus_by_the_admm} for \ac{LS} estimation scenarios,  \cite{he_yuan__2011__on_the_o_1_on_t_convergence_rate_of_alternating_direction_method} for linearly constrained convex programs, and \cite{deng_yin__2012__on_the_global_linear_convergence_of_the_generalized_ADMM} for more general ADMM algorithms. Even if proved to be an effective algorithm, \ac{ADMM} suffers from requiring synchronous communication protocols, although some recent attempts for asynchronous and distributed implementations have appeared~\cite{wei_ozdaglar__2012__distributed_alternating_direction_method_of_multipliers,mota_et_al__2012__distributed_admm_for_model_predictive_control_and_congestion_control,jakovetic_et_al__2011__cooperative_convex_optimization_in_networked_systems__augmented_lagrangian_algorithms_with_directed_gossip_communication}.

On the other hand, among the distributed methods not exploiting Lagrangian formalisms, the most popular ones are the \acp{DSM}~\cite{demyanov_vasilev__1985__nondifferentiable_optimization}. Here the optimization of non-smooth cost functions is performed by means of subgradient based descent/ascent directions. These methods arise in both primal and dual formulations, since sometimes it is better to perform dual optimization. Subgradient methods have been exploited for several practical purposes, e.g., to optimally allocate resources in \acp{WSN}~\cite{johansson__2008__on_distributed_optimization_in_networked_systems}, to maximize the convergence speeds of gossip algorithms~\cite{boyd_et_al__2006__randomized_gossip_algorithms}, to manage optimality criteria defined in terms of ergodic limits~\cite{ribeiro__2010__ergodic_stochastic_optimization_algorithms_for_wireless_communication_and_networking}. Several works focus on the analysis of the convergence properties of the \ac{DSM} basic algorithm~\cite{nedic_bertsekas__2001__incremental_subgradient_methods_for_nondifferentiable_optimization,nedic_et_al__2001__distributed_asynchronous_incremental_subgradient_methods,nedic_ozdaglar__2008__approximate_primal_solutions_and_rate_analysis_for_dual_subgradient_methods} (see~\cite{kiewel__2004__convergence_of_approximate_and_incremental_subgradient_methods_for_convex_optimization} for a unified view of many convergence results). We can also find analyses for several extensions of the original idea, e.g., directions that are computed combining information from other agents~\cite{blatt_et_al__2007__a_convergent_incremental_gradient_method_with_a_constant_step_size,xiao_boyd__2006__optimal_scaling_of_a_gradient_method_for_distributed_resource_allocation} and stochastic errors in the evaluation of the subgradients \cite{ram_et_al__2009__incremental_stochastic_subgradient_algorithms_for_convex_optimization}. Explicit characterizations can also show trade-offs between desired accuracy and number of iterations~\cite{nedic_ozdaglar__2009__distributed_subgradient_methods_for_multi_agent_optimization}. 

These methods have the advantage of being easily distributed, to have limited computational requirements and to be inherently asynchronous as shown in~\cite{johansson_et_al__2009__a_randomized_incremental_subgradient_method_for_distributed_optimization_in_networked_systems,lobel_et_al__2011__distributed_multi_agent_optimization_with_state_dependent_communication,nedic__2010__asynchronous_broadcast_based_convex_optimization_over_a_network}. However they suffer from low convergence rate since they require the update steps to decrease to zero as $1/t$ (being $t$ the time) therefore as a consequence the rate of convergence is sub-exponential. In fact, one of the current trends is to design strategies that improve the convergence rate of \acp{DSM}. For example, a way is to accelerate the convergence of subgradient methods by means of multi-step approaches, exploiting the history of the past iterations to compute the future ones~\cite{ghadimi_et_al__2012__accelerated_gradient_methods_for_networked_optimization}. Another is to use Newton-like methods, when additional smoothness assumptions can be used. These techniques are based on estimating the Newton direction starting from the Laplacian of the communication graph. More specifically, distributed Newton techniques have been proposed in dual ascent scenarios~\cite{jadbabaie_et_al__2009__a_distributed_newton_method_for_network_optimization,zargham_et_al__2011__accelerated_dual_descent_for_network_optimization,wei_et_al__2011__a_distributed_newton_method_for_network_utility_maximization}. Since the Laplacian cannot be computed exactly, the convergence rates of these schemes rely on the analysis of inexact Newton methods~\cite{dembo_et_al__1982__inexact_newton_methods}. These Newton methods are shown to have super-linear convergence under specific assumptions, but can be applied only to specific optimization problems such as network flow problems.

Recently, several alternative approaches to \ac{ADMM} and \ac{DSM} have appeared. For example, in \cite{nedic_et_al__2010__constrained_consensus_and_optimization_in_multi_agent_networks,zhu_martinez__2012__on_distributed_convex_optimization_under_inequality_and_equality_constraints} the authors construct contraction mappings by means of cyclic projections of the estimate of the optimum onto the constraints. A similar idea based on contraction maps is used in F-Lipschitz methods~\cite{fischione__2011__f_lipschitz_optimization_with_wsn_applications} but it requires additional assumptions on the cost functions. Other methods are the control-based approach~\cite{wang_elia__2010__control_approach_to_distributed_optimization} which exploits distributed consensus, the distributed randomized Kaczmarz method~\cite{freris_zouzias__2012__fast_distributed_smoothing_for_network_clock_synchronization} for quadratic cost functions,  and distributed dual sub-gradient methods \cite{necoara_nedelcu__2014__distributed_dual_gradient_methods_and_error_bound_conditions}.

\subsection*{Statement of contributions}

Here we propose a distributed Newton-Raphson optimization procedure, named \ac{NRC}, for the exact minimization of smooth multidimensional convex separable problems, where the global function is a sum of private local costs. With respect to the classification proposed before, the strategy exploits neither Lagrangian formalisms nor Laplacian estimation steps. More specifically, it is based on average consensus techniques~\cite{garin_schenato__2011__a_survey_on_distributed_estimation_and_control_applications_using_linear_consensus_algorithms} and on the principle of separation of time-scales~\cite[Chap.~11]{khalil__2001__nonlinear_systems}. The main idea is that agents compute and keep updated, by means of average consensus protocols, an approximated Newton-Raphson direction that is built from suitable Taylor expansions of the local costs. Simultaneously, agents move their local guesses towards the Newton-Raphson direction. It is proved that, if the costs satisfy some smoothness assumptions and the rate of change of the local update steps is sufficiently slow to allow the consensus algorithm to converge, then the \ac{NRC} algorithm exponentially converges to the global minimizer.

The main contribution of this work is to propose an algorithm that extends Newton-Raphson ideas in a distributed setting, thus being able to exploit second order information to speed up converge rate. By using singular perturbation theory we formally show that under suitable assumptions the convergence of the algorithm is exponential (linear in logspace). Differently, \ac{DSM} algorithms have sublinear convergence rate even if the cost functions are smooth~\cite{nedic_et_al__2010__constrained_consensus_and_optimization_in_multi_agent_networks,nedic_olshevsky__2013__distributed_optimization_over_time_varying_directed_graphs}, although they are easy to implement and can be employed also for non-smooth cost functions and for constrained optimization. We also show by means of numerical simulations on real-world database benchmarks that the proposed algorithm exhibits faster convergence rates (in number of communications) than standard implementations of distributed \ac{ADMM} algorithms \cite{boyd_et_al__2010__distributed_optimization_and_statistical_learning_via_the_admm}, probably due to the second-order information embedded into the Newton-Raphson consensus. Although we have no theoretical guarantee of the superiority of the proposed algorithmic in terms of convergence rate, these simulations suggest that it is at least a potentially competitive algorithm. Moreover, one of the promising features of the \ac{NRC} is that it is essentially based on average consensus algorithms, for which there exist robust implementations that encompass asynchronous communications, time-varying network topologies~\cite{fagnani_zampieri__2008__randomized_consensus_algorithms_over_large_scale_networks}, directed graphs~\cite{dominguez_garcia_et_al__2011__distributed_algorithms_for_consensus_and_coordination_in_the_presence_of_packet-dropping_communication_links}, and packet-losses effects.

\subsection*{Structure of the paper}

The paper is organized as follows: \Section~\ref{sec:notation} collects the notation used through the whole paper, while \Section~\ref{sec:problem_formulation} formulates the considered problem and provides some ancillary results that are then used to study the convergence properties of the main algorithm. \Section~\ref{sec:newton_raphson_consensus} proposes the main optimization algorithm, provides convergence results and describes some strategies to trade-off communication and computational complexities with convergence speed. \Section~\ref{sec:numerical_examples} compares, via numerical simulations, the performance of the proposed algorithm with several distributed optimization strategies available in the literature. Finally, \Section~\ref{sec:conclusions} collects some final observations and suggests future research directions. We collect all the proofs in the Appendix.

\section{Notation}
\label{sec:notation}

We model the communication network as a graph $\mathcal{G}=(\mathcal{N},\mathcal{E})$ whose vertices $\mathcal{N} \DefinedAs \{ 1, 2, \ldots, N \}$ represent the agents and whose edges $(i,j) \in \mathcal{E}$ represent the available communication links. We assume that the graph is undirected and connected, and that the matrix $P\in\Reals^{N\times N}$ is stochastic, i.e., its elements are non-negative, it is s.t.\ $P \OnesVector = \OnesVector$ (where $\OnesVector \DefinedAs [1 \; 1 \; \cdots \; 1]^T \in \Reals^N$), symmetric, i.e., $P=P^T$ and consistent with the graph $\mathcal{G}$, in the sense that each entry $p_{ij}$ of $P$ is $p_{ij} > 0$ only if $(i,j) \in \mathcal{E}$. We recall that if $P$ is stochastic, symmetric, and includes all edges (i.e., $p_{ij} > 0$ if and only if $(i,j) \in \mathcal{E}$) then $\lim_{k \rightarrow \infty} P^k = \frac{1}{N} \OnesVector \OnesVector^T$. Such $P$'s are also often referred to as \emph{average consensus matrices}. We will indicate with $\rho(P) \DefinedAs \max_{i,\lambda_i\neq 1} \left| \lambda_i(P) \right|$ the spectral radius of $P$, with $\sigma(P) \DefinedAs 1 - \rho(P)$ its spectral gap.

We use fraction bars to indicate also Hadamard divisions, e.g., if $\bm{a} = [a_1, \ldots, a_N]^T$ and $\bm{b} = [b_1, \ldots, b_N]^T$ then
$
\displaystyle
\frac{\bm{a}}{\bm{b}}
\DefinedAs 
\left[
	\frac{a_1}{b_1}\,
	\ldots\,
	\frac{a_N}{b_N}
\right]^{T}
$. 
Fraction bars like the previous ones may also indicate pre-multiplication with inverse matrices, i.e., if $b_{i}$ is a matrix then $\displaystyle \frac{a_{i}}{b_{i}}$ indicates $b_{i}^{-1} a_{i}$. We indicate with $n$ the dimensionality of the domains of the cost functions, $k$ a discrete time index, $t$ a continuous time index. For notational simplicity we denote differentiation with $\nabla$ operators, so that $\nabla f = \partial f / \partial x$ and $\nabla^{2} f = \partial^{2} f / \partial x^{2}$. With a little abuse of notation, we will define $\chi = (x,Z)$, where $x\in \Reals^n$ and $ Z\in \Reals^{\ell\times q}$ as the vector obtained by stacking in a column both the vector $x$ and the vectorized matrix $Z$. We indicate with $\| \cdot \|$ Frobenius norms. With an other abuse of notation we also define the norm of the pair $\chi = (x, Z)$ where $x$ is a vector and $Z$ a matrix with $\|\chi\|^2=\|x\|^2+\|Z\|^2$.

When using plain italic fonts with a subscript (usually $i$, e.g., $x_{i} \in \Reals^{n}$) we refer to the local decision variable of the specific agent $i$. When using bold italic fonts, e.g., $\bm{x}$, we instead refer to the collection of the decision variables of all the various agents, e.g., $\bm{x} \DefinedAs \left[ x_{1}^T \;, \ldots, \; x_{N}^T \right]^T \in \Reals^{nN}$. To indicate special variables we will instead consider the following notation: 
\begin{equation*}
\begin{array}{rcll}
	\overline{x}
	& \DefinedAs &
	\displaystyle
	\frac{1}{N} \sum_{i=1}^N x_i
	& \quad \GrayMath{\Reals^{n}} \\
	\bm{x}^{\parallel}
	& \DefinedAs &
	\OnesVector_{N} \otimes \overline{x}
	& \quad \GrayMath{\Reals^{nN}} \\
	\bm{x}^{\perp}
	& \DefinedAs &
	\bm{x} - \bm{x}^{\parallel}
	& \quad \GrayMath{\Reals^{nN}} \\
\end{array}
\end{equation*}
As in~\cite[p.~116]{khalil__2001__nonlinear_systems}, we say that a function $V$ is a \emph{Lyapunov function} for a specific dynamics if $V$ is continuously differentiable and satisfies $V(0) = 0$, $V(x) > 0$ for $x \neq 0$, and $\dot{V}(x) \leq 0$.

\section{Problem formulation and preliminary results}
\label{sec:problem_formulation}

\subsection{Structure of the section}
\label{ssec:structure_of_the_section}

Our main contribution is to characterize the convergence properties of the distributed \ac{NR} scheme proposed in \Section~\ref{sec:newton_raphson_consensus}. In doing so we both exploit standard singular perturbation analysis tools~\cite[Chap.~11]{khalil__2001__nonlinear_systems}~\cite{kokotovic_et_al__1999__singular_perturbation_methods_in_control__analysis_and_design} and a set of ancillary results, collected for readability in this section.

The logical flow of these ancillary results is the following: \Section~\ref{ssec:stability_of_discretized_dynamics} claims that, under suitable assumptions, forward-Euler discretizations of stable continuous dynamics lead to stable discrete dynamics.  This basic result enables reasoning on continuous-time systems. Then, \Sections~\ref{ssec:stability_of_single_agent_nr_dynamics} and~\ref{ssec:stability_of_multi_agent_nr_dynamics} respectively claim that single- and multi-agent continuous-time \ac{NR} dynamics satisfy these discretization assumptions. \Sections~\ref{ssec:multi_agent_nr_dynamics_under_vanishing_perturbations} and~\ref{ssec:multi_agent_nr_dynamics_under_non_vanishing_perturbations} then generalize these dynamics by introducing perturbation terms that mimic the behavior of the proposed main optimization algorithm, and characterize their stability properties. Summarizing, the ancillary results characterize the stability properties of systems that are progressive approximations of the dynamics under investigation.

\subsection{Problem formulation}
\label{ssec:problem_formulation}

We assume that the $N$ agents of the network are endowed with cost functions $f_{i} : \Reals^n \mapsto \Reals$ so that
\begin{equation}
    \overline{f} : \Reals^n \mapsto \Reals ,
	\qquad \qquad
	\overline{f} \left( x \right)
    \DefinedAs
    \frac{1}{N}
    \sum_{i = 1}^{N}
	f_{i} \left( x \right)
\label{equ:definition_of_global_cost_function}
\end{equation}
is a well-defined global cost. We assume that the aim of the agents is to cooperate and distributedly compute the minimizer of $\overline{f}$, namely
\begin{equation}
	x^{\ast}
    \DefinedAs
	\arg \min_{x \in \Reals^{n}}
	\overline{f} \left( x \right) .
\label{equ:definition_of_global_optimum}
\end{equation}
We now enforce the following simplifying assumptions, valid throughout the rest of the paper:
\begin{assumption}[Convexity]
	The local costs $f_i$ in~\eqref{equ:definition_of_global_cost_function} are of class $\mathcal{C}^3$. Moreover the global cost $\overline{f}$ has bounded positive definite Hessian, i.e., $0 < cI \leq \nabla^2 \overline{f}(x) \leq mI$ for some $c, m \in \PositiveReals$ and $\forall x \in \Reals^n$. Moreover, w.l.o.g., we assume $\overline{f}(x^*)=0$, $c\leq 1$ and $m\geq 1$. 
\label{ass:smoothness_of_global_function}
\end{assumption}

The scalar $c$ is assumed to be known by all the agents a-priori. \Assumption~\ref{ass:smoothness_of_global_function} ensures that $x^{\ast}$ in~\eqref{equ:definition_of_global_optimum} exists and is unique. The strictly positive definite Hessian is moreover a mild sufficient condition to guarantee that the minimum $x^*$ defined in~\eqref{equ:definition_of_global_optimum} will be globally exponentially stable under the continuous and discrete Newton-Raphson dynamics described in the following \Theorem~\ref{thm:continuous_NR}. We also notice that, for the subsequent \Theorems~\ref{thm:continuous_to_discrete} and~\ref{thm:continuous_NR}, in principle just the average function $\overline{f}$ needs to have specific properties, and thus no conditions for the single $f_i$'s are required (that for example might be even non convex). For the convergence of the distributed \ac{NR} scheme we will nonetheless enforce the more restrictive \Assumptions~\ref{ass:global_lipschitzianity} and~\ref{ass:phi_p_xi_is_locally_uniformly_lipschitz}, not presented now for readability issues.
 In the rest of this section, in order to simplify notation, we will considerer, without loss of generality, the following translated cost functions:
\begin{equation}
	f'_i(x)=f_i(x+x^*), \ \ \overline{f'}(x)=\frac{1}{N}\sum_{i=1}^N f'_i(x)
\end{equation}
so that the origin becomes the minimizer of the averaged cost function $\overline{f'}(x)$, i.e. $\overline{f'}(0)=0$.

\subsection{Stability of discretized dynamics}
\label{ssec:stability_of_discretized_dynamics}

This subsection aims to show that, under suitable assumptions, forward-Euler discretization of suitable exponentially stable continuous-time dynamics maintains the same global exponential stability properties.

\begin{theorem}
	Let the continuous-time system
	\begin{equation}
		\dot{x} = \phi(x)
	\label{equ:continuous_system}
	\end{equation}
	admit $x = 0 \in \Reals^{n}$ as an equilibrium, and let $V(x) : \Reals^n \mapsto \Reals$ be a Lyapunov function for~\eqref{equ:continuous_system} for which there exist positive scalars $a_{1}, a_{2}, a_{3}, a_{4}$ s.t., $\forall x \in \Reals^{n}$,
	\begin{subnumcases}
		{\label{equ:inequalities_continuous_to_discrete}}
		a_{1} I \leq \nabla^2 V(x) \leq a_{2} I
		\label{equ:inequalities_continuous_to_discrete:nabla2_V} \\
		\frac{\partial V(x)}{\partial x} \phi(x) \leq - a_{3} \|x\|^2
		\label{equ:inequalities_continuous_to_discrete:partial_V} \\
		\|\phi(x)\| \leq a_{4} \|x\| .
		\label{equ:inequalities_continuous_to_discrete:phi_Lipschitz}
	\end{subnumcases}
	Then:
	\begin{enumerate}[label=\emph{\alph*)}]
		\item for system~\eqref{equ:continuous_system} the origin is globally exponentially stable;
		\label{item:continuous_to_discrete:continuous}
		\item for the following forward-Euler discretization of system~\eqref{equ:continuous_system},
		\begin{equation}
			x(k+1) = x(k) + \varepsilon \phi \big( x(k) \big),
		\label{equ:discrete_system}
		\end{equation}
		\label{item:continuous_to_discrete:discrete}
		there exists a positive scalar $\overline\varepsilon$ such that for every $\varepsilon \in (0, \overline\varepsilon)$ the origin is globally exponentially stable.
	\end{enumerate}
\label{thm:continuous_to_discrete}
\end{theorem}

\subsection{Stability of single-agent \ac{NR} dynamics}
\label{ssec:stability_of_single_agent_nr_dynamics}

This subsection shows that the results of \Section~\ref{ssec:stability_of_discretized_dynamics} apply to continuous NR dynamics, i.e., that forward-Euler discretizations maintain global exponential stability properties\footnote{We notice that other asymptotic properties of continuous time \ac{NR} methods are available in the literature, e.g.,~\cite{tanabe__1985__global_analysis_of_continuous_analogues_of_the_levenberg_marquardt_and_newton_raphson_methods,hauser_nedic__2005__the_continuous_newton_raphson_method_can_look_ahead}.}.

\begin{theorem}
	Let
	\begin{equation}
		\phi_{\textrm{NR}}(x) \DefinedAs - \overline{\overline{h'}}(x)^{-1} \nabla \overline{f'}(x)
		\label{equ:function_for_NR}
	\end{equation}
	be defined by a generic function $\overline{\overline{h'}}(x) \in \Reals^{n \times n}$ that satisfies the positive definiteness conditions  $cI \leq \overline{\overline{h'}}(x) = \overline{\overline{h'}}(x)^T \leq m I$ for all $x \in \Reals^n$ where $c$ and $m$ are defined in Assumption~\ref{ass:smoothness_of_global_function}. Let~\eqref{equ:function_for_NR} define both the dynamics
	\begin{equation}
		\dot{x} = \phi_{\textrm{NR}}(x) ,
	\label{equ:continuous_NR_system}
	\end{equation}
	\begin{equation}
		x(k+1) = x(k) + \varepsilon \phi_{\textrm{NR}} \big( x(k) \big) .
	\label{equ:discrete_NR_system}
	\end{equation}
	Then, under \Assumption~\ref{ass:smoothness_of_global_function}:
	\begin{enumerate}[label=\emph{\alph*)}]
		\item
		\begin{equation}
			V_{\textrm{NR}}(x)
			\DefinedAs
			\overline{f'}(x)
			\label{equ:definition_of_lyapunov_for_NR}
		\end{equation}
		is a Lyapunov function for~\eqref{equ:continuous_NR_system};
		\label{item:continuous_NR:lyapunov}
		\item there exist positive scalars $b_1, b_2, b_3, b_4$ s.t., $\forall x \in \Reals^{n}$,
		\begin{subnumcases}
			{\label{equ:inequalities_for_NR}}
			b_1 I \leq \nabla^2 V_{\textrm{NR}}(x) \leq b_2 I
			\label{eqn:NR_2} \\
			\frac{\partial V_{\textrm{NR}}}{\partial x} \phi_{\textrm{NR}}(x)
			\leq
			-b_3 \|x\|^2
			\label{eqn:NR_3} \\
			\|\phi_{\textrm{NR}}(x)\|
			\leq
			b_4 \|x\| ,
			\label{eqn:NR_4}
		\end{subnumcases}
		\label{item:continuous_NR:scalars}
	\end{enumerate}
	i.e., \Theorem~\ref{thm:continuous_to_discrete} applies to dynamics~\eqref{equ:continuous_NR_system} and~\eqref{equ:discrete_NR_system}.
\label{thm:continuous_NR}
\end{theorem}

For suitable choices of $\overline{\overline{h'}}(x)$ the dynamics~\eqref{equ:continuous_NR_system} corresponds to continuous versions of well known descent dynamics. Indeed, the correspondences are
\begin{subnumcases}
	{\overline{\overline{h'}}(x) =}
	\nabla^2 \overline{f'}(x) & \hspace{-0.7cm} $\rightarrow  \text{Newton-Raphson descent} \label{h1}$ \\
	\mathrm{diag} \left[ \nabla^2 \overline{f'}(x) \right]& \hspace{-0.7cm} $\rightarrow  \text{Jacobi descent} \label{h2}$ \\
	I & \hspace{-0.7cm} $\rightarrow  \text{Gradient descent} \label{h3}$
\end{subnumcases}
\noindent where $\mathrm{diag}[A]$ is a diagonal matrix containing the main diagonal of $A$.  Note that for every choice of $\overline{\overline{h'}}(x)$ as in~\eqref{h1}-\eqref{h3}, \Assumption~\ref{ass:smoothness_of_global_function} ensures the hypotheses\footnote{For the Jacobi descent, clearly $\min_{\|x\|=1}x^T\mathrm{diag} \left[ \nabla^2 \overline{f'}(x) \right] x= \min_{x\in\{e_1,\ldots,e_n\}}x^T\mathrm{diag} \left[ \nabla^2 \overline{f'}(x) \right] x =  \min_{x\in\{e_1,\ldots,e_n\}}x^T \nabla^{2} \overline{f'}(x) x \geq \min_{\|x\|=1}x^T\nabla^{2} \overline{f'}(x) x=c $, where $e_i$ is the $n$-dimensional vector with all zeros except for a one in the $i$-th entry.} of \Theorem~\ref{thm:continuous_NR}, therefore by combining \Theorem~\ref{thm:continuous_NR} with \Theorem~\ref{thm:continuous_to_discrete} we are guaranteed that both continuous and discrete  generalized \ac{NR} dynamics induced by~\eqref{equ:function_for_NR} are globally exponentially stable:
\begin{lemma}
	Under \Assumption~\ref{ass:smoothness_of_global_function}, the origin is a globally exponentially stable point for dynamics~\eqref{equ:continuous_NR_system}. Moreover there exists $\overline{\varepsilon}>0$ such that the origin is a globally exponentially stable point also for dynamics~\eqref{equ:discrete_NR_system} for all $\varepsilon < \overline{\varepsilon}$.
\label{lem:NR_system_is_globally_exponentially_Stable}
\end{lemma}

The previous lemma and theorems do not require $\overline{\overline{h'}}(x)$ to be differentiable. However, differentiability may be used to linearize the system dynamics and obtain explicit rates of convergence. In fact, the linearized dynamics around the origin is given by 
$$
	F(0)
	\DefinedAs
	\frac
	{\partial \phi_{\textrm{NR}}(0)}
	{\partial x}
	=
	- \overline{\overline{h'}}(0)^{-1} \nabla^{2}\overline{f'}(0)
	-
	\frac
	{\partial \overline{\overline{h'}}(0)^{-1}}
	{\partial x}
	\nabla \overline{f'}(0) .
$$
In particular, for the NR descent it holds that $\overline{\overline{h'}}(x) = \nabla^2\overline{f'} (x)$. Thus in this case $F(0) = -I$, since $\nabla \overline{f'}(0) = 0$, and this says that the linearized continuous time NR dynamics is $\dot{x}=-x$, independent of the cost $\overline{f'}(x)$ and whose rate of convergence is unitary and uniform along any direction.

\subsection{Stability of multi-agent \ac{NR} dynamics}
\label{ssec:stability_of_multi_agent_nr_dynamics}

We now generalize~\eqref{equ:continuous_NR_system} by considering $N$ coupled dynamical systems that, when starting at the very same initial condition, behave like $N$ decoupled systems~\eqref{equ:continuous_NR_system}. This novel dynamics is the core of the slow-dynamics embedded in the main algorithm presented in \Section~\ref{sec:newton_raphson_consensus}. In this section we also include additional assumptions to show that the generalization of~\eqref{equ:continuous_NR_system} presented here preserves global exponential stability and some other additional properties.

To this aim we introduce some additional notation: let $h'_i(x) : \Reals^{n} \mapsto \Reals^{n \times n}, i = 1, \ldots, N$ be defined according to one of the possible three cases
\begin{subnumcases}{
	h'_{i}(x)
	=}
	\nabla^2 f'_{i}(x) \\
	\mathrm{diag} \left[ \nabla^2 f'_{i}(x) \right] \\
	I 
	\label{equ:choice_of_h_i}
\end{subnumcases}
so that $h'_{i}(x) = h'_{i}(x)^{T}$ for all $x$. Moreover let
\begin{equation*}
\begin{array}{rcll}
	h' \big( \bm{x} \big)
    & \DefinedAs &
    \left[
		h'_1 \big( x_{1} \big) \;,
        \ldots, \;
		h'_N \big( x_{N} \big)
    \right]^T
	& \GrayMath{\Reals^{nN} \mapsto \Reals^{nN \times n}} \\
	\overline{h'} \big( \bm{x} \big)
	& \DefinedAs &
	\displaystyle
	\frac{1}{N} \sum_{i=1}^N h'_i(x_i)
	& \GrayMath{\Reals^{nN} \mapsto \Reals^{n \times n}} \\
	\overline{\overline{h'}} \big( \overline{x} \big)
	& \DefinedAs &
	\displaystyle
	\frac{1}{N} \sum_{i=1}^N h'_{i} \big( \overline{x} \big)
	& \GrayMath{\Reals^{n} \mapsto \Reals^{n \times n}} \\
\end{array}
\end{equation*}
be additional composite functions defined starting from the $h'_i$'s (recall that $\bm{x} \DefinedAs \left[ x_{1}^T \;, \ldots, \; x_{N}^T \right]^T \in \Reals^{nN}$ and that $\overline{x} \DefinedAs \frac{1}{N} \sum_{i=1}^N x_i \in \Reals^{n}$). Let moreover
\begin{equation}
	g'_i(x) \DefinedAs h'_i(x)x - \nabla f'_i(x) \qquad \GrayMath{\Reals^{n} \mapsto \Reals^{n}}
\label{equ:definition_of_gi}
\end{equation}
and $g'(\bm{x}), \overline{g'}(\bm{x}), \overline{\overline{g'}}(\overline{x})$ be defined accordingly as for $h'_i$.

The definitions of $h'_{i}$ and $g'_{i}$ are instrumental to generalize the NR dynamics~\eqref{equ:continuous_NR_system} to the distributed case. Indeed, let 
\begin{equation}
	\psi \big( \bm{x} \big)
	\DefinedAs
	\displaystyle
	\overline{h'}(\bm{x})^{-1}\ \overline{g'}(\bm{x}) 
	\qquad
	\GrayMath{\Reals^{nN} \mapsto \Reals^n}
\label{equ:definition_of_psi}
\end{equation}
(with the existence of $\overline{h'}(\bm{x})^{-1}$ guaranteed by the following \Assumption~\ref{ass:global_lipschitzianity}). It is easy to verify that the previous functions satisfy the following properties:
\begin{subnumcases}
	{\label{equ:important_equivalences}}
	\overline{h'} \big( \bm{x}^{\parallel} \big)
	=
	\overline{\overline{h'}} \big( \overline{x} \big)
	\label{equ:important_equivalences:h} \\
	\overline{g'} \big( \bm{x}^{\parallel} \big)
	=
	\overline{\overline{g'}} \big( \overline{x} \big)
	=
	\overline{\overline{h'}} \big( \overline{x} \big)
	\overline{x}
	-
	\nabla \overline{f'} \big( \overline{x} \big)
	\label{equ:important_equivalences:g} \\
	\psi \big( \bm{x}^{\parallel} \big)
	=
	\overline{x}
	-
	\overline{\overline{h'}} \big( \overline{x} \big)^{-1}
	\nabla \overline{f'} \big( \overline{x} \big) 
	\label{equ:important_equivalences:psi}
\end{subnumcases}
Consider then
\begin{equation}
	\bm{\dot x}= \phi_{\textrm{PNR}}(\bm{x})
	\DefinedAs
	- \bm{x} + \OnesVector_{N} \otimes \psi(\bm{x}),
\label{equ:function_for_coupled_NR}
\end{equation}
that can be also equivalently written as
$$
	\dot{x}_i = -x_i + \psi(\bm{x}), \qquad i = 1, \ldots, N,
$$
i.e., as the combination of $N$ independent dynamical systems that are driven by the same forcing term $\psi(\bm{x})$. 

As mentioned above, this dynamics embeds the centralized generalized \ac{NR} dynamics since, under identical initial conditions $x_i(0) = \overline{x}(0) \in \Reals^{n}$ for all $i$, the trajectories coincide, i.e., $x_i(t) = \overline{x}(t), \forall i, \forall t\geq 0$. Moreover, due to~\eqref{equ:important_equivalences:psi},
\begin{equation}
\begin{array}{rcl}
	\dot{\overline{x}}
	& = &
	- \overline{x} + \psi(\OnesVector_{N} \otimes \overline{x}) \\
	& = &
	- \overline{x} + \overline{x} - \overline{\overline{h'}} \big( \overline{x} \big)^{-1} \nabla \overline{f'}(\overline{x})
	=
	\phi_{\textrm{NR}}(\overline{x}) ,
\end{array}
\label{equ:NR_evolution_for_average_x}
\end{equation}
i.e., we obtain dynamics~\eqref{equ:function_for_NR}, that is, thanks to \Theorem~\ref{thm:continuous_NR} and the assumption that $\overline{h'}(\bm{x})$ is invertible, globally exponentially stable.

The question is then whether dynamics~\eqref{equ:function_for_coupled_NR} is exponentially stable also in the general case where the $x_{i}(0)$'s may not be identical. To characterize this case we assume some additional global properties:
\begin{assumption}[Global properties]
	The local costs $f'_{1}, \ldots, f'_{N}$ in~\eqref{equ:definition_of_global_cost_function} are s.t.\ there exist positive scalars $m_g, a_g, a_h, a_\psi$ s.t., $\forall \, x, x' \in \Reals^{n}$ and $\forall \, \bm{x}, \bm{x}' \in \Reals^{nN}$,
	\begin{subnumcases}
		{\label{equ:inequalities_for_global_lipschitzianity}}
		{c I \leq \overline{h'}(\bm{x}) \leq m} I 
		\label{equ:inequalities_for_global_lipschitzianity:h} \\
		\left\| \overline{g'}(\bm{x}) \right\| \leq m_g
		\label{equ:inequalities_for_global_lipschitzianity:g} \\
		\left\| g'_i(x) - g'_i(x') \right\| \leq a_g \left\| x - x' \right\| \\
		\left\| h'_i(x) - h'_i(x') \right\| \leq a_{h} \left\| x - x' \right\| \\
		\left\| \psi(\bm{x}) - \psi(\bm{x}') \right\| \leq a_{\psi} \left\| \bm{x} - \bm{x}' \right\|
	\end{subnumcases}
	with $c$ and $m$ from \Assumption~\ref{ass:smoothness_of_global_function}.
\label{ass:global_lipschitzianity}
\end{assumption}

Note that \Assumption~\ref{ass:global_lipschitzianity} implies
\begin{subnumcases}
	{\label{equ:implied_inequalities_for_global_lipschitzianity}}
	\left\| \overline{g'}(\bm{x}) - \overline{g'}(\bm{x}') \right\|
	\leq
	a_{g} \left\| \bm{x}-\bm{x}'\right\| \\
	\left\| \overline{h'}(\bm{x}) - \overline{h'}(\bm{x}') \right\|
	\leq
	a_{h} \left\| \bm{x}-\bm{x}'\right\| \\
	\left\| g'(\bm{x}) -g'(\bm{x}') \right\|
	\leq
	a_{g} \left\| \bm{x}-\bm{x}'\right\| \\
	\left\| h'(\bm{x}) -h'(\bm{x}') \right\|
	\leq
	a_{h} \left\| \bm{x}-\bm{x}'\right\| 
\end{subnumcases}

Using the previous assumptions we can now prove global stability of dynamics~\eqref{equ:function_for_coupled_NR}:

\begin{theorem}
	Under \Assumptions~\ref{ass:smoothness_of_global_function} and~\ref{ass:global_lipschitzianity}, and for a suitable positive scalar $\eta$,
	\begin{enumerate}[label=\emph{\alph*)}]
		\item
		\begin{equation}
			V_{\textrm{PNR}}(\bm{x})
			\DefinedAs
			V_{\textrm{NR}}(\overline{x})
			+
			\frac{1}{2} \eta \| \bm{x}^\perp \|^2
			=
			\overline{f'}(\overline{x})
			+
			\frac{1}{2} \eta \| \bm{x}^\perp \|^2
		\label{equ:definition_of_V_PNR}
		\end{equation}
		is a Lyapunov function for~\eqref{equ:function_for_coupled_NR};
		\label{item:continuous_multidimensional_NR:lyapunov}
		\item there exist positive scalars $b_5, b_6, b_7, b_8$ s.t., $\forall \bm{x} \in \Reals^{nN}$,
		\begin{subnumcases}
			{\label{equ:inequalities_for_continuous_multidimensional_NR}}
			b_5 I\leq\nabla^2 V_{\textrm{PNR}}(\bm{x})
			\leq
			b_6 I
			\label{eqn:PNR_2} \\
			\frac{\partial V_{\textrm{PNR}}}{\partial\bm{x}}\phi_{\textrm{PNR}}(\bm{x})
			\leq
			-b_7 \| \bm{x} \|^2
			\label{eqn:PNR_3} \\
			\|\phi_{\textrm{PNR}}(\bm{x})\|
			\leq
			b_8 \|\bm{x}\| .
			\label{eqn:PNR_4}
		\end{subnumcases}
		\label{item:continuous_multidimensional_NR:scalars}
	\end{enumerate}
	\label{thm:continuous_multidimensional_NR}
\end{theorem}

As in \Lemma~\ref{lem:NR_system_is_globally_exponentially_Stable}, combining \Theorem~\ref{thm:continuous_multidimensional_NR} with \Theorem~\ref{thm:continuous_to_discrete} it is possible to claim that~\eqref{equ:function_for_coupled_NR} and its discrete-time counterpart are globally exponentially stable.

\subsection{Multi-agent \ac{NR} dynamics under vanishing perturbations}
\label{ssec:multi_agent_nr_dynamics_under_vanishing_perturbations}

We now aim to generalize the dynamics $\phi_{\textrm{PNR}}(\bm{x})$ by considering some perturbation term, that will be described by the variable $\bm{\chi}$. Let then $\bm{\chi}^{y} \DefinedAs \left( \chi^y_1,\ldots,\chi^y_N \right)$ where $\chi_i^y \in \Reals^{n}$, $\bm{\chi}^{z} \DefinedAs \left( \chi^z_1,\ldots,\chi^z_N \right)$ where $\chi_i^z = \left( \chi_i^z \right)^{T} \in \Reals^{n\times n}$, and  $\bm{\chi} \DefinedAs \left( \bm{\chi}^{y}, \bm{\chi}^{z} \right)$. We also define the operator
$
	\left[ \; \cdot \; \right]_{c}
	\; : \;
	\Reals^{nN \times n} \mapsto \Reals^{nN \times n}
$, 
which indicates the component-wise matrix-operation
\begin{equation}
	\left[ \bm{z} \right]_{c}
	=
	\begin{bmatrix}
		z_{1} \\
		\vdots \\
		z_{N} \\
	\end{bmatrix}_{c}
	\DefinedAs	
	\begin{bmatrix}
		z'_{1} \\
		\vdots \\
		z'_{N} \\
	\end{bmatrix}
	\qquad
	z'_{i}
	=
	\left\{
		\begin{array}{ll}
			z_{i}	& \displaystyle \text{if } z_{i} \ge \frac{c}{2} I \vspace{0.2cm} \\
			\displaystyle \frac{c}{2} I		& \text{otherwise.}
		\end{array}
	\right.%}
\label{equ:definition_of_c_operator}
\end{equation}

Consider then the perturbed version of the multi-agent NR dynamics~\eqref{equ:function_for_coupled_NR},
\begin{equation}
	\dot{\bm{x}}
	=
	\phi_x(\bm{x}, \bm{\chi})
	\DefinedAs
	- \bm{x}
	- \OnesVector_{N} \otimes x^{\ast}
	+ \frac
	{ \bm{\chi}^{y} + \OnesVector_{N} \otimes \left( \overline{g'}(\bm{x}) + \overline{h'}(\bm{x})x^{\ast} \right)}
	{ \left[ \bm{\chi}^{z} + \OnesVector_{N} \otimes \overline{h'}(\bm{x})  \right]_c}
	,
\label{equ:dynamics_of_phi_x}
\end{equation}
where the division is a Hadamard division, as recalled in \Section~\ref{sec:notation}.  Direct inspection of dynamics~\eqref{equ:dynamics_of_phi_x} then shows that
\begin{equation}
	\phi_x(\bm{x}, \bm 0)=\phi_{\textrm{PNR}}(\bm{x})	\label{equ:definition_of_psi_x_of_i} .
\end{equation}
The next lemma provides perturbations interconnection bounds that will be used in \Theorem~\ref{thm:distributed_NR_converges_to_the_global_optimum}.
\begin{lemma}
	Under \Assumptions~\ref{ass:smoothness_of_global_function} and~\ref{ass:global_lipschitzianity} there exist positive scalars $a_{x}$, $a_{\Delta}$ s.t., for all $\bm{x}$ and $\bm{\chi}$,
	\begin{subnumcases}
		{\label{equ:inequalities_for_phi_p_x}}
		\|\phi_x(\bm{x}, \bm{\chi})\|
		\leq
		a_x \big( \| \bm{x} \| + \| \bm{\chi} \| \big)
		\label{eqn:phi_p_x_1} \\
		\| \phi_x( \bm{x}, \bm{\chi} ) - \phi_{\textrm{PNR}}(\bm{x})\|
		\leq
		a_\Delta \| \bm{\chi} \| . 
		\label{eqn:phi_p_x_2}
	\end{subnumcases}
\label{lem:phi_x_is_lipschitz}
\end{lemma}

\subsection{Multi-agent \ac{NR} dynamics under non-vanishing perturbations}
\label{ssec:multi_agent_nr_dynamics_under_non_vanishing_perturbations}

Let us now consider some additional properties of the flow~\eqref{equ:dynamics_of_phi_x} for some specific non-vanishing perturbation.  Consider then the perturbations $\xi^{y} \in \Reals^n$ and $\xi^{z} \in \Reals^{n \times n}$, and their multi-agents versions $\bm{\xi}^y = \OnesVector_N \otimes \xi^{y}$, $\bm{\xi}^y = \OnesVector_N \otimes \xi^{z}$. Consider also the shorthand $\bm{\xi} = (\bm{\xi}^y, \bm{\xi}^z)$. The equilibrium points of the dynamics induced by $\phi_x(\bm{x}, \bm{\xi})$ are characterized by the following theorem:

\begin{theorem}
	Let $\xi^{y} \in \Reals^n$, $\xi^{z} \in \Reals^{n \times n}$, $\xi = (\xi^{y}, \xi^{z})$, $\bm{\xi}^y = \OnesVector_N \otimes \xi^{y}$, $\bm{\xi}^z = \OnesVector_N \otimes \xi^{z}$, $\bm{\xi} = (\bm{\xi}^y, \bm{\xi}^z)$, and consider the equation
	$$ \phi_x(\bm{x}, \bm{\xi}) = 0 ,$$
	defining the equilibrium points of the dynamics $\dot{\bm{x}} = \phi_x(\bm{x}, \bm{\xi})$. Then, under \Assumptions~\ref{ass:smoothness_of_global_function} and~\ref{ass:global_lipschitzianity} there exist a positive scalar $r > 0$ and a unique continuously differentiable function \mbox{$\bm{x}^{eq} : \mathcal{B}_r\rightarrow \mathbb{R}^{nN}$} where $\mathcal{B}_r \DefinedAs \{ \xi \; | \; \|\xi\|\leq r\}$ such that
	\begin{equation}
		\phi_x\big(\bm{x}^{eq}(\xi), \bm{\xi} \big) = 0,
		\quad
		\bm{x}^{eq}(0) = 0;
		\label{equ:phi_x_has_equilibrium_x_eq}
	\end{equation}
	Moreover, $\bm{x}^{eq}(\xi) = \OnesVector_{N} \otimes x^{eq}(\xi)$, with
	\begin{equation}
		x^{eq}(\xi)
		=
		\Big(
			\overline{\overline{h'}} \big( x^{eq}(\xi) \big)
			+
			\xi^{z}
		\Big)^{-1}
		\Big(
			\overline{\overline{g'}} \big( x^{eq}(\xi) \big)
			+
			\xi^{y} -\xi^zx^*
		\Big).
	\label{equ:solution_of_x_of_xi_scalar}
	\end{equation}
\label{the:zeros_of_phi_xi_is_a_smooth_manifold}
\end{theorem}

\Theorem~\ref{the:zeros_of_phi_xi_is_a_smooth_manifold} allows to define
\begin{equation}
	\phi'_x(\bm{x}, \xi)
	\DefinedAs
	\phi_x \big(\bm{x} + \OnesVector_N \otimes x^{eq}(\xi), \bm{\xi} \big)
\label{equ:phi_p_xi}
\end{equation}
and the corresponding dynamics
\begin{equation}
	\bm{\dot x}
	=
	\phi'_x (\bm{x}, \xi)
\label{equ:dynamics_of_phi_p_xi}
\end{equation}
which corresponds to the translated version of the original perturbed system $\phi_x (\bm{x}, \bm{\xi})$, which has now the property that the origin is an equilibrium point, i.e., $\phi'_x (\bm 0, \xi)=0, \forall \|\xi\|\leq r$.

To prove the global exponential stability of~\eqref{equ:dynamics_of_phi_p_xi} we need the flow $\phi'_x$ to satisfy a global Lipschitz condition:
\begin{assumption}[Global Lipschitz perturbation]
	There exist positive scalars $a_{\xi}$ and $r$ such that, for all $\bm{x} \in \Reals^{nN}$ and $\xi$ satisfying $\| \xi \| \leq r$,
	$$
		\left\| \phi'_x(\bm{x}, \xi) - \phi'_x(\bm{x}, 0) \right\|
		\leq
		a_\xi \| \xi \| \| \bm{x}\|.
	$$
\label{ass:phi_p_xi_is_locally_uniformly_lipschitz}
\end{assumption}

With these assumptions we can prove that the origin is a globally exponentially stable equilibrium for dynamics~\eqref{equ:dynamics_of_phi_p_xi}:
\begin{theorem}
	Under \Assumptions~\ref{ass:smoothness_of_global_function},~\ref{ass:global_lipschitzianity} and~\ref{ass:phi_p_xi_is_locally_uniformly_lipschitz},
	\begin{enumerate}[label=\emph{\alph*)}]
		\item $V_{\textrm{PNR}}(\bm{x})$ defined in~\eqref{equ:definition_of_V_PNR} is a Lyapunov function for~\eqref{equ:dynamics_of_phi_p_xi};
		\label{item:perturbed_NR_is_globally_exponentially_stable:V_PNR_is_Lyapunov}
	\item there exist positive scalars $r$, $b_{7}'$, $b_{8}'$ s.t., for all $\bm{x} \in \Reals^{nN}$ and $\xi$ satisfying $\|\xi\|\leq r$,
		\begin{subnumcases}
			{\label{equ:inequalities_for_perturbed_NR}}
			\frac{\partial V_{\textrm{PNR}}}{\partial\bm{x}}\phi'_x(\bm{x},\xi)
			\leq
			-b'_7\| \bm{x}\|^2
			\label{eqn:xi_3} \\
			\|\phi'_x(\bm{x},\xi)\|
			\leq
			b'_8\|\bm{x}\|.
			\label{eqn:xi_4}
		\end{subnumcases}
		\label{item:inequalities_on_xi}
	\end{enumerate}
\label{thm:perturbed_NR_is_globally_exponentially_stable}
\end{theorem}

Again, as in \Lemma~\ref{lem:NR_system_is_globally_exponentially_Stable}, combining \Theorem~\ref{thm:perturbed_NR_is_globally_exponentially_stable} with \Theorem~\ref{thm:continuous_to_discrete} it is possible to claim that~\eqref{equ:dynamics_of_phi_p_xi} and its discrete-time counterpart are globally exponentially stable.

\subsection{Quadratic Functions}

Before presenting the main algorithm, we show that quadratic costs satisfy all the previous assumptions. In fact, let us consider then
$$ f_i(x) = \frac{1}{2}(x-d_i)^TA_i(x-d_i)+e_i, \ \ A_i=A_i^T  $$
%
%and assume that, w.l.o.g.,
%%
%$$ \sum_{i}A_id_i=0, \ \ \sum_{i} (d_i^TA_id_i+e_i)=0$$
%%
%condition that can be met by a simple change of coordinate and a constant translation of the cost functions.
Based on this definition we have the following result:
\begin{theorem}
	Quadratic costs that satisfy 
	$$A \DefinedAs \sum_i A_i>0$$ 
	satisfy \Assumptions~\ref{ass:smoothness_of_global_function},~\ref{ass:global_lipschitzianity} and~\ref{ass:phi_p_xi_is_locally_uniformly_lipschitz} for $h'_i(x)=\nabla^2  f'_i(x)$.
	\label{thm:quadratic_costs_that_satisfy}
\end{theorem}

% -----------------------------------------------------
\section{Newton-Raphson Consensus}
\label{sec:newton_raphson_consensus}

In this section we provide an algorithm to distributively compute the minimizer of the function $x^*$ defined in~\eqref{equ:definition_of_global_optimum}. The algorithm will be shown to converge to $x^*$ even if $x^*\neq 0$. The proof of convergence will be based on the results derived in the previous sections via a suitable translation of the argument of the cost functions, which basically reduces the problem to the special case $x^*=0$.

Consider then \Algorithm~\ref{alg:distributed_NR}, where $g \big( \bm{x}(-1) \big) = \bm{0}$ and $h \big( \bm{x}(-1) \big) = \bm{0}$ in the initialization step should be intended as initialization of suitable registers and not as operations involving the quantity $\bm{x}(-1)$.

\begin{algorithm}
\caption{\acf{NRC}}
    \begin{algorithmic}[1]
        \Statex \GrayText{(storage allocation and constraints on the parameters)}
		\State $x_i(k), y_i(k) \in \Reals^{n}$ and $z_i(k) \in \Reals^{n \times n}$ for all $k$ and $i=1,\ldots,N$; $\varepsilon \in ( 0, 1 ], c>0$
		\Statex \GrayText{(initialization)}
        \State
		$x_i(0) =0$; $y_i(0) = g_i \big( x_i(-1) \big) = 0$; $z_i(0) = h_i \big( x_i(-1) \big) =0$
            \label{step:NR:scalar:initialization}
        \Statex \GrayText{(main algorithm)}
        \For{$k = 1, 2, \ldots$}
            \For{$i = 1, \ldots, N$}
				\State
				$
					\displaystyle
					x_i(k)
					=
					(1 - \varepsilon) x_i(k-1)
					+
					\varepsilon \cmax{z_i(k-1)}^{-1} y_i(k-1)
				$
				\label{step:NR:scalar:update_of_x}
				\State $\displaystyle y_i(k) \!= \!\sum_{j=1}^N p_{ij} \Big( y_j(k \!- \!1) \!+ \!g_j \big( x_j(k\!-\!1) \big) \!-\! g_j \big( x_j(k\!-\!2) \big) \Big)$
				\label{step:NR:scalar:communication_of_y}
				\State $\displaystyle z_i(k) \!=\! \sum_{j=1}^N p_{ij} \Big( z_j(k\! -\! 1) \!+ \!h_j \big( x_j(k\!-\!1) \big)\! - \!h_j \big( x_j(k\!-\!2) \big) \Big)$
				\label{step:NR:scalar:communication_of_z}
			\EndFor
        \EndFor
    \end{algorithmic}
\label{alg:distributed_NR}
\end{algorithm} 

Intuitively, the algorithm functions as follows: if the dynamics of the $x_{i}(k)$s is sufficiently slow w.r.t.\ the dynamics of the $y_{i}(k)$s and $z_{i}(k)$s, then the two latter quantities tend to reach consensus. Then, the more these quantities reach consensus, the more the products $\cmax{z_i(k)}^{-1} y_i(k)$ exhibit these two specific characteristics: \emph{i)} being the same among the various agent; \emph{ii)} representing Newton descent directions. Thus, the more the $y_{i}(k)$s and $z_{i}(k)$s in \Algorithm~\ref{alg:distributed_NR} are sufficiently close, the more the various $x_{i}(k)$s are driven by the same forcing term, that makes them converge to the same value, equal to the optimum $x^{\ast}$.

We now characterize the convergence properties of \Algorithm~\ref{alg:distributed_NR}. Let us define
$$
\begin{array}{c}
	\displaystyle
	\xi^{y} \DefinedAs \frac{1}{N} \sum_{i=1}^N \left( y_i(0)-g_i \big( x_i(-1) \big) \right) \\
	\displaystyle
	\xi^{z} \DefinedAs \frac{1}{N} \sum_{i=1}^N \left( z_i(0)-h_i \big( x_i(-1) \big) \right) ,
\end{array}
$$
then we have the following theorem:

\begin{theorem}
	Consider the dynamics defined by \Algorithm~\ref{alg:distributed_NR} with possibly nonzero initial conditions. If $\xi^{y} = 0$ and $\xi^{z} = 0$, then under \Assumptions~\ref{ass:smoothness_of_global_function} and~\ref{ass:global_lipschitzianity} there exists a positive scalar $\overline{\varepsilon} > 0$ such that \Theorem~\ref{thm:continuous_to_discrete} holds, i.e., the algorithm can be considered a forward-Euler discretization of a globally exponentially stable continuous dynamics. Thus the local estimates $x_i(k)$ produced by the algorithm exponentially converge to the global minimizer, i.e., 
    $$ \lim_{k \to \infty} x_i(k) = x^* \quad \forall i = 1, \ldots, N $$
	for all $\varepsilon \in(0, \overline{\varepsilon})$ and $x_{i}(0) \in \Reals^{n}$.
\label{thm:distributed_NR_converges_to_the_global_optimum}
\end{theorem}

Consider now that, due to finite-precision issues, the quantities $\xi^{y}$ and $\xi^{z}$ may be non-null. Non-null initial $\xi^{y}$ and $\xi^{z}$ will make the proposed algorithm converge to a point that, in general does not coincide with the global optimum $x^{\ast}$. Nonetheless in this case the computed solution, as a function of the initial conditions, is a smooth function and thus small errors in the initial conditions do not produce dramatic errors in the computation of the optimum:
\begin{theorem}
	Consider the dynamics defined by \Algorithm~\ref{alg:distributed_NR} with possibly nonzero initial $\xi^{y}$ and $\xi^{z}$ but generic $x_{i}(0)$'s. Under \Assumptions~\ref{ass:smoothness_of_global_function},~\ref{ass:global_lipschitzianity} and~\ref{ass:phi_p_xi_is_locally_uniformly_lipschitz} there exist positive scalars $a, r,\overline{\varepsilon}$ and a continuously differentiable function $\Psi : \Reals^n\times\Reals^{n\times n}\mapsto \Reals^n$ satisfying
	$$
		\left\| \Psi(\xi^{y},\xi^{z})-x^* \right\| \leq a \left( \left\| \xi^{y} \right\| + \left\| \xi^{z} \right\| \right)
	$$
	s.t.\ the local estimates exponentially converge to it, i.e., 
	$$
		\lim_{k \to \infty} x_i(k) = \Psi \left( \xi^{y}, \xi^{z} \right)
		\quad \forall i = 1, \ldots, N
	$$
	for all $\varepsilon \in(0, \overline{\varepsilon})$, initial conditions $x_{i}(0) \in \Reals^{n}$ and $\left( \left\| \xi^{y} \right\| + \left\| \xi^{z} \right\| \right) \leq r$.
\label{thm:distributed_NR_converges_to_neighborhood_of_the_global_optimum}
\end{theorem}

We notice that \Theorem~\ref{thm:distributed_NR_converges_to_neighborhood_of_the_global_optimum} ensures global convergence properties w.r.t.\ the initial conditions $x_{i}(0)$'s by requiring \Assumptions~\ref{ass:smoothness_of_global_function},~\ref{ass:global_lipschitzianity} and~\ref{ass:phi_p_xi_is_locally_uniformly_lipschitz}, while for the same convergence properties \Theorem~\ref{thm:distributed_NR_converges_to_the_global_optimum} requires only \Assumptions~\ref{ass:smoothness_of_global_function} and~\ref{ass:global_lipschitzianity}. The difference is that \Theorem~\ref{thm:distributed_NR_converges_to_neighborhood_of_the_global_optimum} considers a non-null perturbation $\xi$ and \Assumption~\ref{ass:phi_p_xi_is_locally_uniformly_lipschitz} is needed to cope with this additional perturbation term.

The \Assumptions~\ref{ass:smoothness_of_global_function},~\ref{ass:global_lipschitzianity} and~\ref{ass:phi_p_xi_is_locally_uniformly_lipschitz} are not needed if only local convergence is ought. In fact, local differentiability, and therefore local Lipschitzianity,  of the cost functions $f_i(x)$ at the minimizer $x^*$ is sufficient to guarantee that \Assumptions~\ref{ass:global_lipschitzianity} and~\ref{ass:phi_p_xi_is_locally_uniformly_lipschitz} are locally valid. As so, the proof that the equilibrium point is a locally exponentially stable point is exactly the same, with the difference that all bounds and inequalities are local. This observation is summarized in the following theorem.
\begin{theorem}
	Consider the dynamics defined by \Algorithm~\ref{alg:distributed_NR} with possibly nonzero initial conditions. Under the assumptions that the $f_i$'s are $\mathcal{C}^3$ and that $\nabla^2 \overline{f}(x^*) \geq c I$, there exist positive scalars $a, r, \overline{\varepsilon}$ and a continuously differentiable function $\Psi : \Reals^n\times\Reals^{n\times n}\mapsto \Reals^n$ s.t.
	$$
		\lim_{k \to \infty} x_i(k) = \Psi \left( \xi^{y}, \xi^{z} \right)
		\quad \forall i = 1, \ldots, N
	$$
	and satisfying
	$$
		\left\| \Psi(\xi^{y},\xi^{z})-x^* \right\|
		\leq
		a \big( \left\| \xi^{y} \right\| + \left\| \xi^{z} \right\| \big)
	$$
	for all $\varepsilon \in(0, \overline{\varepsilon})$ and initial conditions
	$$
		\left\| x_i(0) - x^* \right\| \leq r, \quad
		\left\| y_i - \overline{\overline{g}}(x^*) \right\| \leq r, \quad
		\left\| z_i - \overline{\overline{h}}(x^*) \right\| \leq r
	$$
	$$
		\left\| g_i(x_i(-1)) - \overline{\overline{g}}(x^*) \right\| \leq r, \quad
		\left\| h_i(x_i(-1)) - \overline{\overline{h}}(x^*) \right\| \leq r .
	$$
\label{thm:local_stability}
\end{theorem}

Numerical simulations suggest that the algorithm is robust w.r.t.\ numerical errors and quantization noise. We also notice that \Theorem~\ref{thm:distributed_NR_converges_to_the_global_optimum} guarantees the existence of a critical value $\overline{\varepsilon}$ but does not provide indications on its value. This is a known issue in all the systems dealing with separation of time scales. A standard rule of thumb is then to let the rate of convergence of the fast dynamics be sufficiently faster than the one of the slow dynamics, typically 2-10 times faster. In our algorithm the fast dynamics inherits the rate of convergence of the consensus matrix $P$, given by its spectral gap $\sigma(P)$, i.e., its spectral radius $\rho(P) = 1 - \sigma(P)$. The rate of convergence of the slow dynamics is instead governed by~\eqref{equ:NR_evolution_for_average_x}, which is nonlinear and therefore possibly depending on the initial conditions. However, close to the equilibrium point the dynamic behavior is approximately given by $\dot{\overline{x}}(t) \approx - \big( \overline{x}(t) - x^* \big)$, thus, since $x_i(k) \approx \overline{x}(\varepsilon k)$, then the convergence rate of the algorithm approximately given by $1 - \varepsilon$.

Thus we aim to let
$
    1 - \rho(P) \gg 1 - ( 1 - \varepsilon ) ,
$
which provides the rule of thumb
\begin{equation}
    \varepsilon \ll \sigma(P) \; .
\label{eqn:rule}
\end{equation} 
which is suitable for generic cost functions. We then notice that, although the spectral gap $\sigma(P)$ might not be known in advance, it is possible to distributedly estimate it, see, e.g.,~\cite{sahai_et_al__2012__hearing_the_clusters_of_a_graph__a_distributed_algorithm}. However, such rule of thumb might be very conservative. In fact, if all the $f_i$'s are quadratic and are, w.l.o.g.\ s.t.\ $\nabla^{2} f_{i} \geq c I$, then one can set $\varepsilon = 1$ and neglect the thresholding $\left[ \cdot \right]_{c}$, so that the procedure reduces to
\begin{equation}
	\begin{array}{rcl}
		\bm{x}(k+1) & = & \displaystyle \frac{\bm{y}(k)}{\bm{z}(k)} \\
		\bm{y}(k+1) & = & ( P \otimes I_{n} ) \bm{y}(k) \\
		\bm{z}(k+1) & = & ( P \otimes I_{n} ) \bm{z}(k) \; .
	\end{array}
\label{equ:minimization_strategy:quadratic_case}
\end{equation}
where $\bm{x}(k) \DefinedAs \left[ x_{1}^T(k) \;, \ldots, \; x_{N}^T(k) \right]^T$, $\bm{y}(k) \DefinedAs \left[ y_{1}^T(k) \;, \ldots, \; y_{N}^T(k) \right]^T$, $\bm{z}(k) \DefinedAs \left[ z_{1}(k) \;, \ldots, \; z_{N}(k) \right]^T$. Thus:
\begin{theorem}
	Consider \Algorithm~\ref{alg:distributed_NR} with arbitrary initial conditions $x_{i}(0)$, quadratic cost functions $f_i = \frac{1}{2} \left( x - d_{i} \right)^{T} A_i \left( x - d_{i} \right)$ with $A_{i} > 0$ and $\varepsilon = 1$. Then
    $
		\left\| x_{i}(k) - x^{\ast} \right\|
        \leq
        \alpha \left( \rho(P) \right)^{k}
    $
    for all $k, i$ and for a suitable positive $\alpha$.
\label{thm:convergence_for_quadratic_functions}
\end{theorem}
Thus, if the cost functions are close to be quadratic then the overall rate of convergence is limited by the rate of convergence of the embedded consensus algorithm. Moreover, the values of $\varepsilon$ that still guarantee convergence can be much larger than those dictated by the rule of thumb~\eqref{eqn:rule}.

\subsection{On the selection of the structure of $h(x)$}
\label{ssec:on_the_selection_of_the_structure_of_hx}

As introduced in \Section~\ref{ssec:stability_of_single_agent_nr_dynamics}, by selecting different structures for $h_i(x)$ one can obtain different procedures with different convergence properties and different computational/communication requirements. Plausible choices for $h_i$ are the ones in~\eqref{equ:choice_of_h_i}, and the correspondences are the following:
	
\vspace{0.2cm}

\noindent $\bullet$ $h_{i}(x) = \nabla^{2} f_{i}(x)$ $\rightarrow$ \acf{NRC}: in this case it is possible to rewrite the main algorithm and show that, for sufficiently small $\varepsilon$, $x_i(k) \approx \overline{x}(\varepsilon k)$, where $\overline{x}(t)$ evolves according to the continuous-time Newton-Raphson dynamics
$$%\begin{equation}
    \dot{\overline{x}}(t)
    =
    -
    \Big[
        \nabla^{2} \overline{f} \big( \overline{x}(t) \big)
    \Big]^{-1}
    \nabla \overline{f} \big( \overline{x}(t) \big) \; .
    %
%\label{equ:NR:dynamics_of_the_average_component:multidimensional_case}
$$%\end{equation}

\vspace{0.2cm}

\noindent $\bullet$ $h_{i}(x) = \mathrm{diag} \left[ \nabla^2 f_{i}(x) \right]$ $\rightarrow$ \acf{JC}: choice~$h_{i}(x) = \nabla^{2} f_{i}(x)$ requires agents to exchange information on $\BigOOf{n^2}$ scalars, and this could pose problems under heavy communication bandwidth constraints and large $n$'s. Choice $h_{i}(x) = \mathrm{diag} \left[ \nabla^2 f_{i}(x) \right]$ instead reduces the amount of information to be exchanged via the underlying diagonalization process, also called Jacobi approximation\footnote{In centralized approaches, nulling the Hessian's off-diagonal terms is a well-known procedure, see, e.g.,~\cite{becker_le_cun__1988__improving_the_convergence_of_back_propagation_learning_with_second_order_models}. See also~\cite{athuraliya_low__2000__optimization_flow_control_with_newton_like_algorithm,zargham_et_al__2011__accelerated_dual_descent_for_network_optimization} for other Jacobi algorithms with different communication structures.}. In this case, for sufficiently small $\varepsilon$, $x_i(k) \approx \overline{x}(\varepsilon k)$, where $\overline{x}(t)$ evolves according to the continuous-time dynamics
$$%\begin{equation}
    \dot{\overline{x}}(t)
    =
    -
    \left(
        \mathrm{diag}
        \Big[
			\nabla^{2}	  \overline{f}	  \big(    \overline{x}(t)    \big)
        \Big]
    \right)^{-1}
    \nabla \overline{f} \big( \overline{x}(t) \big) \; ,
    %
%\label{equ:NR:dynamics_of_the_average_component:multidimensional_case}
$$%\end{equation}
which can be shown to converge to the global optimum $x^{\ast}$ with a convergence rate that in general is slower than the Newton-Raphson when the global cost function is skewed.

\vspace{0.2cm}

\noindent $\bullet$ $h_{i}(x) = I$ $\rightarrow$ \acf{GDC}: this choice is motivated in frameworks where the computation of the local second derivatives $\displaystyle \left. \frac{\partial^2 f_i}{\partial x_{m}^2} \right|_{x}$ is expensive (with $x_{m}$ indicating here the $m$-th component of $x$), or where the second derivatives simply might not be continuous. With this choice the main algorithm reduces to a distributed gradient-descent procedure. In fact, for sufficiently small $\varepsilon$, $x_i(k) \approx \overline{x}(\varepsilon k)$ with $\overline{x}(t)$ evolving according to the continuous-time dynamics
$$%\begin{equation}
    \dot{\overline{x}}(t)
    =
    -
    \nabla \overline{f} \big( \overline{x}(t) \big) \; ,
    %
%\label{equ:NR:dynamics_of_the_average_component:multidimensional_case}
$$%\end{equation}
which one again is guaranteed to converge to the global optimum $x^{\ast}$.

\vspace{0.15cm}

The following \Table~\ref{tab:computational_communication_and_memory_costs} summarizes the various costs of the previously proposed strategies.

\begin{table}[!htbp]
	\centering
\begin{tabular}
{
	>{\arraybackslash  \columncolor{black!05!white}}   p{0.32\columnwidth}	 <{}
    >{\centering \arraybackslash}                               p{0.13\columnwidth} <{}
    >{\centering \arraybackslash \columncolor{black!05!white}}  p{0.13\columnwidth} <{}
    >{\centering \arraybackslash}                               p{0.13\columnwidth} <{}
}
\toprule % --------------------------------------------

Choice & 
NRC, $h_{i}(x) = \nabla^{2} f_{i}(x)$ &
JC, $h_{i}(x) = \mathrm{diag} \left[ \nabla^2 f_{i}(x) \right]$ &
GDC, $h_{i}(x) = I$ \\

\midrule

Computational Cost  &	$\BigOOf{n^3}$	&	$\BigOOf{n}$	&   $\BigOOf{n}$	\\
Communication Cost  &	$\BigOOf{n^2}$	&	$\BigOOf{n}$	&   $\BigOOf{n}$	\\
Memory Cost			&	$\BigOOf{n^2}$	&	$\BigOOf{n}$ 	&	$\BigOOf{n}$	\\

\bottomrule % --------------------------------------------
\end{tabular}
\caption{Computational, communication and memory costs of NRC, JC, GDC per single unit and single step.}
\label{tab:computational_communication_and_memory_costs}
\end{table}

We remark that $\overline{\varepsilon}$ in \Theorem~\ref{thm:distributed_NR_converges_to_the_global_optimum} depends also on the particular choice for $h_i$. The list of choices for $h_i$ given above is not exhaustive. For example, future directions are to implement distributed quasi-Newton procedures. To this regard, we recall that approximations of the Hessians that do not maintain symmetry and positive definiteness or are bad conditioned require additional modification steps, e.g., through Cholesky factorizations~\cite{golub_van_loan__1996__matrix_computations}.

Finally, we notice that in scalar scenarios \ac{JC} and \ac{NRC} are equivalent, while \ac{GDC} corresponds to algorithms requiring just the knowledge of first derivatives.

% -----------------------------------------------------
\section{Numerical Examples}
\label{sec:numerical_examples}

    In \Section~\ref{ssec:effects_of_varepsilon} we analyze the effects of different choices of $\varepsilon$ on the \ac{NRC} on regular graphs and exponential cost functions. We then propose two machine learning problems in \Section~\ref{ssec:optimization_problems}, used in \Sections~\ref{ssec:comparison_of_NRC_JC_GD} and~\ref{ssec:comparison_with_other_algorithms}, and numerically compare the convergence performance of the \ac{NRC}, \ac{JC}, \ac{GDC} algorithms and other distributed convex optimization algorithms on random geometric graphs.

	Notice that we will use cost functions that may not satisfy \Assumptions~\ref{ass:smoothness_of_global_function},~\ref{ass:global_lipschitzianity} and~\ref{ass:phi_p_xi_is_locally_uniformly_lipschitz} to highlight the fact that the algorithm seems to have favorable numerical properties and large basins of stability even if the assumptions needed for global stability are not satisfied.

\subsection{Effects of the choice of $\varepsilon$}
\label{ssec:effects_of_varepsilon}

Consider a ring network of $S = 30$ agents that communicate only to their left and right neighbors through the consensus matrix
\begin{equation}
    P =
    \begin{bmatrix}
        0.5     & 0.25      &           &           & 0.25  \\
        0.25    & 0.5       & 0.25      &           &       \\
                & \ddots    & \ddots    & \ddots    &       \\
                &           & 0.25      & 0.5       & 0.25  \\
        0.25    &           &           & 0.25      & 0.5   \\
    \end{bmatrix} ,
\label{equ:P:symmetric_circulant_with_weight_05}
\end{equation}
so that the spectral radius $\rho(P) \approx 0.99$, implying a spectral gap $\sigma(P) \approx 0.01$. Consider also scalar costs of the form
$
    f_i(x)
    =
        c_i e^{ a_i x }
    +   d_i e^{ -b_i x },
$
$
    i = 1, \ldots, N,
$ 
with $a_i, b_i \sim \UniformDistribution{0}{0.2}$, $c_i, d_i \sim \UniformDistribution{0}{1}$ and where $\mathcal{U}$ indicates the uniform distribution.

\begin{figure}[!htbp]
	\centering
    \includegraphics{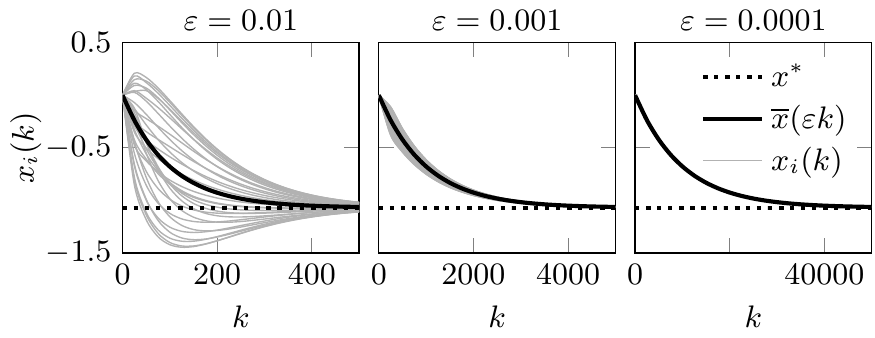}
\caption{Temporal evolution of system~\eqref{equ:augmented_continuous_time_system} for different values of $\varepsilon$, with $N = 30$. The black dotted line indicates $x^{\ast}$. The black solid line indicates the slow dynamics $\overline{x}(\varepsilon k)$ of \Equation~\eqref{equ:NR_evolution_for_average_x}. As $\varepsilon$ decreases, the difference between the time scale of the slow and fast dynamics increases, and the local states $x_i(k)$ converge to the manifold of $\overline{x}(\varepsilon k)$.}
\label{fig:epsiloncomparison:time_evolution_of_states}
\end{figure}

\Figure~\ref{fig:epsiloncomparison:time_evolution_of_states} compares the evolution of the local states $x_i$ of the continuous system~\eqref{equ:augmented_continuous_time_system} for different values of $\varepsilon$. When $\varepsilon$ is not sufficiently small, then the trajectories of $x_i(t)$ are different even if they all start from the same initial condition $x_i(0) = 0$. As $\varepsilon$ decreases, the difference between the two time scales becomes more evident and all the trajectories $x_i(k)$ become closer to the trajectory given by the slow \ac{NR} dynamics $\overline{x}(\varepsilon k)$ given in~\eqref{equ:NR_evolution_for_average_x} and guaranteed to converge to the global optimum $x^{\ast}$.

\begin{figure}[!htbp]
	\centering
    \subfigure
    [Time evolution of the local states $x_i(k)$ with $v(0) = w(0) = y(0) = z(0) = 0$ and $x_i(0) \sim \UniformDistribution{-2}{2}$.]
    {
        \includegraphics{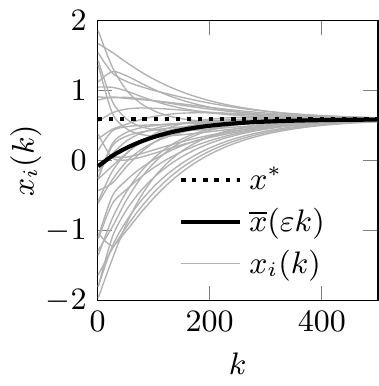}
        \label{fig:stability:time_evolution_of_states_not_null_initialization}
    }
    $\quad$
    \subfigure
    [Empirical distribution of the errors $x_i(+\infty) - x^{\ast}$ under artificially perturbed initial conditions $\alpha(0),\beta(0) \sim \UniformDistribution{-\sigma}{\sigma}$ for different values of $\sigma$.]
    {
      	\includegraphics{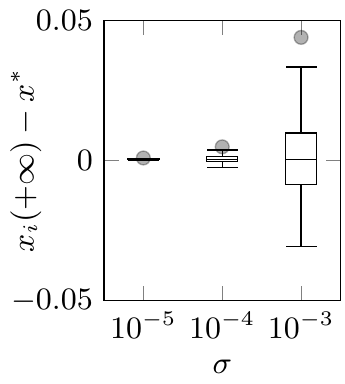}
        \label{fig:robustness:boxplots}
    }
\caption{Characterization of the dependency of the performance of \Algorithm~\ref{alg:distributed_NR} on the initial conditions. In all the experiments $\varepsilon = 0.01$ and $N = 30$.}
\label{fig:robustness}
\end{figure}

In \Figure~\ref{fig:robustness} we address the robustness of the proposed algorithm w.r.t.\ the choice of the initial conditions. In particular, \Figure~\ref{fig:stability:time_evolution_of_states_not_null_initialization} shows that if $\alpha = \beta = 0$ then the local states $x_i(t)$ converge to the optimum $x^{\ast}$ for arbitrary initial conditions $x_i(0)$. \Figure~\ref{fig:robustness:boxplots} considers, besides different initial conditions $x_i(0)$, also perturbed initial conditions $v(0)$, $w(0)$, $y(0)$, $z(0)$ leading to non null $\alpha$'s and $\beta$'s. More precisely we apply \Algorithm~\ref{alg:distributed_NR} to different random initial conditions s.t.
$
    \alpha, \beta
    \sim
    \UniformDistribution{-\sigma}{\sigma}
$. 
\Figure~\ref{fig:robustness:boxplots} shows the boxplots of the errors $x_i(+\infty) - x^{\ast}$ for different $\sigma$'s based on 300 Monte Carlo runs with $\varepsilon = 0.01$ and $N = 30$.

\begin{figure}[!htbp]
	\centering
    \subfigure
    [Relative \ac{MSE} at a given time $k$ as a function of the parameter $\varepsilon$ for classification problem~\eqref{equ:spam_nonspam_classification}.]
    {
		\includegraphics[scale = 0.9]{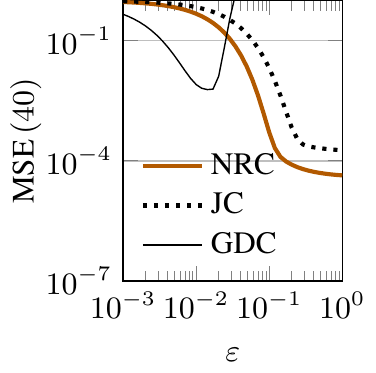}
       \label{fig:spam_nonspam__RMSE_at_40_vs_parameters__NRC_JC_GD}
    }
    $\quad$
    \subfigure
    [Relative \ac{MSE} as a function of the time $k$, with the parameter $\varepsilon$ chosen as the best from \Figure~\ref{fig:spam_nonspam__RMSE_at_40_vs_parameters__NRC_JC_GD} for classification problem~\eqref{equ:spam_nonspam_classification}.] 
    {
		\includegraphics[scale = 0.9]{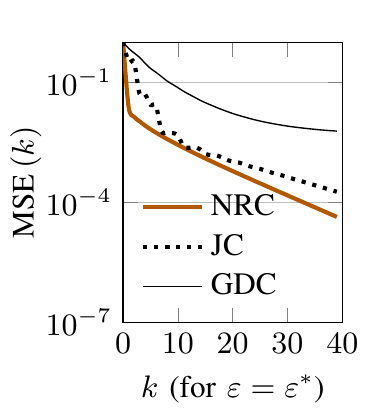}
       \label{fig:spam_nonspam__RMSE_evolution__NRC_JC_GD}
    }
    \\
    \subfigure
    [Relative \ac{MSE} at a given time $k$ as a function of the parameter $\varepsilon$ for regression problem~\eqref{equ:housing_regression}.]
    {
		\includegraphics[scale = 0.9]{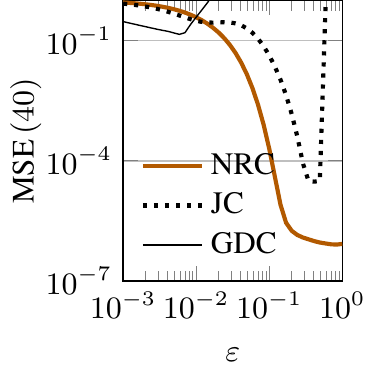}
       \label{fig:housing__RMSE_at_40_vs_parameters__NRC_JC_GD}
    }
    $\quad$
    \subfigure
    [Relative \ac{MSE} as a function of the time $k$, with the parameter $\varepsilon$ chosen as the best from \Figure~\ref{fig:housing__RMSE_at_40_vs_parameters__NRC_JC_GD} for regression problem~\eqref{equ:housing_regression}.] 
    {
		\includegraphics[scale = 0.9]{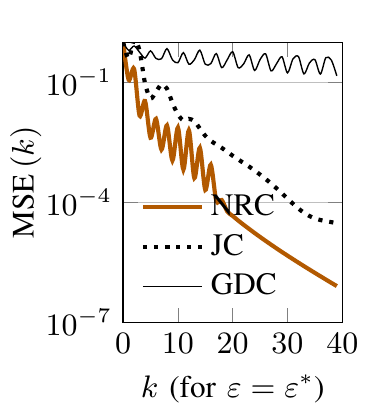}
       \label{fig:housing__RMSE_evolution__NRC_JC_GD}
    }
\caption{Convergence properties of \Algorithm~\ref{alg:distributed_NR} for the problems described in \Section~\ref{ssec:optimization_problems} and for different choices of $h_{i}(\cdot)$. Choice $h_{i}(x) = \nabla^{2} f_{i}(x)$ corresponds to the \ac{NRC} algorithm, $h_{i}(x) = \mathrm{diag} \left[ \nabla^2 f_{i}(x) \right]$ to the \ac{JC}, $h_{i}(x) = I$ to the \ac{GDC}.}
\label{fig:comparison_of_NRC_JC_and_GD}
\end{figure}

\subsection{Optimization problems}
\label{ssec:optimization_problems}

\begin{figure}[!htbp]
    \centering
    \includegraphics{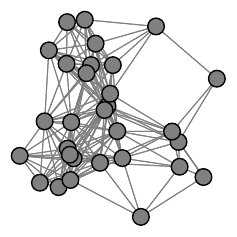}
    \caption{Random geometric graph exploited in the simulations relative to the optimization problem~\eqref{equ:spam_nonspam_classification}. For this graph $\rho(P) \approx 0.9338$, with $P$ the matrix of Metropolis weights.}
    \label{fig:random_geometric_for_spamnonspam}
\end{figure}

The first problem considered is the distributed training of a Binomial-Deviance based classifier, to be used, e.g., for spam-nonspam classification tasks \cite[Chap.~10.5]{hastie_et_al__2001__the_elements_of_statistical_learning}. More precisely, we consider a database of emails $E$, where $j$ is the email index, $y_{j} = -1, 1$ denotes if the email $j$ is considered spam or not, $\chi_{j} \in \Reals^{n-1}$ numerically summarizes the $n-1$ features of the $j$-th email (how many times the words ``money'', ``dollars'', etc., appear). If the $E$ emails come from different users that do not want to disclose their private information, then it is meaningful to exploit the distributed optimization algorithms described in the previous sections. More specifically, letting $x=\left( x', x_{0} \right) \in \Reals^{n-1} \times \Reals$ represents a generic classification hyperplane, training a Binomial-Deviance based classifier corresponds to solve a distributed optimization problem where the local cost functions are given by:
\begin{equation}
    f_i\left( x \right)
    \DefinedAs
    \sum_{j \in E_i}
    \log
    \Big(
        1 
        +
        \exp
        \left(
			-	y_{j}	\left(	 \chi_{j}^T   x'	+	x_{0}	 \right)
        \right)
    \Big)
    +
    \gamma \left\| x' \right\|^{2}_{2} .
    \label{equ:spam_nonspam_classification}
\end{equation}
where $E_i$ is the set of emails available to agent $i$, $E=\cup_{i=1}^N E_i$, and $\gamma$ is a global regularization parameter. In the following numerical experiments we consider $|E| = 5000$ emails from the spam-nonspam UCI repository, available at \texttt{http://archive.ics.uci.edu/ml/datasets/Spambase}, randomly assigned to 30 different users communicating as in graph of \Figure~\ref{fig:random_geometric_for_spamnonspam}. For each email we consider $3$ features (the frequency of words ``make'', ``address'', ``all'') so that the corresponding optimization problem is 4-dimensional.

The second problem considered is a regression problem inspired by the UCI Housing dataset available at \texttt{http://archive.ics.uci.edu/ml/datasets/Housing}. In this task, an example $\chi_{j}\in\Reals^{n-1}$ is a vector representing some features of a house (e.g., per capita crime rate by town, index of accessibility to radial highways, etc.), and $y_{j}\in\Reals$ denotes the corresponding median monetary value of of the house. The objective is to obtain a predictor of house value based on these data. Similarly as the previous example, if the datasets come from different users that do not want to disclose their private information, then it is meaningful to exploit the distributed optimization algorithms described in the previous sections. This problem can be formulated as a convex regression problem on the local costs
\begin{equation}
    f_i\left( x \right)
    \DefinedAs
    \sum_{j \in E_i}
    \frac
    {
        \left(
            y_{j} - \chi_{j}^T x' - x_{0} 
        \right)^{2}
    }
    {
        \left|
            y_{j} - \chi_{j}^T x' - x_{0} 
        \right|
        +
        \beta
    }
    +
    \gamma \left\| x' \right\|^{2}_{2} .
    \label{equ:housing_regression}
\end{equation}
where $x= (x', x_{0}^{\ast})\in\Reals^{n-1}\times\Reals$ is the vector of coefficient for the linear predictor $\widehat{y}=\chi^T x' + x_0$ and $\gamma$ is a common regularization parameter. The loss function $\frac{(\cdot)^{2}}{|\cdot| + \beta}$ corresponds to a smooth $\mathcal{C}^{2}$ version of the Huber robust loss, a loss that is usually employed to minimize the effects of outliers. In our case $\beta$ dictates for which arguments the loss is pseudo-linear or pseudo-quadratic and has been manually chosen to minimize the effects of outliers. In our experiments we used $4$ features, $\beta = 50$, $\gamma=1$, and $|E| = 506$ total number of examples in the dataset randomly assigned to the $N=30$ users communicating as in the graph of \Figure~\ref{fig:random_geometric_for_spamnonspam}.

In both the previous problems the optimum, in the following indicated for simplicity with $x^{\ast}$, has been computed with a centralized \ac{NR} with the termination rule ``stop when in the last 5 steps the norm of the guessed $x^{\ast}$ changed less than $10^{-9}\%$''.

\subsection{Comparison of the \ac{NRC}, \ac{JC} and \ac{GDC} algorithms}
\label{ssec:comparison_of_NRC_JC_GD}

In \Figure~\ref{fig:comparison_of_NRC_JC_and_GD} we analyze the performance of the three proposed \ac{NRC}, \ac{JC} and \ac{GDC} algorithms defined by the various choices for $h_{i}(x)$ in \Algorithm~\ref{alg:distributed_NR} in terms of the relative \ac{MSE}
$$
    \MSEOf{k} 
    \DefinedAs 
    \frac{1}{N}
    \sum_{i = 1}^{N} 
    {\left\| x_{i}(k) - x^{*} \right\|^{2}} / {\left\| x^{*} \right\|^{2}}
$$
for the classification and regression optimization problem described above. The consensus matrix $P$ has been by selecting the Metropolis-Hastings weights which are consistent with the communication graph \cite{xiao_et_al__2007__distributed_average_consensus_with_least_mean_square_deviation}. Panels~\ref{fig:spam_nonspam__RMSE_at_40_vs_parameters__NRC_JC_GD} and~\ref{fig:housing__RMSE_at_40_vs_parameters__NRC_JC_GD} report the MSE obtained at a specific iteration ($k = 40$) by the various algorithms, as a function of $\varepsilon$. These plots thus inspect the sensitivity w.r.t.\ the choice of the tuning parameters. Consistently with the theorems in the previous section, the \ac{GDC} and \ac{JC} algorithms are stable only for $\varepsilon$ sufficiently small, while \ac{NRC} exhibit much larger robustness and best performance for $\varepsilon=1$. Panels~\ref{fig:spam_nonspam__RMSE_evolution__NRC_JC_GD} and~\ref{fig:housing__RMSE_evolution__NRC_JC_GD} instead report the evolutions of the relative \ac{MSE} as a function of the number of iterations $k$ for the optimally tuned algorithms.

We notice that the differences between \ac{NRC} and \ac{JC} are evident but not resounding, due to the fact that the Jacobi approximations are in this case a good approximation of the analytical Hessians. Conversely, \ac{GDC} presents a slower convergence rate which is a known drawback of gradient descent algorithms.

\subsection{Comparisons with other distributed convex  optimization  algorithms}
\label{ssec:comparison_with_other_algorithms}

We now compare \Algorithm~\ref{alg:distributed_NR} and its accelerated version, referred as \ac{FNRC} and described in detail below in \Algorithm~\ref{alg:distributed_FNRC}), with three popular distributed convex optimization methods, namely the \ac{DSM}, the \ac{DCM} and the \ac{ADMM}, described respectively in \Algorithm~\ref{alg:distributed_subgradient},~\ref{alg:DCM} and~\ref{alg:ADMM}. The following discussion provides some details about these strategies.

$\bullet$ \ac{FNRC} is an accelerated version of \Algorithm~\ref{alg:distributed_NR} that inherits the structure of the so called \emph{second order diffusive schedules}, see, e.g.,~\cite{muthukrishnan_et_al__1998__first_and_second_order_diffusive_methods_for_rapid_coarse_distributed_load_balancing}, and exploits an additional level of memory to speed up the convergence properties of the consensus strategy. Here the weights multiplying the $g_{i}$'s and $h_i$'s are necessary to guarantee exact tracking of the current average, i.e., $\sum_{i} y_i(k) = \sum_{i} g_{i}\big(x(k-1)\big)$ for all $k$. As suggested in~\cite{muthukrishnan_et_al__1998__first_and_second_order_diffusive_methods_for_rapid_coarse_distributed_load_balancing}, we set the $\varphi$ that weights the gradient and the memory to $\displaystyle \varphi = \frac{2}{1 + \sqrt{1 - \rho(P)^2}}$. This guarantees second order diffusive schedules to be faster than first order ones (even if this does not automatically imply the \ac{FNRC} to be faster than the \ac{NRC}). This setting can be considered a valid heuristic to be used when $\rho(P)$ is known. For the graph in \Figure~\ref{fig:random_geometric_for_spamnonspam}, $\varphi \approx 1.4730$.
\begin{algorithm}
\caption{Fast Newton-Raphson Consensus}
    \begin{algorithmic}[1]
        \State storage allocation, constraints on the parameters and initialization as in \Algorithm~\ref{alg:distributed_NR}
        \For{$k = 1, 2, \ldots$}
        \For{$i = 1, \ldots, N$}
			\State
			$
				\displaystyle
				x_{i}(k)
				=
				(1 - \varepsilon) x_{i}(k-1)
				+
				\varepsilon
				\cmax{z_{i}(k-1)}^{-1}
				y_{i}(k-1)
			$
			\State
			$
				\displaystyle
				\widetilde{y}_{i}(k)
				=
												y_{i}(k - 1)
				+			\frac{1}{\varphi}	g_{i} \big(	x_{i}(k - 1) \big)
				-                               g_{i} \big( x_{i}(k - 2) \big)
				- \frac{1 - \varphi}{\varphi}   g_{i} \big( x_{i}(k - 3) \big)
			$
			\State
			$
				\displaystyle
				\widetilde{z}_{i}(k)
				=
												z_{i}(k - 1)
				+			\frac{1}{\varphi}	h_{i} \big(	x_{i}(k - 1) \big)
				-                               h_{i} \big( x_{i}(k - 2) \big)
				- \frac{1 - \varphi}{\varphi}   h_{i} \big( x_{i}(k - 3) \big)
			$
			\State  
			$
				\displaystyle
				y_{i}(k)
				=
				\varphi
				\sum_{j = 1}^{N}
				\big(
					p_{ij} \widetilde{y_{j}}(k)
				\big)
				+
				(1 - \varphi)
				y_{i}(k-2)
			$
			\State  
			$
				\displaystyle
				z_{i}(k)
				=
				\varphi
				\sum_{j = 1}^{N}
				\big(
					p_{ij} \widetilde{z_{j}}(k)
				\big)
				+
				(1 - \varphi)
				z_{i}(k-2)
			$
        \EndFor
        \EndFor
    \end{algorithmic}
\label{alg:distributed_FNRC}
\end{algorithm}

$\bullet$ \ac{DSM}, as proposed in~\cite{nedic_ozdaglar__2009__distributed_subgradient_methods_for_multi_agent_optimization}, alternates consensus steps on the current estimated global minimum $x_i(k)$ with subgradient updates of each $x_i(k)$ towards the local minimum. To guarantee the convergence, the amplitude of the local subgradient steps should appropriately decrease. \Algorithm~\ref{alg:distributed_subgradient} presents a synchronous \ac{DSM} implementation, where $\varrho$ is a tuning parameter and $P$ is the matrix of Metropolis-Hastings weights.
\begin{algorithm}
\caption{\ac{DSM}~\cite{nedic_ozdaglar__2009__distributed_subgradient_methods_for_multi_agent_optimization}}
    \begin{algorithmic}[1]

		\Statex \GrayText{(storage allocation and  constraints	on	parameters)}
		\State $x_{i}(k) \in \Reals^{n}$ for all $i$. $\varrho \in \PositiveReals$
        \Statex \GrayText{(initialization)}
        \State $x_{i}(0) = 0$
        \Statex \GrayText{(main algorithm)}
        \For{$k = 0, 1, \ldots$}
            \For{$i = 1, \ldots, N$}
                \State
                $
                    \displaystyle
                    x_i(k+1)
                    =
                    \sum_{j = 1}^{N}
                    p_{ij}
                    \left(
                        x_j(k)
                        -
						\frac{\varrho}{k} \nabla f_{j} \big(  x_j(k)  \big)
                    \right)
                $
            \EndFor
        \EndFor

    \end{algorithmic}
\label{alg:distributed_subgradient}
\end{algorithm}

$\bullet$ \ac{DCM}, as proposed in~\cite{wang_elia__2010__control_approach_to_distributed_optimization}, differentiates from the gradient searching because it forces the states to the global optimum by controlling the subgradient of the global cost. This approach views the subgradient as an input/output map and uses small gain theorems to guarantee the convergence property of the system. Again, each agents $i$ locally computes and exchanges information with its neighbors, collected in the set $\mathcal{N}_i \DefinedAs \{ j \,| \,(i,j) \in \mathcal{E} \}$. \ac{DCM} is summarized in \Algorithm~\ref{alg:DCM}, where $\mu, \nu > 0$ are parameters to be designed to ensure the stability property of the system. Specifically, $\mu$ is chosen in the interval $\displaystyle 0 < \mu < \frac{2}{2 \max_{i = \{1, \ldots, N\}}{|\mathcal{N}_i|} + 1}$ to bound the induced gain of the subgradients. Also here the parameters have been manually tuned for best convergence rates.
\begin{algorithm}
\caption{\ac{DCM}~\cite{wang_elia__2010__control_approach_to_distributed_optimization}}
    \begin{algorithmic}[1]

		\Statex \GrayText{(storage allocation and  constraints	on	parameters)}
        \State
        $x_{i}(k), z_{i}(k) \in \Reals^{n}$, for all $i$. $\mu, \nu \in \PositiveReals$
        \Statex \GrayText{(initialization)}
		\State	$x_{i}(0)	=  z_{i}(0)   =   0$	 for   all	 $i$
        \Statex \GrayText{(main algorithm)}
        \For{$k = 0, 1, \ldots$}
            \For{$i = 1, \ldots, N$}
                \State
                $
                    \displaystyle
                    z_i(k+1)
                    =
                    z_i(k)
                    +
                    \mu \sum_{j \in \mathcal{N}_i} \big( x_i(k) - x_j(k) \big)
                $
                \State
                $
                    \displaystyle
                    x_i(k+1)
                    =
                    x_i(k)
                    +
                    \mu \sum_{j \in \mathcal{N}_i} \big( x_j(k) - x_i(k) \big)
                    +
                    \mu \sum_{j \in \mathcal{N}_i} \big( z_j(k) - z_i(k) \big)
                    -
					\mu  \,  \nu  \,  \nabla   f_i	 \big(	 x_i(k)   \big)
                $
            \EndFor
        \EndFor

    \end{algorithmic}
\label{alg:DCM}
\end{algorithm}

$\bullet$ \ac{ADMM}, instead, requires the augmentation of the system through additional constraints that do not change the optimal solution but allow the Lagrangian formalism. There exist different implementations of \ac{ADMM} in distributed contexts, see, e.g.,~\cite[pp.~253-261]{bertsekas_tsitsiklis__1997__parallel_and_distributed_computation,schizas_et_al__2008__consensus_in_ad_hoc_wsns_with_noisy_links__part_i__distributed_estimation_of_deterministic_signals,boyd_et_al__2010__distributed_optimization_and_statistical_learning_via_the_admm}. For simplicity we consider the following formulation, 
\begin{equation*}
    \begin{array}{l}
        \displaystyle
        \min_{x_1, \ldots, x_N} \;\;
        \sum_{i = 1}^{N} f_i(x_i) \\
        \text{s.t.\ } z_{(i,j)} = x_{i},
        \quad \forall i \in \mathcal{N},
        \quad \forall (i,j) \in \mathcal{E},
    \end{array}
\end{equation*}
where the auxiliary variables $z_{(i,j)}$ correspond to the different links in the network, and where the local Augmented Lagrangian is given by
$$%\begin{equation*}
    %
    %\displaystyle
    %\begin{array}{ll}
        %
        L_i(x_i,k)
        \DefinedAs
        f_i \left( x_i \right)
        +
        \sum_{j \in \mathcal{N}_i} 
        y_{(i,j)} \big( x_i - z_{(i,j)} \big)
        +
        \sum_{j \in \mathcal{N}_i} 
        \frac{\delta}{2}
        \big\|
            x_i - z_{(i,j)}
        \big\|^2
        ,
        %
    %\end{array}
    %
$$%\end{equation*}
with $\delta$ a tuning parameter (see~\cite{ghadimi_et_al__2012__on_the_optimal_step_size_selection_for_the_admm} for a discussion on how to tune it) and the $y_{(i,j)}$'s Lagrange multipliers.

\begin{algorithm}
\caption{\ac{ADMM}~\cite[pp.~253-261]{bertsekas_tsitsiklis__1997__parallel_and_distributed_computation}}
    \begin{algorithmic}[1]
		\Statex \GrayText{(storage allocation and  constraints	on	parameters)}
        \State $x_i(k), z_{(i,j)}(k), y_{(i,j)}(k) \in \Reals^{n}$, $\delta \in \left( 0, 1 \right)$
        \Statex \GrayText{(initialization)}
		\State $x_i(k) = z_{(i,j)}(k) =  y_{(i,j)}(k)  =  0$
        \Statex \GrayText{(main algorithm)}
        \For{$k = 0, 1, \ldots$}
            \For{$i = 1, \ldots, N$}
                \State
                $
                    \displaystyle
                    x_i(k+1)
                    =
                    \arg \min_{x_i} L_i(x_i,k)
                $
                \For{$j \in \mathcal{N}_i$}
                    \State \vspace{0.2cm}
                    $
                        \displaystyle
                        z_{(i,j)}(k+1)
                        =
                        \frac{1}{2 \delta}
                        \Big(
                            y_{(i,j)}(k)
                            +
                            y_{(j,i)}(k)
                        \Big)
                        +
                        \frac{1}{2}
                        \Big(
                            x_{i}(k + 1)
                            +
                            x_{j}(k+1)
                        \Big)
                    $
                    \State \vspace{0.1cm}
                    $
                        y_{(i,j)}(k + 1)
                        =
                        y_{(i,j)}(k)
                        +
                        \delta \Big( x_i(k + 1) - z_{(i,j)}(k + 1) \Big)
                    $
                \EndFor
            \EndFor
        \EndFor

    \end{algorithmic}
\label{alg:ADMM}
\end{algorithm}

	The computational, communication and memory costs of these algorithms is reported in \Table~\ref{tab:computational_communication_and_memory_costs_of_competing_algorithms}. Notice that the computational and memory costs of \ac{ADMM} algorithms depends on how nodes minimize the local augmented Lagrangian $L_{i}(x_i,k)$. E.g., in our simulations the step has been performed through a dedicated Newton-Raphson procedure with associated $\BigOOf{n^{3}}$ computational costs and $\BigOOf{n^{2}}$ memory costs.

\begin{table}[!htbp]
	\centering
\begin{tabular}
{
	>{\arraybackslash \columncolor{black!05!white}}   			p{0.28\columnwidth}	 <{}
    >{\centering \arraybackslash}                               p{0.11\columnwidth} <{}
    >{\centering \arraybackslash \columncolor{black!05!white}}  p{0.11\columnwidth} <{}
    >{\centering \arraybackslash}                               p{0.30\columnwidth} <{}
}
\toprule % --------------------------------------------

Choice & 
DSM &
DCM &
ADMM \\

\midrule

Computational Cost  &	$\BigOOf{n}$	&	$\BigOOf{n}$	&   from $\BigOOf{n}$ to $\BigOOf{n^{3}}$	\\
Communication Cost  &	$\BigOOf{n}$	&	$\BigOOf{n}$	&   $\BigOOf{n}$	\\
Memory Cost			&	$\BigOOf{n}$	&	$\BigOOf{n}$ 	&	from $\BigOOf{n}$ to $\BigOOf{n^{2}}$	\\

\bottomrule % --------------------------------------------
\end{tabular}
\caption{Computational, communication and memory costs of DSM, DCM, and ADMM per single unit and single step.}
\label{tab:computational_communication_and_memory_costs_of_competing_algorithms}
\end{table}

\begin{figure}[!htbp]
	\centering
    \subfigure
    [Relative \ac{MSE} at a given time $k$ as a function of the algorithms parameters for problem~\eqref{equ:spam_nonspam_classification}. For the \ac{DCM}, $\nu = 1.7$.]
    {
		\includegraphics{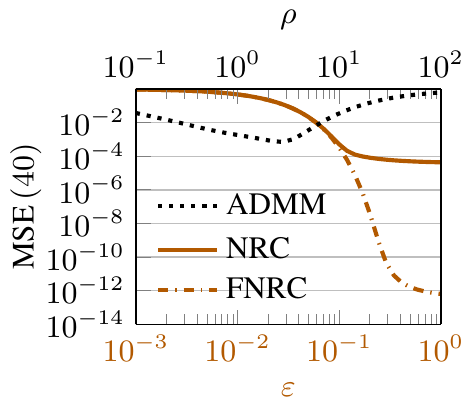}
		\includegraphics{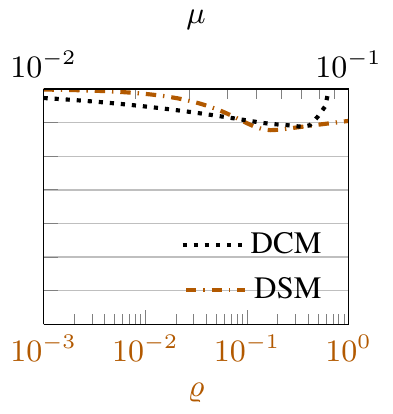}
        \label{fig:spam_nonspam__RMSE_at_40_vs_parameters__NRC_FNRC_ADMM_SG_CM}
    }
    \subfigure
    [Relative \ac{MSE} as a function of the time $k$ for the three fastest algorithms for problem~\eqref{equ:spam_nonspam_classification}. Their parameters are chosen as the best ones from \Figure~\ref{fig:spam_nonspam__RMSE_at_40_vs_parameters__NRC_FNRC_ADMM_SG_CM}.] 
    {
		\includegraphics{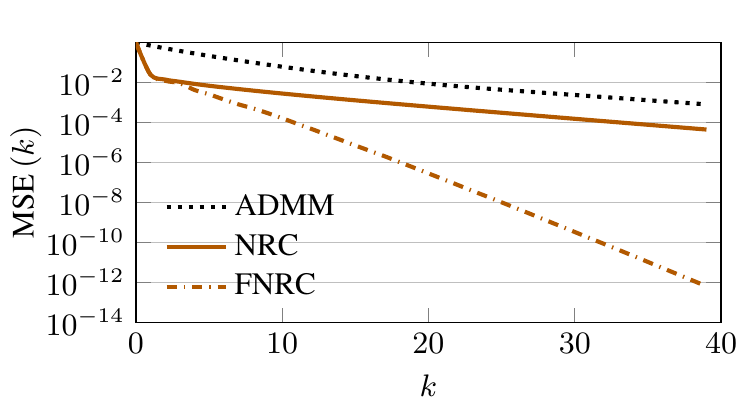}
        \label{fig:spam_nonspam__RMSE_evolution__NRC_FNRC_ADMM_SG_CM}
    }
    \\
    \subfigure
    [Relative \ac{MSE} at a given time $k$ as a function of the algorithms parameters for problem~\eqref{equ:housing_regression}. For the \ac{DCM}, $\nu = 1.7$.]
    {
       \includegraphics{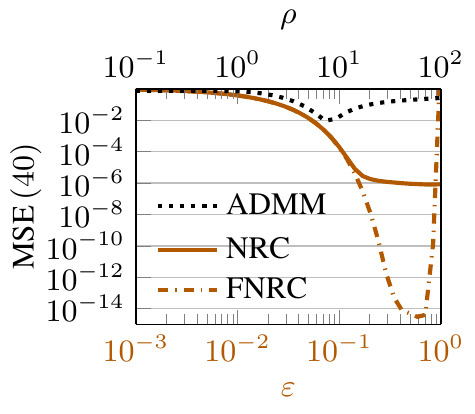}
       \includegraphics{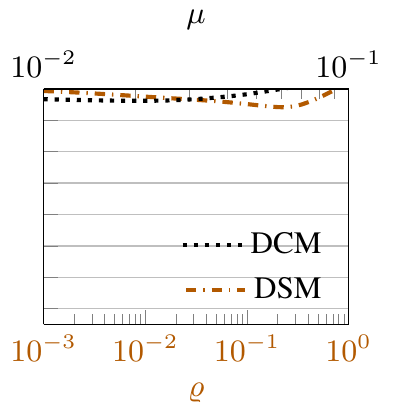}
        \label{fig:housing__RMSE_at_40_vs_parameters__NRC_FNRC_ADMM_SG_CM}
    }
    \subfigure
    [Relative \ac{MSE} as a function of the time $k$ for the three fastest algorithms for problem~\eqref{equ:housing_regression}. Their parameters are chosen as the best ones from \Figure~\ref{fig:housing__RMSE_at_40_vs_parameters__NRC_FNRC_ADMM_SG_CM}.] 
    {
       \includegraphics{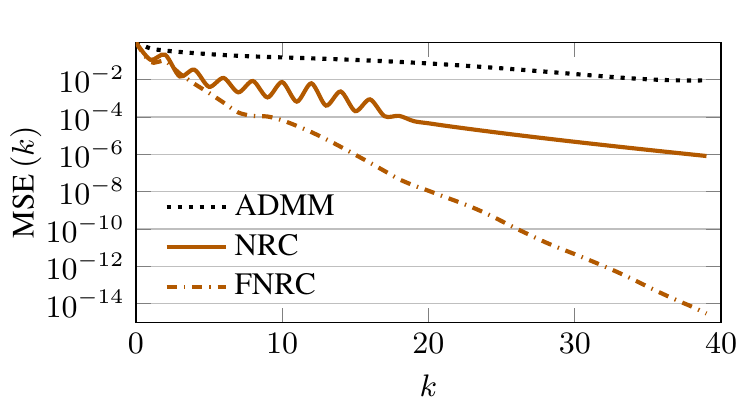}
        \label{fig:housing__RMSE_evolution__NRC_FNRC_ADMM_SG_CM}
    }
\caption{Convergence properties of the various algorithms for the problems described in \Section~\ref{ssec:optimization_problems}.}
\label{fig:comparison_of_NRC_FNRC_ADMM_SG_CM}
\end{figure}

\Figure~\ref{fig:comparison_of_NRC_FNRC_ADMM_SG_CM} then compares the previously cited algorithms as did in \Figure~\ref{fig:comparison_of_NRC_JC_and_GD}. The first panel thus reports the relative \ac{MSE} of the various algorithms at a given number of iterations ($k = 40$) as a function of the parameters. The second panel instead reports the temporal evolution of the relative \ac{MSE} for the case of optimal tuning. 

We notice that the \ac{DCM} and the \ac{DSM} are both much slower, in terms of communications iterations, than the \ac{NRC}, \ac{FNRC} and \ac{ADMM}. Moreover, both the \ac{NRC} and its accelerated version converge faster than the \ac{ADMM}, even if not tuned at their best. These numerical examples seem to indicate that the proposed \ac{NRC} might be a viable alternative to the \ac{ADMM}, although further comparisons are needed to strengthen this claim. Moreover, a substantial potential advantage of \ac{NRC} as compared to \ac{ADMM} is that the former can be readily adapted to asynchronous and time-varying graphs, as preliminary made in~\cite{zanella_et_al__2012__asynchronous_newton_raphson_consensus_for_distributed_convex_optimization}. Moreover, as in the case of the \ac{FNRC}, the strategy can implement any improved linear consensus algorithm.

% -----------------------------------------------------
\section{Conclusion}
\label{sec:conclusions}

We proposed a novel distributed optimization strategy suitable for convex, unconstrained, multidimensional, smooth and separable cost functions. The algorithm does not rely on Lagrangian formalisms and acts as a distributed \acl{NR} optimization strategy by repeating the following steps: agents first locally compute and update second order Taylor expansions around the current local guesses and then they suitably combine them by means of average consensus algorithms to obtain a sort of approximated Taylor expansion of the global cost. This allows each agent to infer a local Newton direction, used to locally update the guess of the global minimum.

Importantly, the average consensus protocols and the local updates steps have different time-scales, and the whole algorithm is proved to be convergent only if the step-size is sufficiently slow. Numerical simulations based on real-world databases show that, if suitably tuned, the proposed algorithm is faster then \acp{ADMM} in terms of number of communication iterations, although no theoretical proof is provided.

The set of open research paths is extremely vast. We envisage three main avenues. The first one is to study how the agents can dynamically and locally tune the speed of the local updates w.r.t.\ the consensus process, namely how to tune their local step-size $\varepsilon_i$. In fact large values of $\varepsilon$ gives faster convergence but might lead to instability. A second one is to let the communication protocol be asynchronous: in this regard we notice that some preliminary attempts can be found in~\cite{zanella_et_al__2012__asynchronous_newton_raphson_consensus_for_distributed_convex_optimization}. A final branch is about the analytical characterization of the rate of convergence of the proposed strategies, a theoretical comparison with \acp{ADMM}, and the extensions to non-smooth convex functions.

% ~~~~~~~~~~~~~~~~~~~~~~~~~~~~~~~~~~~~~~~~~~~~~~~~~~~~~~~~~~~~~~~~~ %
%                                                                   %
\addcontentsline    {toc}{section}{References}                      %
% \bibliographystyle  {../Bibliography/ieeetransactions_no_urls}      %
% \bibliography       {../Bibliography/NewtonRaphsonConsensus}        %

\begin{thebibliography}{10}
\providecommand{\url}[1]{#1}
\csname url@rmstyle\endcsname
\providecommand{\newblock}{\relax}
\providecommand{\bibinfo}[2]{#2}
\providecommand\BIBentrySTDinterwordspacing{\spaceskip=0pt\relax}
\providecommand\BIBentryALTinterwordstretchfactor{4}
\providecommand\BIBentryALTinterwordspacing{\spaceskip=\fontdimen2\font plus
\BIBentryALTinterwordstretchfactor\fontdimen3\font minus
  \fontdimen4\font\relax}
\providecommand\BIBforeignlanguage[2]{{%
\expandafter\ifx\csname l@#1\endcsname\relax
\typeout{** WARNING: IEEEtran.bst: No hyphenation pattern has been}%
\typeout{** loaded for the language `#1'. Using the pattern for}%
\typeout{** the default language instead.}%
\else
\language=\csname l@#1\endcsname
\fi
#2}}


\bibitem{zanella_et_al__2011__newton_raphson_consensus_for_distributed_convex_optimization}
F.~Zanella, D.~Varagnolo, A.~Cenedese, G.~Pillonetto, and L.~Schenato,
  ``{Newton-Raphson consensus for distributed convex optimization},'' in
  \emph{IEEE Conference on Decision and Control and European Control
  Conference}, Dec. 2011, pp. 5917--5922.

\bibitem{zanella_et_al__2012__multidimensional_newton_raphson_consensus_for_distributed_convex_optimization}
------, ``{Multidimensional Newton-Raphson consensus for distributed convex
  optimization},'' in \emph{American Control Conference}, 2012.

\bibitem{shor__1985__minimization_methods_for_non_differentiable_functions}
N.~Z. Shor, \emph{{Minimization Methods for Non-Differentiable
  Functions}}.\hskip 1em plus 0.5em minus 0.4em\relax Springer-Verlag, 1985.

\bibitem{bertsekas_et_al__2003__convex_analysis_and_optimization}
D.~P. Bertsekas, A.~Nedi\'{c}, and A.~E. Ozdaglar, \emph{{Convex Analysis and
  Optimization}}.\hskip 1em plus 0.5em minus 0.4em\relax Athena Scientific,
  2003.

\bibitem{boyd_vandenberghe__2004__convex_optimization}
S.~Boyd and L.~Vandenberghe, \emph{{Convex Optimization}}.\hskip 1em plus 0.5em
  minus 0.4em\relax Cambridge University Press, 2004.

\bibitem{tsitsiklis__1984__problems_in_decentralized_decision_making_and_computation}
\BIBentryALTinterwordspacing
J.~N. Tsitsiklis, ``{Problems in decentralized decision making and
  computation},'' Ph.D. dissertation, Massachusetts Institute of Technology,
  1984.
\BIBentrySTDinterwordspacing

\bibitem{bertsekas_tsitsiklis__1997__parallel_and_distributed_computation}
D.~P. Bertsekas and J.~N. Tsitsiklis, \emph{{Parallel and Distributed
  Computation: Numerical Methods}}.\hskip 1em plus 0.5em minus 0.4em\relax
  Athena Scientific, 1997.

\bibitem{bertsekas__1998__network_optimization__continuous_and_discrete_models}
\BIBentryALTinterwordspacing
D.~P. Bertsekas, \emph{{Network Optimization: Continuous and Discrete
  Models}}.\hskip 1em plus 0.5em minus 0.4em\relax Belmont, Massachusetts:
  Athena Scientific, 1998.
\BIBentrySTDinterwordspacing

\bibitem{burger_et_al__2012__a_distributed_simplex_algorithm_for_degenerate_linear_programs_and_multi_agent_assignments}
\BIBentryALTinterwordspacing
M.~B\"{u}rger, G.~Notarstefano, F.~Bullo, and F.~Allg\"{o}wer, ``{A distributed
  simplex algorithm for degenerate linear programs and multi-agent
  assignments},'' \emph{Automatica}, vol.~48, no.~9, pp. 2298 -- 2304, 2012.
\BIBentrySTDinterwordspacing

\bibitem{bertsekas__1982__constrained_optimization_and_lagrange_multiplier_methods}
D.~P. Bertsekas, \emph{{Constrained Optimization and Lagrange Multiplier
  Methods}}.\hskip 1em plus 0.5em minus 0.4em\relax Boston, MA: Academic Press,
  1982.

\bibitem{hestenes__1969__multiplier_and_gradient_methods}
\BIBentryALTinterwordspacing
M.~R. Hestenes, ``{Multiplier and gradient methods},'' \emph{Journal of
  Optimization Theory and Applications}, vol.~4, no.~5, pp. 303--320, 1969.
\BIBentrySTDinterwordspacing

\bibitem{boyd_et_al__2010__distributed_optimization_and_statistical_learning_via_the_admm}
S.~Boyd, N.~Parikh, E.~Chu, B.~Peleato, and J.~Eckstein, ``{Distributed
  Optimization and Statistical Learning via the Alternating Direction Method of
  Multipliers},'' Stanford Statistics Dept., Tech. Rep., 2010.

\bibitem{erseghe_et_al__2011__fast_consensus_by_the_admm}
\BIBentryALTinterwordspacing
T.~Erseghe, D.~Zennaro, E.~Dall'Anese, and L.~Vangelista, ``{Fast Consensus by
  the Alternating Direction Multipliers Method},'' \emph{IEEE Transactions on
  Signal Processing}, vol.~59, no.~11, pp. 5523--5537, Nov. 2011.
\BIBentrySTDinterwordspacing

\bibitem{he_yuan__2011__on_the_o_1_on_t_convergence_rate_of_alternating_direction_method}
\BIBentryALTinterwordspacing
B.~He and X.~Yuan, ``{On the O(1/t) convergence rate of alternating direction
  method},'' \emph{SIAM Journal on Numerical Analysis (to appear)}, 2011.
\BIBentrySTDinterwordspacing

\bibitem{deng_yin__2012__on_the_global_linear_convergence_of_the_generalized_ADMM}
W.~Deng and W.~Yin, ``On the global and linear convergence of the generalized
  alternating direction method of multipliers,'' DTIC Document, Tech. Rep.,
  2012.

\bibitem{wei_ozdaglar__2012__distributed_alternating_direction_method_of_multipliers}
E.~Wei and A.~Ozdaglar, ``{Distributed Alternating Direction Method of
  Multipliers},'' in \emph{IEEE Conference on Decision and Control}, 2012.

\bibitem{mota_et_al__2012__distributed_admm_for_model_predictive_control_and_congestion_control}
J.~a. Mota, J.~Xavier, P.~Aguiar, and M.~P\"{u}schel, ``{Distributed ADMM for
  Model Predictive Control and Congestion Control},'' in \emph{IEEE Conference
  on Decision and Control}, 2012.

\bibitem{jakovetic_et_al__2011__cooperative_convex_optimization_in_networked_systems__augmented_lagrangian_algorithms_with_directed_gossip_communication}
\BIBentryALTinterwordspacing
D.~Jakoveti\'{c}, J.~a. Xavier, and J.~M.~F. Moura, ``{Cooperative convex
  optimization in networked systems: Augmented lagrangian algorithms with
  directed gossip communication},'' \emph{IEEE Transactions on Signal
  Processing}, vol.~59, no.~8, pp. 3889 -- 3902, Aug. 2011.
\BIBentrySTDinterwordspacing

\bibitem{demyanov_vasilev__1985__nondifferentiable_optimization}
V.~F. Dem'yanov and L.~V. Vasil'ev, \emph{{Nondifferentiable
  Optimization}}.\hskip 1em plus 0.5em minus 0.4em\relax Springer - Verlag,
  1985.

\bibitem{johansson__2008__on_distributed_optimization_in_networked_systems}
\BIBentryALTinterwordspacing
B.~Johansson, ``{On Distributed Optimization in Networked Systems},'' Ph.D.
  dissertation, KTH Royal Institute of Technology, 2008.
\BIBentrySTDinterwordspacing

\bibitem{boyd_et_al__2006__randomized_gossip_algorithms}
\BIBentryALTinterwordspacing
S.~Boyd, A.~Ghosh, B.~Prabhakar, and D.~Shah, ``{Randomized Gossip
  Algorithms},'' \emph{IEEE Transactions on Information Theory / ACM
  Transactions on Networking}, vol.~52, no.~6, pp. 2508--2530, June 2006.
\BIBentrySTDinterwordspacing

\bibitem{ribeiro__2010__ergodic_stochastic_optimization_algorithms_for_wireless_communication_and_networking}
\BIBentryALTinterwordspacing
A.~Ribeiro, ``{Ergodic stochastic optimization algorithms for wireless
  communication and networking},'' \emph{IEEE Transactions on Signal
  Processing}, vol.~58, no.~12, pp. 6369 -- 6386, Dec. 2010.
\BIBentrySTDinterwordspacing

\bibitem{nedic_bertsekas__2001__incremental_subgradient_methods_for_nondifferentiable_optimization}
\BIBentryALTinterwordspacing
A.~Nedi\'{c} and D.~P. Bertsekas, ``{Incremental subgradient methods for
  nondifferentiable optimization},'' \emph{SIAM Journal on Optimization},
  vol.~12, no.~1, pp. 109--138, 2001.
\BIBentrySTDinterwordspacing

\bibitem{nedic_et_al__2001__distributed_asynchronous_incremental_subgradient_methods}
\BIBentryALTinterwordspacing
A.~Nedi\'{c}, D.~Bertsekas, and V.~Borkar, ``{Distributed asynchronous
  incremental subgradient methods},'' \emph{Studies in Computational
  Mathematics}, vol.~8, pp. 381--407, 2001.
\BIBentrySTDinterwordspacing

\bibitem{nedic_ozdaglar__2008__approximate_primal_solutions_and_rate_analysis_for_dual_subgradient_methods}
\BIBentryALTinterwordspacing
A.~Nedi\'{c} and A.~Ozdaglar, ``{Approximate primal solutions and rate analysis
  for dual subgradient methods},'' \emph{SIAM Journal on Optimization},
  vol.~19, no.~4, pp. 1757 -- 1780, 2008.
\BIBentrySTDinterwordspacing

\bibitem{kiewel__2004__convergence_of_approximate_and_incremental_subgradient_methods_for_convex_optimization}
\BIBentryALTinterwordspacing
K.~C. Kiwiel, ``{Convergence of approximate and incremental subgradient methods
  for convex optimization},'' \emph{SIAM Journal on Optimization}, vol.~14,
  no.~3, pp. 807--840, 2004.
\BIBentrySTDinterwordspacing

\bibitem{blatt_et_al__2007__a_convergent_incremental_gradient_method_with_a_constant_step_size}
\BIBentryALTinterwordspacing
D.~Blatt, A.~Hero, and H.~Gauchman, ``{A convergent incremental gradient method
  with a constant step size},'' \emph{SIAM Journal on Optimization}, vol.~18,
  no.~1, pp. 29--51, 2007.
\BIBentrySTDinterwordspacing

\bibitem{xiao_boyd__2006__optimal_scaling_of_a_gradient_method_for_distributed_resource_allocation}
\BIBentryALTinterwordspacing
L.~Xiao and S.~Boyd, ``{Optimal scaling of a gradient method for distributed
  resource allocation},'' \emph{Journal of optimization theory and
  applications}, vol. 129, no.~3, pp. 469 -- 488, 2006.
\BIBentrySTDinterwordspacing

\bibitem{ram_et_al__2009__incremental_stochastic_subgradient_algorithms_for_convex_optimization}
\BIBentryALTinterwordspacing
S.~S. Ram, A.~Nedi\'{c}, and V.~V. Veeravalli, ``{Incremental stochastic
  subgradient algorithms for convex optimzation},'' \emph{SIAM Journal on
  Optimization}, vol.~20, no.~2, pp. 691--717, 2009.
\BIBentrySTDinterwordspacing

\bibitem{nedic_ozdaglar__2009__distributed_subgradient_methods_for_multi_agent_optimization}
\BIBentryALTinterwordspacing
A.~Nedi\'{c} and A.~Ozdaglar, ``{Distributed Subgradient Methods for
  Multi-Agent Optimization},'' \emph{IEEE Transactions on Automatic Control},
  vol.~54, no.~1, pp. 48--61, 2009.
\BIBentrySTDinterwordspacing

\bibitem{johansson_et_al__2009__a_randomized_incremental_subgradient_method_for_distributed_optimization_in_networked_systems}
\BIBentryALTinterwordspacing
B.~Johansson, M.~Rabi, and M.~Johansson, ``{A randomized incremental
  subgradient method for distributed optimization in networked systems},''
  \emph{SIAM Journal on Optimization}, vol.~20, no.~3, pp. 1157--1170, 2009.
\BIBentrySTDinterwordspacing

\bibitem{lobel_et_al__2011__distributed_multi_agent_optimization_with_state_dependent_communication}
\BIBentryALTinterwordspacing
I.~Lobel, A.~Ozdaglar, and D.~Feijer, ``{Distributed multi-agent optimization
  with state-dependent communication},'' \emph{Mathematical Programming}, vol.
  129, no.~2, pp. 255 -- 284, 2011.
\BIBentrySTDinterwordspacing

\bibitem{nedic__2010__asynchronous_broadcast_based_convex_optimization_over_a_network}
\BIBentryALTinterwordspacing
A.~Nedi\'{c}, ``{Asynchronous Broadcast-Based Convex Optimization over a
  Network},'' \emph{IEEE Transactions on Automatic Control}, vol.~56, no.~6,
  pp. 1337 -- 1351, June 2010.
\BIBentrySTDinterwordspacing

\bibitem{ghadimi_et_al__2012__accelerated_gradient_methods_for_networked_optimization}
E.~Ghadimi, I.~Shames, and M.~Johansson, ``{Accelerated Gradient Methods for
  Networked Optimization},'' \emph{IEEE Transactions on Signal Processing
  (under review)}, 2012.

\bibitem{jadbabaie_et_al__2009__a_distributed_newton_method_for_network_optimization}
A.~Jadbabaie, A.~Ozdaglar, and M.~Zargham, ``{A Distributed Newton Method for
  Network Optimization},'' in \emph{IEEE Conference on Decision and
  Control}.\hskip 1em plus 0.5em minus 0.4em\relax IEEE, 2009, pp. 2736--2741.

\bibitem{zargham_et_al__2011__accelerated_dual_descent_for_network_optimization}
\BIBentryALTinterwordspacing
M.~Zargham, A.~Ribeiro, A.~Ozdaglar, and A.~Jadbabaie, ``{Accelerated Dual
  Descent for Network Optimization},'' in \emph{American Control Conference},
  2011.
\BIBentrySTDinterwordspacing

\bibitem{wei_et_al__2011__a_distributed_newton_method_for_network_utility_maximization}
\BIBentryALTinterwordspacing
E.~Wei, A.~Ozdaglar, and A.~Jadbabaie, ``{A Distributed Newton Method for
  Network Utility Maximization},'' in \emph{IEEE Conference on Decision and
  Control}, 2010, pp. 1816 -- 1821.
\BIBentrySTDinterwordspacing

\bibitem{dembo_et_al__1982__inexact_newton_methods}
\BIBentryALTinterwordspacing
R.~S. Dembo, S.~C. Eisenstat, and T.~Steihaug, ``{Inexact newton methods},''
  \emph{SIAM Journal on Numerical}, vol.~19, no.~2, pp. 400--408, 1982.
\BIBentrySTDinterwordspacing

\bibitem{nedic_et_al__2010__constrained_consensus_and_optimization_in_multi_agent_networks}
\BIBentryALTinterwordspacing
A.~Nedi\'{c}, A.~Ozdaglar, and P.~A. Parrilo, ``{Constrained Consensus and
  Optimization in Multi-Agent Networks},'' \emph{IEEE Transactions on Automatic
  Control}, vol.~55, no.~4, pp. 922--938, Apr. 2010.
\BIBentrySTDinterwordspacing

\bibitem{zhu_martinez__2012__on_distributed_convex_optimization_under_inequality_and_equality_constraints}
M.~Zhu and S.~Mart\'{\i}nez, ``{On Distributed Convex Optimization Under
  Inequality and Equality Constraints},'' \emph{IEEE Transactions on Automatic
  Control}, vol.~57, no.~1, pp. 151--164, 2012.

\bibitem{fischione__2011__f_lipschitz_optimization_with_wsn_applications}
\BIBentryALTinterwordspacing
C.~Fischione, ``{F-Lipschitz Optimization with Wireless Sensor Networks
  Applications},'' \emph{IEEE Transactions on Automatic Control}, vol.~56,
  no.~10, pp. 2319 -- 2331, 2011.
\BIBentrySTDinterwordspacing

\bibitem{wang_elia__2010__control_approach_to_distributed_optimization}
\BIBentryALTinterwordspacing
J.~Wang and N.~Elia, ``{Control approach to distributed optimization},'' in
  \emph{Forty-Eighth Annual Allerton Conference}, vol.~1, no.~1.\hskip 1em plus
  0.5em minus 0.4em\relax Allerton, Illinois, USA: IEEE, Sept. 2010, pp.
  557--561.
\BIBentrySTDinterwordspacing

\bibitem{freris_zouzias__2012__fast_distributed_smoothing_for_network_clock_synchronization}
\BIBentryALTinterwordspacing
N.~Freris and A.~Zouzias, ``{Fast Distributed Smoothing for Network Clock
  Synchronization},'' in \emph{IEEE Conference on Decision and Control}, 2012.
\BIBentrySTDinterwordspacing

\bibitem{necoara_nedelcu__2014__distributed_dual_gradient_methods_and_error_bound_conditions}
I.~Necoara and V.~Nedelcu, ``Distributed dual gradient methods and error bound
  conditions,'' \emph{arXiv preprint arXiv:1401.4398}, 2014.

\bibitem{garin_schenato__2011__a_survey_on_distributed_estimation_and_control_applications_using_linear_consensus_algorithms}
\BIBentryALTinterwordspacing
F.~Garin and L.~Schenato, \emph{{A survey on distributed estimation and control
  applications using linear consensus algorithms}}.\hskip 1em plus 0.5em minus
  0.4em\relax Springer, 2011, vol. 406, ch.~3, pp. 75--107.
\BIBentrySTDinterwordspacing

\bibitem{khalil__2001__nonlinear_systems}
H.~K. Khalil, \emph{{Nonlinear Systems}}, 3rd~ed.\hskip 1em plus 0.5em minus
  0.4em\relax Prentice Hall, 2001.

\bibitem{nedic_olshevsky__2013__distributed_optimization_over_time_varying_directed_graphs}
A.~Nedic and A.~Olshevsky, ``Distributed optimization over time-varying
  directed graphs,'' in \emph{Decision and Control (CDC), 2013 IEEE 52nd Annual
  Conference on}.\hskip 1em plus 0.5em minus 0.4em\relax IEEE, 2013, pp.
  6855--6860.

\bibitem{fagnani_zampieri__2008__randomized_consensus_algorithms_over_large_scale_networks}
\BIBentryALTinterwordspacing
F.~Fagnani and S.~Zampieri, ``{Randomized consensus algorithms over large scale
  networks},'' \emph{IEEE Journal on Selected Areas in Communications},
  vol.~26, no.~4, pp. 634--649, May 2008.
\BIBentrySTDinterwordspacing

\bibitem{dominguez_garcia_et_al__2011__distributed_algorithms_for_consensus_and_coordination_in_the_presence_of_packet-dropping_communication_links}
\BIBentryALTinterwordspacing
A.~D. Dom\'{\i}nguez-Garc\'{\i}a, C.~N. Hadjicostis, and N.~H. Vaidya,
  ``{Distributed Algorithms for Consensus and Coordination in the Presence of
  Packet-Dropping Communication Links Part I : Statistical Moments Analysis
  Approach},'' Coordinated Sciences Laboratory, University of Illinois at
  Urbana-Champaign, Tech. Rep., 2011.
\BIBentrySTDinterwordspacing

\bibitem{kokotovic_et_al__1999__singular_perturbation_methods_in_control__analysis_and_design}
P.~Kokotovi\'{c}, H.~K. Khalil, and J.~O'Reilly, \emph{{Singular Perturbation
  Methods in Control: Analysis and Design}}, ser. Classics in applied
  mathematics.\hskip 1em plus 0.5em minus 0.4em\relax SIAM, 1999, no.~25.

\bibitem{tanabe__1985__global_analysis_of_continuous_analogues_of_the_levenberg_marquardt_and_newton_raphson_methods}
\BIBentryALTinterwordspacing
K.~Tanabe, ``{Global analysis of continuous analogues of the
  Levenberg-Marquardt and Newton-Raphson methods for solving nonlinear
  equations},'' \emph{Annals of the Institute of Statistical Mathematics},
  vol.~37, no.~1, pp. 189--203, 1985.
\BIBentrySTDinterwordspacing

\bibitem{hauser_nedic__2005__the_continuous_newton_raphson_method_can_look_ahead}
\BIBentryALTinterwordspacing
R.~Hauser and J.~Nedi\'{c}, ``{The Continuous Newton-Raphson Method Can Look
  Ahead},'' \emph{SIAM Journal on Optimization}, vol.~15, pp. 915--925, 2005.
\BIBentrySTDinterwordspacing

\bibitem{sahai_et_al__2012__hearing_the_clusters_of_a_graph__a_distributed_algorithm}
\BIBentryALTinterwordspacing
T.~Sahai, A.~Speranzon, and A.~Banaszuk, ``{Hearing the clusters of a graph: A
  distributed algorithm},'' \emph{Automatica}, vol.~48, no.~1, pp. 15--24, Jan.
  2012.
\BIBentrySTDinterwordspacing

\bibitem{becker_le_cun__1988__improving_the_convergence_of_back_propagation_learning_with_second_order_models}
\BIBentryALTinterwordspacing
S.~Becker and Y.~{Le Cun}, ``{Improving the convergence of back-propagation
  learning with second order methods},'' University of Toronto, Tech. Rep.,
  Sept. 1988.
\BIBentrySTDinterwordspacing

\bibitem{athuraliya_low__2000__optimization_flow_control_with_newton_like_algorithm}
\BIBentryALTinterwordspacing
S.~Athuraliya and S.~H. Low, ``{Optimization flow control with Newton like
  algorithm},'' \emph{Telecommunication Systems}, vol.~15, pp. 345--358, 2000.
\BIBentrySTDinterwordspacing

\bibitem{golub_van_loan__1996__matrix_computations}
\BIBentryALTinterwordspacing
G.~H. Golub and C.~F. {Van Loan}, \emph{{Matrix computations}}, 3rd~ed.\hskip
  1em plus 0.5em minus 0.4em\relax John Hopkins University Press, 1996.
\BIBentrySTDinterwordspacing

\bibitem{hastie_et_al__2001__the_elements_of_statistical_learning}
T.~Hastie, R.~Tibshirani, and J.~Friedman, \emph{{The Elements of Statistical
  Learning: Data Mining, Inference, and Prediction}}, 2nd~ed.\hskip 1em plus
  0.5em minus 0.4em\relax New York: Springer, 2001.

\bibitem{xiao_et_al__2007__distributed_average_consensus_with_least_mean_square_deviation}
\BIBentryALTinterwordspacing
L.~Xiao, S.~Boyd, and S.-J. Kim, ``{Distributed average consensus with
  least-mean-square deviation},'' \emph{Journal of Parallel and Distributed
  Computing}, vol.~67, no.~1, pp. 33--46, Jan. 2007.
\BIBentrySTDinterwordspacing

\bibitem{muthukrishnan_et_al__1998__first_and_second_order_diffusive_methods_for_rapid_coarse_distributed_load_balancing}
\BIBentryALTinterwordspacing
S.~Muthukrishnan, B.~Ghosh, and M.~H. Schultz, ``{First and Second Order
  Diffusive Methods for Rapid, Coarse, Distributed Load Balancing},''
  \emph{Theory of Computing Systems}, vol.~31, no.~4, pp. 331 -- 354, 1998.
\BIBentrySTDinterwordspacing

\bibitem{schizas_et_al__2008__consensus_in_ad_hoc_wsns_with_noisy_links__part_i__distributed_estimation_of_deterministic_signals}
\BIBentryALTinterwordspacing
I.~D. Schizas, A.~Ribeiro, and G.~B. Giannakis, ``{Consensus in Ad Hoc WSNs
  With Noisy Links - Part I: Distributed Estimation of Deterministic
  Signals},'' \emph{IEEE Transactions on Signal Processing}, vol.~56, pp.
  350--364, Jan. 2008.
\BIBentrySTDinterwordspacing

\bibitem{ghadimi_et_al__2012__on_the_optimal_step_size_selection_for_the_admm}
E.~Ghadimi, A.~Teixeira, I.~Shames, and M.~Johansson, ``{On the Optimal
  Step-size Selection for the Alternating Direction Method of Multipliers},''
  in \emph{Necsys}, 2012.

\bibitem{zanella_et_al__2012__asynchronous_newton_raphson_consensus_for_distributed_convex_optimization}
F.~Zanella, D.~Varagnolo, A.~Cenedese, G.~Pillonetto, and L.~Schenato,
  ``{Asynchronous Newton-Raphson Consensus for Distributed Convex
  Optimization},'' in \emph{Necsys 2012}, 2012.

\bibitem{Kudryavtsev:01}
L.~Kudryavtsev, \emph{Encyclopedia of Mathematics}.\hskip 1em plus 0.5em minus
  0.4em\relax Springer, 2001, ch. Implicit Fucntion.


\end{thebibliography}
%                                                                   %
% ~~~~~~~~~~~~~~~~~~~~~~~~~~~~~~~~~~~~~~~~~~~~~~~~~~~~~~~~~~~~~~~~~ %

\appendix

\begin{proof}[of \Theorem~\ref{thm:continuous_to_discrete}]
	\emph{proof of}~\ref{item:continuous_to_discrete:continuous}: integrating \eqref{equ:inequalities_continuous_to_discrete:nabla2_V} twice implies
	$$
		\frac{1}{2} a_{1} \|x\|^2
		\leq
		V(x)
		\leq
		\frac{1}{2} a_{2} \|x\|^2
	$$ 
	that, jointly with~\eqref{equ:inequalities_continuous_to_discrete:partial_V}, immediately guarantee global exponential stability for~\eqref{equ:continuous_system}~\cite[Thm.~4.10]{khalil__2001__nonlinear_systems}.

	\emph{proof of}~\ref{item:continuous_to_discrete:discrete}: consider
	\begin{eqnarray}
		\Delta V \big( x(k) \big)
		\DefinedAs
		V\big(x(k+1)\big) - V\big(x(k)\big).
		\label{equ:definition_of_Delta}
	\end{eqnarray}
	To prove the claim we show that $\Delta V \left( x(k) \right) \leq - d\|x(k)\|^2$ for some positive scalar $d$. To this aim, expand $V\big(x(k+1)\big)$ with a second order Taylor expansion around $x(k)$ with remainder in Lagrange form, to obtain
	$$
		V \big( x + \varepsilon \phi (x)\big)
		=
		V(x)
		+
		\varepsilon
		\frac
		{\partial V}{\partial x} \phi(x)
		+
		\frac{1}{2}
		\varepsilon^2
		\phi^T(x)
		\nabla^2 V(x_\varepsilon)
		\phi(x)			
	$$
	with $x_\varepsilon = x + \varepsilon'\phi(x)$ for $\varepsilon'\in[0,\varepsilon]$. Using inequalities~\eqref{equ:inequalities_continuous_to_discrete} we then obtain
	$$
		\begin{array}{rcl}
			\Delta V(x(k)) & = &
			V(x(k+1))-V(x(k)) \\
			& \leq &
			-\varepsilon a_{3} \|x(k)\|^2 + \frac{1}{2}\varepsilon^2 a_{2} a_{4}^2 \|x(k)\|^2 \\
			& = &
			-\varepsilon (a_{3} - \varepsilon \frac{1}{2} a_{2} a_{4}^2)\|x(k)\|^2.
		\end{array}
	$$
	Thus, for all $\varepsilon < \overline{\varepsilon} = \frac{2 a_{3}}{a_{2} a_{4}^{2}}$ the origin is globally exponentially stable.
\end{proof}

\begin{proof}[of \Theorem~\ref{thm:continuous_NR}]
	\emph{proof of}~\ref{item:continuous_NR:lyapunov}: \Assumption~\ref{ass:smoothness_of_global_function} guarantees that $V_{\textrm{NR}}(0)=0$ and $V_{\textrm{NR}}(x) > 0$ for $x \neq 0$. Moreover, for $x \neq 0$,
	$$
		\begin{array}{ll}
			\displaystyle
			\frac{\partial V_{\textrm{NR}}}{\partial x} \phi_{\textrm{NR}}(x) &
			= - \left( \nabla \overline{f'}(x) \right)^T  \overline{\overline{h'}}(x)^{-1} \nabla \overline{f'}(x) \\
			&
			=- \left\| \overline{\overline{h'}}(x)^{-\frac{1}{2}} \, \nabla \overline{f'}(x) \right\|^2 < 0.
		\end{array}
	$$ 

	\emph{proof of}~\ref{item:continuous_NR:scalars}: \Assumption~\ref{ass:smoothness_of_global_function} guarantees that~\eqref{eqn:NR_2} is satisfied with $b_1=c$ and $b_2=m$. To prove~\eqref{eqn:NR_4} we start by considering that~\eqref{eqn:NR_2} guarantees $c \|x\|\leq \| \nabla \overline{f'}(x)\| \leq m \|x\|$. This in its turn implies
	$$
		\left\| \phi_{\textrm{NR}}(x) \right\|
		=
		\left\|\overline{\overline{h'}}^{-1}(x) \nabla \overline{f'}(x) \right\|
		\leq
		\frac{1}{c} \left\| \nabla \overline{f'}(x) \right\|
		\leq
		\frac{m}{c} \| x \|
		=
		b_4 \| x \|.
	$$

	To prove~\eqref{eqn:NR_3} eventually consider then that~\eqref{eqn:NR_4} implies
	$$
		\begin{array}{ll}
			\displaystyle
			\frac{\partial V_{\textrm{NR}}}{\partial x} \phi_{\textrm{NR}}(x)
			& =
			- \left( \nabla \overline{f'}(x) \right)^T  \overline{\overline{h'}}(x)^{-1} \nabla \overline{f'}(x) \\
			& \leq
			\displaystyle
			- \frac{c^2}{m} \|x\|^2 = - b_3 \|x\|^2.
		\end{array}
	$$ 
\end{proof}

\begin{proof}[of \Theorem~\ref{thm:continuous_multidimensional_NR}]
	In the interest of clarity we analyze the case where the local costs $f'_{i}$ are scalar, i.e., $n=1$. The multivariable case is indeed a straightforward extension with just a more involved notation. We also recall the following equivalences:
	$$
	\begin{array}{c}
		\bm{x} = \bm{x}^\parallel + \bm{x}^\perp ,
		\qquad
		\left( \bm{x}^\perp \right)^T \bm{x}^\parallel = 0 , \\
		\|\bm{x}\|^2 = \left\| \bm{x}^\parallel \right\|^2 + \left\| \bm{x}^\perp \right\|^2 = N |\overline{x}|^2 + \left\| \bm{x}^\perp \right\|^2 .
	\end{array}
	$$
	
	\emph{proof of}~\ref{item:continuous_multidimensional_NR:lyapunov}: $V_{\textrm{PNR}}(\bm{0}) = 0$ and $V_{\textrm{PNR}}(\bm{x}) > 0$ for $\bm{x} \neq \bm{0}$ follow immediately from the fact that $V_{\textrm{NR}}(0) = 0$ and $V_{\textrm{NR}}(\overline{x}) > 0$ for $\overline{x} \neq 0$. $\dot{V}_{\textrm{PNR}} < 0$ is instead proved by proving~\eqref{eqn:PNR_3}.

	\emph{proof of inequality~\eqref{eqn:PNR_2}:} given~\eqref{equ:definition_of_V_PNR},
	$$
		\frac{\partial^{2} V_{\textrm{PNR}}(\bm{x})}{\partial \bm{x}^{2}}
		=
		\frac
		{\partial^{2} \left( V_{\textrm{NR}}(\overline{x}) + \frac{1}{2} \eta \left\| \bm{x}^{\perp} \right\|^{2} \right)}
		{\partial \bm{x}^{2}}.
	$$
	Since
	$
		0
		\leq
		\left\| \bm{x}^{\perp} \right\|^{2}
		\leq
		\| \bm{x} \|^{2}
	$
	and
	$$
		\frac{\partial^2 V_{\textrm{NR}}(\overline{x})}{\partial \bm{x}^2}
		=
		\frac{1}{N^2}\OnesVector \OnesVector^T
		\nabla^{2} V_{\textrm{NR}}(\overline{x}),
	$$
	thanks to~\eqref{eqn:NR_2} it follows immediately that~\eqref{eqn:PNR_2} holds with  
	$$
		b_5 \DefinedAs \min \left\{ \frac{b_1}{N}, \eta \right\},
		\quad
		b_6 \DefinedAs \max \left\{ \frac{b_2}{N}, \eta \right\}.
	$$

	\emph{proof of inequality}~\eqref{eqn:PNR_4}: since the origin of $\overline{f'}$ is a minimum, it follows that $\nabla \overline{f'}(0) = 0$, and thus $\overline{g'}(\bm 0) = 0$ (cf.~\eqref{equ:definition_of_gi}). Thus also $\psi(\bm 0) = 0$, that in turn implies $\| \psi(\bm{x}) \|\leq a_{\psi} \| \bm{x}\|$ by \Assumption~\ref{ass:global_lipschitzianity}. Therefore
	$$
	\|\phi_{\textrm{PNR}}( \bm{x})\| \leq \| \bm{x}\|+ N\|\psi(\bm{x})\|\leq (1+Na_\psi)\|\bm{x}\|=b_8  \|\bm{x}\| .
	$$

	\emph{proof of inequality}~\eqref{eqn:PNR_3}: since
	$$
	\frac{\partial \overline{x}}{\partial \bm{x}}
		=
		\frac{1}{N} \OnesVector_{N}^T,
		\qquad
		\frac{\partial \bm{x}^\perp}{\partial \bm{x}}
		=
		I-\frac{1}{N} \OnesVector_{N} \OnesVector_{N}^T
		\IDefinedAs
		\Pi ,
	$$ 
	it follows that
	$$
	\begin{array}{lll}
		\displaystyle
		\frac
		{\partial V_{\textrm{PNR}}}
		{\partial \bm{x}} \phi_{\textrm{PNR}}(\bm{x})
		&=&
		\displaystyle
		\left(
			\frac
			{\partial V_{\textrm{PNR}}}
			{\partial \overline{x}}
			\frac
			{\partial \overline{x}}
			{\partial \bm{x}}
			+
			\frac
			{\partial V_{\textrm{PNR}}}
			{\partial \bm{x}^\perp }
			\frac
			{\partial \bm{x}^\perp}
			{\partial \bm{x}}
		\right)
		\phi_{\textrm{PNR}}(\bm{x})
		\vspace{0.1cm} \\
		&=&
		\displaystyle
		\left(
			\frac
			{\partial V_{\textrm{NR}}(\overline{x})}
			{\partial \overline{x}}
			\frac{1}{N} \OnesVector_{N}^T
			+
			\eta (\bm{x}^\perp )^T \Pi
		\right)
		\phi_{\textrm{PNR}}(\bm{x}) .
	\end{array}
	$$
	Considering then~\eqref{equ:function_for_coupled_NR}, the definition of $\overline{x}$ and $\bm{x}^{\perp}$, and the fact that $\Pi \OnesVector_{N} = 0$, it follows that
	$$
	\begin{array}{lll}
		\displaystyle
		\frac
		{\partial V_{\textrm{PNR}}}
		{\partial \bm{x}} \phi_{\textrm{PNR}}(\bm{x})
		&=&
		\displaystyle
		\frac
		{\partial V_{\textrm{NR}}(\overline{x})}
		{\partial \overline{x}}
		\big( -\overline{x} + \psi(\bm{x}) \big)
		+
		\eta (\bm{x}^\perp)^T(-\bm{x}^\perp)\\
	\end{array}
	$$
	Adding and subtracting 
	$
		\displaystyle
		\frac
		{\partial V_{\textrm{NR}}(\overline{x})}
		{\partial \overline{x}}
		\psi(\bm{x}^\parallel)
	$, 
	and recalling definition~\eqref{equ:function_for_NR} and equivalence~\eqref{equ:important_equivalences:psi}, since $\big( -\overline{x} + \psi(\bm{x}^\parallel) \big) = \phi_{\textrm{NR}}(\overline{x})$ it then follows that
	$$
	\begin{array}{lll}
		\displaystyle
		\frac
		{\partial V_{\textrm{PNR}}}
		{\partial \bm{x}} \phi_{\textrm{PNR}}(\bm{x})
		&=&
		\displaystyle
		\frac
		{\partial V_{\textrm{NR}}(\overline{x})}
		{\partial \overline{x}}
		\phi_{\textrm{NR}}(\overline{x})
		-
		\eta \|\bm{x}^\perp\|^2 +
		\vspace{0.1cm} \\
		&&
		\displaystyle
		+
		\frac
		{\partial V_{\textrm{NR}}(\overline{x})}
		{\partial \overline{x}}
		\left( \psi(\bm{x}) - \psi(\bm{x}^\parallel) \right)
		\vspace{0.1cm} \\
		& \leq &
		\displaystyle
		- b_3 \overline{x}^2
		- \eta
		  \| \bm{x}^\perp \|^2 + b_2 \big| \overline{x} \big| a_\psi \|\bm{x}-\bm{x}^\parallel \|
		\vspace{0.1cm} \\
		& = &
		\displaystyle
		- b_3 \overline{x}^2 - \eta \|\bm{x}^\perp\|^2+b_2a_\psi |\overline{x}|  \|\bm{x}^\perp\| \\
		& \leq &
		\displaystyle
		- \frac{b_3+\eta}{2}\big( |\overline{x}|^2+\|\bm{x}^\perp\|^2\big)\\ 
		&\leq&
		\displaystyle
		- \frac{b_3+\eta}{2}\big(N |\overline{x}|^2+\|\bm{x}^\perp\|^2\big) \\
		&=&
		\displaystyle
		- \frac{b_3+\eta}{2N} \|\bm{x}\|^2= - b_7 \|\bm{x}\|^2 
	\end{array}
	$$
	where for obtaining the various inequalities we used the various assumptions and where the second inequality is valid for $\displaystyle \eta>\frac{b_2^2a_\psi^2}{b_3}$.
\end{proof}

\begin{proof}[of \Lemma~\ref{lem:phi_x_is_lipschitz}]
	
	\emph{proof of}~\eqref{eqn:phi_p_x_1}: notice that $\phi_x( \bm{x}, \bm{\chi} )$ is globally defined since $[ \; \cdot \;]_c$ ensures that the matrix inverse exists. Also note that, since $\overline{h'}(\bm{x}) \geq c I>\frac{c}{2}I$ by \Assumption~\ref{ass:global_lipschitzianity}, then there exists $r>0$ such that, for $\|\bm{x}\| + \|\bm{\chi}\| \leq r$,
	$$
	\phi_x(\bm{x}, \bm{\chi})
		=
		- \bm{x}
		- \OnesVector_{N} \otimes x^{\ast}
		+ \frac{\bm{\chi}^{y} + \OnesVector_{N} \otimes \left( \overline{g'}(\bm{x}) + \overline{h'}(\bm{x}) x^{\ast} \right)}
		{ \bm{\chi}^{z} + \OnesVector_{N} \otimes \overline{h'}(\bm{x})} .
	$$
	The differentiability of the elements defining $\phi_{x}$, plus the fact that $\left[ \cdot \right]_{c}$ acts as the identity in the neighborhood under consideration implies that $\phi_{x}$ is locally differentiable in $\|\bm{x}\| + \|\bm{\chi}\| \leq r$. In addition to this local differentiability, also observe that $\phi_x(\bm 0, \bm 0)=0$, therefore there must exist $a_1>0$ s.t.
	\begin{equation}
		\|\phi_x(\bm{x},\bm{\chi})\| \leq a_1 (\| \bm{x} \| + \| \bm{\chi} \|),
		\qquad
		\forall (\|\bm{x}\| + \| \bm{\chi} \|) \leq r.
		\label{equ:first_inequality_for_phi_x}
	\end{equation}
	
	To extend the linear inequality~\eqref{equ:first_inequality_for_phi_x} for $\left( \bm{x}, \bm{\chi} \right)$ s.t.\ $(\|\bm{x}\| + \| \bm{\chi} \|) \geq r$ we then prove that $\phi_x( \bm{x}, \bm{\chi} )$ cannot grow more than linearly globally. In fact,
	\begin{equation}
		\begin{array}{l}
			\|\phi_x(\bm{x},\bm{\chi})\|
			\leq \\
			\quad
			\leq
			\|\bm{x}\| + N \|x^*\| +
			\frac{2}{c} \left\| \bm{\chi}^{y} + \OnesVector \otimes \left( \overline{g'} (\bm{x}) + \overline{h'} (\bm{x}) x^{\ast} \right) \right\| \\
			\quad
			\leq
			\|\bm{x}\| + N \|x^*\| +
			\frac{2}{c} \| \bm{\chi} \| + \frac{2N}{c} \left( \| \overline{g'} (\bm{x}) \| + \| x^* \| \|\overline{h'}(\bm{x}) \| \right) \\
			\quad
			\leq
			\|\bm{x}\| + N \|x^*\| +
			\frac{2}{c} \| \bm{\chi} \|+\frac{2N}{c} a_g \| \bm{x} \|+\\
			\qquad
			+ \frac{2N}{c} \| x^* \| \left( a_h \| \bm{x} \| + \| \overline{h'}(0) \| \right) \\
			\quad
			\leq
			a_2 + a_3 \left( \|\bm{x} \| + \| \bm{\chi} \| \right), \qquad \forall \bm{x},\bm{\chi}
		\end{array}
		\label{equ:second_inequality_for_phi_x}
	\end{equation}
	where we used \Assumption~\ref{ass:global_lipschitzianity} and where $a_2$, $a_3$ are suitable positive scalars. In particular inequality~\eqref{equ:second_inequality_for_phi_x} is valid for $(\|\bm{x}\| + \| \bm{\chi} \|) > r$. As depicted in \Figure~\ref{fig:two_cones}, inequalities~\eqref{equ:first_inequality_for_phi_x} and~\eqref{equ:second_inequality_for_phi_x} define two cones, one affine (shifted by $a_{2}$) and one proper.
	
	\begin{figure}[!htbp]
		\centering
		\includegraphics{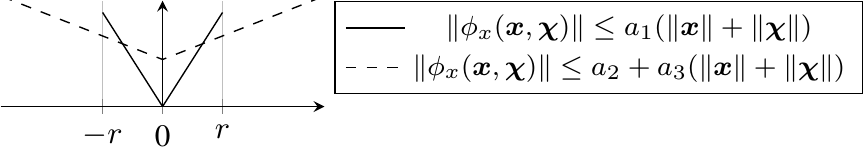}
		\caption{Inequality~\eqref{equ:first_inequality_for_phi_x} represents a proper cone defined in the neighborhood of radius $r$, while inequality~\eqref{equ:second_inequality_for_phi_x} represents an improper cone defined in the whole domain.}
		\label{fig:two_cones}
	\end{figure}

	Therefore, combining the geometry of the two cones leads to an inequality that is defined in the whole domain. In other words, it follows that
	$$
	\|\phi_x(\bm{x}, \bm{\chi})\|
		\leq
		a_x \big( \|\bm{x}\|+\|\bm{\chi}\| \big)
		\qquad
		\forall \bm{x}, \bm{\chi}
	$$
	where
	$$
		a_x \DefinedAs \max \left\{ a_1, \frac{a_2 + a_3 r}{r} \right\} .
	$$

	\emph{proof of}~\eqref{eqn:phi_p_x_2}: Let $\Delta(\bm{x}, \bm{\chi}) \DefinedAs \phi_x(\bm{x}, \bm{\chi}) - \phi_{\textrm{PNR}}(\bm{x})$, with $\phi_{\textrm{PNR}}$ as in~\eqref{equ:function_for_coupled_NR}. Then there exists a positive scalar $r>0$ such that, for all $\| \bm{\chi}\| + \|\bm{x}\| \leq r$,
	\begin{small}
	$$
	\begin{array}{l}
		\Delta(\bm{x},\bm{\chi}) = \\ 
		\quad =
		\displaystyle
		- \OnesVector_{N} \otimes x^{\ast}
		+ \frac
			{	  \bm{\chi}^{y}
				+ \OnesVector_{N} \otimes
					\left(
						\overline{g'} (\bm{x})
						+
						\overline{h'} (\bm{x}) x^{\ast}
					\right)
			}
			{
				\bm{\chi}^{z}
				+
				\OnesVector_{N} \otimes \overline{h'} (\bm{x})
			}
			- \OnesVector_N \otimes \psi(\bm{x}) \\
		\quad = 
		\displaystyle
		\frac
		{	  \bm{\chi}^{y}
			+ \OnesVector_{N} \otimes
				\left(
					\overline{g'} (\bm{x})
					+
					\overline{h'} (\bm{x}) x^{\ast}
				\right)
		}
		{
			\bm{\chi}^{z}
			+
			\OnesVector_{N} \otimes \overline{h'} (\bm{x})
		}
		-
		\frac
		{
			\OnesVector_{N} \otimes
			\left(
				\overline{g'} (\bm{x})
				+
				\overline{h'} (\bm{x}) x^{\ast}
			\right)
		}
		{
			\OnesVector_{N} \otimes \overline{h'} (\bm{x})
		} .
	\end{array}
	$$
	\end{small}
	Considerations similar to the ones that led us claim the differentiability of $\phi_{x}$ in the proof of \Lemma~\ref{lem:phi_x_is_lipschitz} imply that $\Delta(\bm{x},\bm{\chi})$ is continuously differentiable for $\| \bm{\chi} \| + \| \bm{x} \| \leq r$. Moreover, since $\Delta( \bm{x}, \bm 0 ) = 0$, then there exists a positive scalar $a_4 > 0$ s.t.
	\begin{equation}
		\| \Delta(\bm{x}, \bm{\chi}) \|
		\leq
		a_4 \|\bm{\chi}\|
		\qquad
		\|\bm{\chi}\| + \|\bm{x}\| \leq r .
	\label{eqn:phi_p_x_2:ancillary_1}
	\end{equation}
	By using~\eqref{equ:inequalities_for_global_lipschitzianity:h} and~\eqref{equ:inequalities_for_global_lipschitzianity:g} we can then show that $\Delta(\bm{x},\bm{\chi})$ cannot grow more than linearly in the variable $\bm{\chi}$, since
	\begin{equation}
	\begin{array}{l}
		\| \Delta(\bm{x},\bm{\chi}) \| = \\
		\quad =
		\displaystyle
		\left\|
			\frac
			{\bm{\chi}^{y} + \OnesVector_{N} \otimes \left( \overline{g'} (\bm{x}) + \overline{h'} (\bm{x}) x^{\ast} \right)}
			{ \left[ \bm{\chi}^{z} + \OnesVector_{N} \otimes \overline{h'} (\bm{x}) \right]_c }
			-
			\OnesVector_N \otimes
			\left(
				x^*
				+
				\frac
				{\overline{g'}(\bm{x})}
				{\overline{h'}(\bm{x})}
			\right)
		\right\| \\
		\quad \leq
		\displaystyle
		\frac{2}{c}
		\left(
			\| \bm{\chi} \|
			+	2 N \| \overline{g'} (\bm{x}) \|
				+	N \| x^* \| \| \overline{h} '(\bm{x}) \|
		\right)
		+
		N \| x^* \| \\
		\quad \leq
		a_5 + a_6 \| \bm{\chi} \|,
		\qquad \forall \bm{x}, \bm{\chi}
	\end{array}
	\label{eqn:phi_p_x_2:ancillary_2}
	\end{equation}
	for suitable positive scalars $a_5$ and $a_6$. Repeating the same geometrical arguments used above we then obtain
	$$
	\|\Delta(\bm{x}, \bm{\chi})\| \leq a_\Delta \|\bm{\chi}\|,
	\qquad \forall \bm{x}, \bm{\chi}
	$$
	with
	$$
		a_\Delta \DefinedAs \max\left\{ a_3, \frac{a_5 + a_6 r}{r} \right\} .
	$$

\end{proof}

\begin{proof}[of \Theorem~\ref{the:zeros_of_phi_xi_is_a_smooth_manifold}]
	For notational brevity we omit the dependence on $\xi$, i.e., let $\bm{x}^{eq} = \bm{x}^{eq}(\xi)$ and $x^{eq} = x^{eq}(\xi)$.

	We start by assuming that there exists a $\bm{x}^{eq}(\xi)$ satisfying~\eqref{equ:phi_x_has_equilibrium_x_eq} for $\|\xi\|\leq r$ and prove that $\bm{x}^{eq}(\xi)$ must satisfy $\bm{x}^{eq}(\xi) = \OnesVector_{N} \otimes x^{eq}(\xi)$ and~\eqref{equ:solution_of_x_of_xi_scalar}. Consider then $r$ sufficiently small. Then, since $\overline{h'}(\bm{x})>cI$ by \Assumption~\ref{ass:smoothness_of_global_function},
	$$
		\cmax{\bm{\xi}^z +\OnesVector_N \otimes \overline{h'}(\bm{x})}
		=
		\bm{\xi}^z +\OnesVector_N \otimes  \overline{h'}(\bm{x})
		=
		\OnesVector_N \otimes \left( \overline{h'}(\bm{x}) +\xi^{z} \right) .
	$$
	This implies that for $\| \xi \| \leq r$ we have
	$$
	\begin{array}{l}
		\phi_x (\bm{x}^{eq}, \bm{\xi}) =
		- \bm{x}^{eq} \\
		\quad - \OnesVector_N
		\otimes
		\Big(
		x^* -\!	\big( \xi^{z} + \overline{h'}(\bm{x}^{eq}) \big)^{-1}
				\big( \xi^{y} + \overline{g'}(\bm{x}^{eq}) + \overline{h'}(\bm{x}^{eq}) x^* \big)
		\Big)
	\end{array}
	$$
	Therefore $\phi_x (\bm{x}^{eq}, \bm{\xi}) = 0$ if and only if
	$$
		x_i^{eq}
		=
		- x^*
		+	\big( \xi^{z} + \overline{h'}(\bm{x}^{eq}) \big)^{-1}
		\big( \xi^{y} + \overline{g'}(\bm{x}^{eq}) + \overline{h'}(\bm{x}^{eq}) x^* \big) .
	$$
	Since the right-hand-side is independent of $i$, this implies both that the $\bm{x}^{eq}(\xi)$ satisfying~\eqref{equ:phi_x_has_equilibrium_x_eq} must satisfy $\bm{x}^{eq} = \OnesVector \otimes x^{eq}$, and that its expression is given by~\eqref{equ:solution_of_x_of_xi_scalar} (indeed~\eqref{equ:solution_of_x_of_xi_scalar} can be retrieved immediately from the equivalence above since
	$
		- x^*
		=
		\big( \xi^{z} + \overline{h'}(\bm{x}^{eq}) \big)^{-1}
		\big( - \xi^{z} x^* - \overline{h'}(\bm{x}^{eq}) x^* \big) 
	$).

	We now prove~\eqref{equ:phi_x_has_equilibrium_x_eq} by exploiting the Implicit Function Theorem \cite{Kudryavtsev:01}. If we indeed substitute the necessary condition $\bm{x}^{eq} = \OnesVector_{N} \otimes x^{eq}$ into the definition of $ \phi_{x}(\bm{x}^{eq},\xi)$, we obtain the parallelization of $N$ equivalent equations of the form
	$$
		x^{eq} + x^*
		=
		\left( \overline{\overline{h'}}(x^{eq}) + \xi^{z} \right)^{-1}
		\left(\overline{\overline{g}}(x^{eq}) + \xi^{y} + \overline{\overline{h'}}(x^{eq})x^* \right)
	$$
	where we used properties~\eqref{equ:important_equivalences:h} and~\eqref{equ:important_equivalences:g} that lead to
	$
		\overline{h'} \big( \OnesVector_{N} \otimes x \big)
		=
		\overline{\overline{h'}} \big( x \big)
	$
	and
	$
		\overline{g'} \big( \OnesVector_{N} \otimes x \big)
		=
		\overline{\overline{g'}} \big( x \big)
	$.

	Moreover, \Assumption~\ref{ass:global_lipschitzianity} ensures that $\overline{\overline{h'}}(x^*) \geq c I$. Thus, for the continuity assumptions in \Assumption~\ref{ass:smoothness_of_global_function}, there exists a sufficiently small $r > 0$ s.t.\ if $\|\xi^{z}\| \leq \|\xi\| \leq r$ then $\overline{\overline{h'}}(x^*) + \xi^{z}$ is still invertible. Therefore
	$$
		\overline{\overline{g'}} \big( x^{eq} \big)
		+
		\xi^{y} + \overline{\overline{h'}} \big( x^{eq} \big)x^*
		=
		\overline{\overline{h'}} \big( x^{eq} \big)
		(x^{eq}+x^*)
		+
		\xi^{z}
		(x^{eq}+x^*).
	$$

	Exploiting now the equivalence
	$
		\overline{\overline{g'}} \big( x^{eq} \big)
		=
		\overline{\overline{h'}} \big( x^{eq} \big)
		x^{eq}
		-
		\nabla \overline{f'} \big( x^{eq} \big)
	$,
	it follows that $x^{eq}$ must satisfy the following condition:
	$$
	\nabla \overline{f'}(x^{eq}) - \xi^{y} + \xi^{z} (x^{eq}+x^*) = 0 .
	$$
	Given \Assumption~\ref{ass:smoothness_of_global_function}, the left-hand side of the previous inequality is a continuously differentiable function, since
	$$
		\frac
		{\partial \big( \nabla \overline{f'}(x^{eq}) - \xi^{y} + \xi^{z} (x^{eq}+x^*) \big)}
		{\partial x^{eq}}
		=
		\nabla^2 \overline{f'}(x^{eq}) + \xi^{z}.
	$$
	Notice moreover that if $r$ is sufficiently small (i.e., $\|\xi^{z}\|$ is sufficiently small) then the differentiation is an invertible matrix, since once again $\nabla^2 \overline{f'}(x^*) \geq cI$ by assumption. Therefore, by the Implicit Function Theorem, $x^{eq}(\xi)$ exists, is unique and continuously differentiable.

\end{proof}

\begin{proof}[of \Theorem~\ref{thm:perturbed_NR_is_globally_exponentially_stable}]
	\emph{proof of}~\ref{item:perturbed_NR_is_globally_exponentially_stable:V_PNR_is_Lyapunov}: $V_{\textrm{PNR}}(\bm{0}) = 0$ and $V_{\textrm{PNR}}(\bm{x}) > 0$ for $\bm{x} \neq \bm{0}$ have been proved before. $\dot{V}_{\textrm{PNR}} < 0$ is instead proved by proving~\eqref{eqn:xi_3}.

	\emph{proof of}~\ref{item:inequalities_on_xi}: as for~\eqref{eqn:xi_3}, consider that, $\forall \bm{x}\in \Reals^{nN}$,
	$$
	\begin{array}{l}
		\displaystyle
		\frac
		{\partial V_{\textrm{PNR}}}
		{\partial\bm{x}}\phi'_x(\bm{x},\xi) = \\
		\quad
		=
		\displaystyle
		\frac
		{\partial V_{\textrm{PNR}}}
		{\partial\bm{x}}\phi'_x(\bm{x},0)
		+
		\frac
		{\partial V_{\textrm{PNR}}}
		{\partial\bm{x}}
		\big( \phi'_x(\bm{x},\xi)-\phi'_x(\bm{x},0) \big) \\
		\quad
		\leq
		\displaystyle
		\frac
		{\partial V_{\textrm{PNR}}}
		{\partial\bm{x}}\phi_{\textrm{PNR}}(\bm{x})
		+
		\left\|
			\frac
			{\partial V_{\textrm{PNR}}}
			{\partial\bm{x}}
		\right\|
		\left\|
		\phi'_x(\bm{x},\xi)
			-
			\phi'_x(\bm{x},0)
		\right\| \\
		\quad
		\leq
		- b_7\| \bm{x} \|^2
		+
		b_6 \|\bm{x} \| a_\xi \|\xi\| \|\bm{x}\| \\
		\quad
		\leq
		- ( b_7 - b_6 a_\xi r) \| \bm{x} \|^2
		\leq
		- b_7' \|\bm{x} \|^2.
	\end{array}
	$$
	Notice that this inequality is meaningful for $r < \frac{b_7}{b_6a_\xi}$.

	As for~\eqref{eqn:xi_4}, consider that, $\forall \bm{x}\in \Reals^{nN}$,
	$$
	\begin{array}{rcl}
		\| \phi'_x (\bm{x},\xi) \|
		& \leq &
		\big\|
		\phi'_x(\bm{x},0)
		\big\|
			+
		\big\|
		\phi'_x(\bm{x},\xi)
			-
			\phi'_x(\bm{x},0)
		\big\| \\
		& \leq &
		( b_8 + a_\xi r) \| \bm{x} \| \leq b'_8 \|\bm{x} \| .
	\end{array}
	$$

\end{proof}

\begin{proof}[of \Theorem~\ref{thm:quadratic_costs_that_satisfy}]
	 The miminizer of the global cost function is easily seen to be $x^* = \left( \sum_i A_i \right)^{-1} \left( \sum_i A_i d_i \right)$ from which it follows that $\overline{f'}(x) = \frac{1}{N} x^T A x$. Clearly $\overline{f}(x)$ satisfies \Assumption~\ref{ass:smoothness_of_global_function} since $ \nabla^2 \overline{f}(x)=\frac{1}{N}A>0$ is independent of $x$. Considering then $h'_i(x) = \nabla^2 f'_i(x) = A_i$ it follows after some suitable simplifications that:
	\begin{eqnarray*}
		\overline{h'}(\bm{x}) &=& \frac{1}{N} A \\
		g'_i(x) &=& A_i x- A_i(x+x^*-d_i)=A_i(d_i-x^*) \\
		g'(\bm{x})-g'(\bm{x}')&=&0 \\
		\overline{g'}(\bm{x})&=& \frac{1}{N}\left(\sum_i A_id_i -\sum_i A_i x^*\right)=0\\
		h'(\bm{x})-h'(\bm{x}')&=&0 \\
		\psi(\bm{x})&=& \overline{h}^{-1}(\bm{x})\overline{g}(\bm{x})=0\\
		x^{eq}(\xi)&=& \left( \frac{1}{N} A + \xi^{z} \right)^{-1} \left( \xi^{y}-\xi^zx^* \right) \\
		\phi'_x(\bm{x},\xi)&=& \phi'_x(\bm{x},0) = -\bm{x}
	\end{eqnarray*}
	where in the last equivalence we exploited definition~\eqref{equ:solution_of_x_of_xi_scalar}. Thus also the other assumptions are satisfied.
\end{proof}

\begin{proof}[of \Theorem~\ref{thm:distributed_NR_converges_to_the_global_optimum}]
	The proof considers the system as an autonomous singularly perturbed system, and proceeds as follows: \emph{a)} show that $x^{\ast}$ is an equilibrium; \emph{b)} perform a change of variables; \emph{c)} construct a Lyapunov function for the boundary layer system; \emph{d)} construct a Lyapunov function for the reduced system; \emph{e)} join the two Lyapunov functions into one, and show (by cascading the previously introduced Lemmas and Theorems) that the complete system~\eqref{equ:augmented_continuous_time_system} converges to $x^{\ast}$ while satisfying the hypotheses of \Theorem~\ref{thm:continuous_to_discrete}. By doing so it follows that~\eqref{equ:discrete_time_system}, i.e., \Algorithm~\ref{alg:distributed_NR}, is exponentially stable.

	For notational simplicity we let $\bm{x}^{\ast} \DefinedAs \OnesVector_{N} \otimes x^{\ast}$. We also use all the notation collected in \Section~\ref{sec:notation}.

	$\bullet$ \emph{Discrete to continuous dynamics)}  The dynamics of \Algorithm~\ref{alg:distributed_NR} can be written in state space as
	\begin{equation}
		\left\lbrace
			\begin{array}{rcl}
				\bm{v}(k)	&	=	&	g	 \big(	  \bm{x}(k-1)	 \big)	  \\
				\bm{w}(k)	&	=	&	h	 \big(	  \bm{x}(k-1)	 \big)	  \\
				\bm{y}(k) & = & P \Big[ \bm{y}(k-1) + g \big( \bm{x}(k-1) \big) - \bm{v}(k-1) \Big] \\
				\bm{z}(k) & = & P \Big[ \bm{z}(k-1) + h \big( \bm{x}(k-1) \big) - \bm{w}(k-1) \Big] \\
				\displaystyle
				\bm{x}(k) & = & \displaystyle (1 - \varepsilon) \bm{x}(k-1) + \varepsilon \frac{\bm{y}(k-1)}{\cmax{\bm{z}(k-1)}} \\
				%
	%             \bm{v}(0) & = & g\big( \bm{x}(-1)	 \big), \bm{w}(0) = h	 \big(	  \bm{x}(-1)	 \big), \bm{y}(0), \bm{z}(0),  \bm{x}(0)
				%
			\end{array}
		\right.
	\label{equ:discrete_time_system}
	\end{equation}
	with suitable initial conditions. \eqref{equ:discrete_time_system} can then be interpreted as the forward-Euler discretization of
	\begin{equation}
		\left\lbrace
			\begin{array}{rcl}
				\varepsilon \dot{\bm{v}}(t) & = & - \bm{v}(t) + g \big( \bm{x}(t) \big) \\%, \ \ \bm{v}(0)=0 \\
				\varepsilon \dot{\bm{w}}(t) & = & - \bm{w}(t) + h \big( \bm{x}(t) \big) \\%, \ \ \bm{w}(0)=0 \\
				\varepsilon \dot{\bm{y}}(t) & = & - K \bm{y}(t) + (I - K) \Big[ g \big( \bm{x}(t) \big) - \bm{v}(t) \Big] \\%, \ \ \bm{y}(0)=0 \\
				\varepsilon \dot{\bm{z}}(t) & = & - K \bm{z}(t) + (I - K) \Big[ h \big( \bm{x}(t) \big) - \bm{w}(t) \Big] \\%, \ \ \bm{z}(0)=0 \\
				\dot{\bm{x}}(t) & = & \displaystyle - \bm{x}(t) + \frac{\bm{y}(t)}{\cmax{\bm{z}(t)}} \\%, \ \ \bm{x}(0)=0 \\
				%
	%             \bm{v}(0) &=& g\big(	  \bm{x}(-1)	 \big), \bm{w}(0) = h	 \big(	  \bm{x}(-1)	 \big), \bm{y}(0), \bm{z}(0),  \bm{x}(0)
				%
			\end{array}
		\right.
	\label{equ:augmented_continuous_time_system}
	\end{equation}
	with null initial conditions, where $\varepsilon$ is the discretization time interval and $K \DefinedAs I - P$. Notice that, as for $P$, if $n$ is the dimension of the local costs then $P = P' \otimes I_{n}$ with $P'$ a doubly-stochastic average consensus matrix. Nonetheless for brevity we will omit the superscripts $'$.

	$\bullet$ \emph{b)} let
	$$
        \begin{array}{rcl}
			\bm{v}' &\DefinedAs& \bm{v} - g(\bm{x}) \\
			\bm{w}' &\DefinedAs& \bm{w} - h(\bm{x}) \\
			\bm{y}' &\DefinedAs& \bm{y} - \bm{v}' - \OnesVector_{N} \otimes \overline{g}(\bm{x}) \\
			\bm{z}' &\DefinedAs& \bm{z} - \bm{w}' - \OnesVector_{N} \otimes \overline{h}(\bm{x}) \\
			\bm{x}' &\DefinedAs& \bm{x} - \bm{x}^{\ast}
        \end{array}
	$$ 
	and
	\begin{equation*}
	\begin{array}{rcl}
		\vspace{0.1cm}
		\phi_g(\bm{x}')
		&\DefinedAs&
		\displaystyle
		\frac{\partial g}{\partial \bm{x}'}
		-
		\OnesVector_{N}
		\otimes
		\frac{\partial \overline{g}}{\partial \bm{x}'} \\
		\vspace{0.1cm}
		\phi_h(\bm{x}')
		&\DefinedAs&
		\displaystyle
		\frac{\partial h}{\partial \bm{x}'}
		-
		\OnesVector_{N}
		\otimes
		\frac{\partial \overline{h}}{\partial \bm{x}'} \\
		\phi_x ( \bm{x}', \chi )
		&\DefinedAs&
		\displaystyle
		- \bm{x}'(t)
		- \bm{x}^* + \\
		& &
		\displaystyle
		+ \frac
		{ \bm{y}'(t) + \bm{v}'(t) + \OnesVector_{N} \otimes \overline{g} \big( \bm{x}'(t) + \bm{x}^* \big) }
			{ \cmax{ \bm{z}'(t) + \bm{w}'(t) + \OnesVector_{N} \otimes \overline{h} \big( \bm{x}'(t) + \bm{x}^* \big)} } \\
	\end{array}
	\end{equation*}
	with $\bm{\chi} \DefinedAs \left( \bm{v}', \bm{w}', \bm{y}', \bm{z}' \right)$, so that~\eqref{equ:augmented_continuous_time_system} becomes
    \begin{equation}
        \left\lbrace
            \begin{array}{rcl}
				\vspace{0.1cm}
				\varepsilon \dot{\bm{v}}'(t)
				& = &
				\displaystyle
				- \bm{v}'(t) - \varepsilon \frac{\partial g}{\partial \bm{x}'} \dot{\bm{x}}'(t) \\
				\vspace{0.1cm}
				\varepsilon \dot{\bm{w}}'(t)
				& = &
				\displaystyle
				- \bm{w}'(t) - \varepsilon \frac{\partial h}{\partial \bm{x}'} \dot{\bm{x}}'(t)  \\
				\varepsilon \dot{\bm{y}}'(t)
				& = &
				\displaystyle
				- K \bm{y}'(t) + \varepsilon \phi_g(\bm{x}') \dot{\bm{x}}'(t)\\
				\varepsilon \dot{\bm{z}}'(t)
				& = &
				\displaystyle
				- K \bm{z}'(t) + \varepsilon \phi_h(\bm{x}') \dot{\bm{x}}'(t)\\
				\dot{\bm{x}}'(t)
                & = &
                \displaystyle
				\phi_x ( \bm{x}', \bm{\chi} ) \\
            \end{array}
        \right.
    \label{equ:augmented_continuous_time_system_changed_variables}
	\end{equation}
	with initial conditions
    \begin{equation*}
        \left\lbrace
            \begin{array}{rcl}
				\bm{v}'(0) &=& \bm{v}(0) - g \big(\bm{x}(0)\big) \\
				\bm{w}'(0) &=& \bm{w}(0) - h \big(\bm{x}(0)\big) \\
				\bm{y}'(0) &=& \bm{y}(0) - \bm{v}(0) + g^\perp \big(\bm{x}(0)\big) \\
				\bm{z}'(0) &=& \bm{z}(0) - \bm{w}(0) + h^\perp \big(\bm{x}(0)\big) \\
				\bm{x}'(0) &=& \bm{x}(0) - \bm{x}^\ast
            \end{array}
        \right.
	\end{equation*}
	where $g^{\perp}\left( \bm{x} \right) \DefinedAs g(\bm{x}) - \OnesVector_{N} \otimes \overline{g}(\bm{x})$ (equivalent definition for $h^{\perp}$). Notice that~\eqref{equ:augmented_continuous_time_system_changed_variables} has the origin as an equilibrium point. Moreover this dynamics exploits the function $\phi_x$ defined in~\eqref{equ:dynamics_of_phi_x}, with $\bm{\chi}^{y} = \bm{y}' + \bm{v}'$, and $\bm{\chi}^{z} = \bm{z}' + \bm{w}'$.

	The next step is to exploit the structure of $K$ (more precisely, the fact that it contains an average consensus matrix) to reduce the dynamics, i.e., to eliminate the dynamics of the average since the latter does not change in time. To this aim, we analyze the behavior of the average of the $y_{i}'$s, i.e., the behavior of $\left( \OnesVector_{N}^{T} \otimes I_{n} \right) \dot{\bm{y}}'$. To this point, consider the third equation in~\eqref{equ:augmented_continuous_time_system_changed_variables}. Recalling that $(A \otimes B) (C \otimes D) = AB \otimes CD$, and exploiting the fact that $\OnesVector_{N}^{T} P' = 0$, we notice that $\left( \OnesVector_{N}^{T} \otimes I_{n} \right) K = 0$. Moreover, from the definitions of $g$ and $\overline{g}$,
	$$
		\left( \OnesVector_{N}^T \otimes I_{n} \right)
		\frac{\partial g(\bm{x}')}{\partial \bm{x}'}
		=
		N \frac{\partial \overline{g}(\bm{x}')}{\partial \bm{x}'} .
	$$
	Since $N = \OnesVector_{N}^{T} \OnesVector_{N}$, it follows also that 
	$$
		\left( \OnesVector_{N}^T \otimes I_{n} \right)
		\phi_{g} \left( \bm{x}' \right)
		=
		0
	$$
	for all $t \geq 0$, i.e.,
	$
		\OnesVector^T \bm{y}'(t)
		=
		\OnesVector^T \bm{y}'(0)
		\equiv
		0
	$.
	Similarly it is possible to show that $\bm{z}'(t) \equiv 0$. This eventually implies that
	$$
		\bm{y}'^{\parallel}(t) = 0
		\qquad
		\bm{z}'^{\parallel}(t) = 0
		\qquad
		\forall t
	$$ 
	that means, recalling that $\bm{y}' = \bm{y}'^{\parallel} + \bm{y}'^{\perp}$ and $\bm{z}' = \bm{z}'^{\parallel} + \bm{z}'^{\perp}$, that~\eqref{equ:augmented_continuous_time_system_changed_variables} can be equivalently rewritten as
    \begin{equation}
        \left\lbrace
            \begin{array}{rcl}
				\vspace{0.1cm}
				\varepsilon \dot{\bm{v}}'(t)
				& = &
				\displaystyle
				- \bm{v}'(t) - \varepsilon \frac{\partial g}{\partial \bm{x}'} \phi_x ( \bm{x}', \bm{\chi}' ) \\
				\varepsilon \dot{\bm{w}}'(t)
				& = &
				\displaystyle
				- \bm{w}'(t) - \varepsilon \frac{\partial h}{\partial \bm{x}'} \phi_x ( \bm{x}', \bm{\chi}' ) \\
				\varepsilon \dot{\bm{y}}'^{\perp}(t)
				& = &
				\displaystyle
				- K \bm{y}'^{\perp}(t) + \varepsilon\phi_g(\bm{x}')\phi_x(\bm{x}',\bm{\chi}') \\
				\varepsilon \dot{\bm{z}}'^{\perp}(t)
				& = &
				\displaystyle
				- K \bm{z}'^{\perp}(t) + \varepsilon \phi_h(\bm{x}')\phi_x(\bm{x}',\bm{\chi}') \\
                \displaystyle
				\dot{\bm{x}}'(t)
                & = &
				\phi_x(\bm{x}',\bm{\chi}') \\
            \end{array}
        \right.
	\label{equ:augmented_continuous_time_system_final}
    \end{equation}
	where now $\bm{\chi}' \DefinedAs \left( \bm{v}, \bm{w}, \bm{y}'^{\perp}, \bm{z}'^{\perp} \right)$ and where the novel initial conditions for the changed variables are
    \begin{equation*}
        \left\lbrace
            \begin{array}{rcl}
				\bm{y}'^{\perp}(0) &=& \bm{y}^{\perp}(0) - \bm{v}^{\perp}(0) + g^\perp \big(\bm{x}(0)\big) \\
				\bm{z}'^{\perp}(0) &=& \bm{z}^{\perp}(0) - \bm{w}^{\perp}(0) + h^\perp \big(\bm{x}(0)\big) \\
            \end{array}
        \right.
	\end{equation*}

	$\bullet$ \emph{c)} the boundary layer of~\eqref{equ:augmented_continuous_time_system_final} is computed by setting $\bm{x}'(t) = \bm{x}'$. Since a constant $\bm{x}'$ implies $\dot{\bm{x}}' = \phi_x = 0$, this boundary layer reduces to a linear system globally exponentially converging to the origin. Notice that this implies that, in the original coordinates system,
	$$
		v = g(\bm{x}), \;\;
		w = h(\bm{x}), \;\;
		y = \OnesVector_{N} \otimes \overline{g}(\bm{x}), \;\;
		z = \OnesVector_{N} \otimes \overline{h}(\bm{x}).
	$$
	In the novel coordinates system we thus consider, as a Lyapunov function, $\frac{1}{2} \| \bm{\chi}' \|^{2}$.

	$\bullet$ \emph{d)} the reduced system of~\eqref{equ:augmented_continuous_time_system_final} is computed by plugging $\bm{\chi}' = \bm{0}$ into the equations (i.e., by setting $\bm{v}'(t) = \bm{0}$, $\bm{w}'(t) = \bm{0}$, $\bm{y}'^{\perp}(t) = \bm{0}$, $\bm{z}'^{\perp}(t) = \bm{0}$). Defining then
	$$
		f'_i(\bm{x}') \DefinedAs f_i(\bm{x}'+\bm{x}^*),
		\qquad
		h'_i(\bm{x}') \DefinedAs h_i(\bm{x}'+\bm{x}^*),
	$$
	we obtain
	\begin{small}
	$$
	\begin{array}{rcl}
		\dot{\bm{x}}'(t)
		&\!\!\!\!\!=\!\!\!\!\!&
		\displaystyle
		- \bm{x}'(t)
		- \bm{x}^*
		+
		\OnesVector_{N} \! \otimes \!
		\frac
		{\overline{g'} \big( \bm{x}'(t) \big)}
		{\overline{h'} \big( \bm{x}'(t) \big)} \\
		&\!\!\!\!\!=\!\!\!\!\!&
		\displaystyle
		- \bm{x}'(t)
		- \bm{x}^*
		+
		\OnesVector_{N} \! \otimes \!
		\frac
		{\overline{h'} \big( \bm{x}'(t) \big) \big( \bm{x}'(t) + \bm{x}^* \big) - \nabla f' \big( \bm{x}'(t) \big)}
		{\overline{h'} \big( \bm{x}'(t) \big)} \\
		&\!\!\!\!\!=\!\!\!\!\!&
		\displaystyle
		- \bm{x}'(t)
		+
		\OnesVector_{N} \otimes
		\frac
		{\overline{h'} \big( \bm{x}'(t) \big) \bm{x}'(t) - \nabla f' \big( \bm{x}'(t) \big)}
		{\overline{h'} \big( \bm{x}'(t) \big)} \\
		&\!\!\!\!\!=\!\!\!\!\!&
		\displaystyle
		- \bm{x}'(t)
		+
		\OnesVector_{N} \otimes
		\psi \big( \bm{x}'(t) \big) \\
		&\!\!\!\!\!=\!\!\!\!\!&
		\phi_{\textrm{PNR}} (\bm{x}')
	\end{array}
	$$
	\end{small}
	$\!\!\!\!\!\!\!$ where $\psi$ and $\phi_{\textrm{PNR}}$ are the functions defined  in~\eqref{equ:definition_of_psi} and~\eqref{equ:function_for_coupled_NR}, respectively. Thus the reduced system, thanks to \Theorem~\ref{thm:continuous_multidimensional_NR}, admits $\bm{x}^{\ast}$ as a global exponentially stable equilibrium, and admits $V_{\textrm{PNR}}$ in~\eqref{equ:definition_of_V_PNR} as a Lyapunov function.

	$\bullet$ \emph{e)} we now notice that the interconnection of the boundary layer and reduced systems maintains the global stability, since their Lyapunov functions are quadratic type. Thus (see~\cite[pp.~453]{khalil__2001__nonlinear_systems}) the global system is asymptotically globally stable. To check that forward-Euler discretizations of the system preserve these stability properties we then consider as a global Lyapunov function the function
	$$
		V \big( \bm{x}', \bm{\chi}' \big)
		=
		(1 - d) V_{\textrm{PNR}}(\bm{x}')
		+
		d \frac{1}{2} \| \bm{\chi}'\|^2 ,
	$$
	that is clearly positive definite for every $d \in (0, 1)$, and prove that inequalities~\eqref{equ:inequalities_continuous_to_discrete} of \Theorem~\ref{thm:continuous_to_discrete} are satisfied.

	\emph{proof that~\eqref{equ:inequalities_continuous_to_discrete:nabla2_V} holds}: from~\eqref{eqn:PNR_2} and the structure of $V$ it follows immediately that
	$$
		\left( \left( 1 - d \right) b_{5} + d \right) I
		\leq
		\nabla^{2} V \big( \bm{x}', \bm{\chi}' \big)
		\leq
		\left( \left( 1 - d \right) b_{6} + d \right) I .
	$$

	\emph{proof that~\eqref{equ:inequalities_continuous_to_discrete:phi_Lipschitz} holds}: applying~\eqref{equ:implied_inequalities_for_global_lipschitzianity} and~\eqref{eqn:phi_p_x_1} to~\eqref{equ:augmented_continuous_time_system_final} it follows that~\eqref{equ:inequalities_continuous_to_discrete:phi_Lipschitz} holds with
	$$
		a_{4}
		=
		a_{V}
		\DefinedAs
		\max \left\{ 1 + 2 \varepsilon a_{g} a_{x}, \; 1 + 2 \varepsilon a_{h} a_{x}, \; a_{x} \right\}.
	$$

	\emph{proof that~\eqref{equ:inequalities_continuous_to_discrete:partial_V} holds}: the part relative to the slow dynamics is already characterized by~\eqref{eqn:xi_3}. For the part relative to the fast dynamics, since $\frac{\partial \frac{1}{2} \| \bm{\chi} \|^{2}}{\partial \bm{\chi}} = \bm{\chi}^{T}$ to check that~\eqref{equ:inequalities_continuous_to_discrete:partial_V} corresponds to check the negativity of the terms
	$$
		\begin{array}{l}
		- \bm{v}'^{T} \bm{v}' - \varepsilon \bm{v}'^{T} \frac{\partial g}{\partial \bm{x}'} \phi_x ( \bm{x}', \bm{\chi}' ) \\
		- \bm{w}'^{T} \bm{w}' - \varepsilon \bm{w}'^{T} \frac{\partial h}{\partial \bm{x}'} \phi_x ( \bm{x}', \bm{\chi}' ) \\
		- \left( \bm{y}'^{\perp} \right)^{T} K \bm{y}'{\perp} + \varepsilon \bm{y}'^{\perp T} \phi_{g} (\bm{x}') \phi_x ( \bm{x}', \bm{\chi}' ) \\
		- \left( \bm{z}'^{\perp} \right)^{T} K \bm{z}'{\perp} + \varepsilon \bm{z}'^{\perp T} \phi_{h} (\bm{x}') \phi_x ( \bm{x}', \bm{\chi}' )
		\end{array}
	$$
	These terms can then be majorized using~\eqref{equ:implied_inequalities_for_global_lipschitzianity} and~\eqref{eqn:phi_p_x_1}. E.g., the third term can be majorized with
	$$
	- \sigma(P) \| \bm{y}'^{\perp} \|^{2} + 2 \varepsilon a_{g} a_{x} \| \bm{y}'^{\perp} \| \left( \| \bm{x} \| + \| \bm{\chi} \| \right)
	$$
	where $\sigma(P)$ is the spectral gap of $P$. Applying similar concepts also to the other terms it follows that~\eqref{equ:inequalities_continuous_to_discrete:partial_V} holds with
	$$
		a_{3}
		= 
		\min \left\{	\sigma(P) - 2 \varepsilon a_{g} a_{x}, \;
						\sigma(P) - 2 \varepsilon a_{h} a_{x}	\right\} .
	$$
\end{proof}

\begin{proof}[of \Theorem~\ref{thm:distributed_NR_converges_to_neighborhood_of_the_global_optimum}]
	The proof is identical to the one of \Theorem~\ref{thm:distributed_NR_converges_to_the_global_optimum} with the exception that the substitution is now $\bm{x}'' = \bm{x} - \bm{x}^{\ast} - \OnesVector_{N} \otimes \Psi(\xi^{y}, \xi^{z})$. Indeed one can prove the stability of the novel system using the same Lyapunov function of \Theorem~\ref{thm:distributed_NR_converges_to_the_global_optimum}. Notice that we are ensured that there exists a sufficiently small neighborhood of the origin for which the function $\Psi$ exists due to the smoothness conditions in \Assumption~\ref{ass:smoothness_of_global_function}.
\end{proof}

\begin{proof}[of \Theorem~\ref{thm:local_stability}]
	The proof is the local version of the one in \Theorem~\ref{thm:distributed_NR_converges_to_neighborhood_of_the_global_optimum}. Indeed the local versions of \Assumptions~\ref{ass:smoothness_of_global_function},~\ref{ass:global_lipschitzianity} and~\ref{ass:phi_p_xi_is_locally_uniformly_lipschitz} always hold, i.e., they hold when considering $\bm{x}$ s.t.\ $\| \bm{x} \| \leq r'$, and one can thus repeat that reasonings using local perspectives.
\end{proof}

\begin{proof}[of \Theorem~\ref{thm:convergence_for_quadratic_functions}]
	Consider for simplicity the scalar case. Let $y^* \DefinedAs \frac{1}{N}\sum_i A_i d_i$ and $z^* \DefinedAs \frac{1}{N} \sum_i A_i$, so that $\displaystyle x^* = \frac{y^*}{z^*}$. Since $\bm{y}(k+1) = P \bm{y}(k)$ and $\bm{z}(k+1) = P \bm{z}(k)$, given the assumptions on $P$, there exist positive $\alpha_{y}, \alpha_{z}$ independent of $\bm{x}(0)$ s.t.\ $|y_i(k) - y^*| \leq \alpha_{y} \left( \rho(P) \right)^{k}$ and $|z_i(k)-z^*|\leq \alpha_{z} \left( \rho(P) \right)^{k}$. The claim thus follows considering that $x_i(k)=\frac{y_i(k)}{\cmax{z_i(k)}}$ and that, since the elements of $P$ are non negative, all the $z_{i}(k)$ are non smaller than $c$ for all $k \geq 0$ (i.e., the operator $\cmax{\cdot}$ is always performing as the identity operator).
\end{proof}

\end{document}